# The Empirical Metamathematics of Euclid and Beyond

Stephen Wolfram*

*As an example of empirical metamathematics, we present a detailed study of the dependency structure of the 465 theorems in Euclid's Elements, finding empirical signatures of concepts such as the power of a theorem. We apply similar methods to a more exhaustive study of possible theorems in logic, as well as to the analysis of dependency structures in projects to formalize modern pure mathematics. We discuss the process of identifying both intrinsic features of metamathematical space, and features of its exploration through the historical progress of mathematics.*

## Towards a Science of Metamathematics

One of the many surprising things about our Wolfram Physics Project is that it seems to have implications even beyond physics. In our effort to develop a fundamental theory of physics it seems as if the tower of ideas and formalism that we've ended up inventing are actually quite general, and potentially applicable to all sorts of areas.

One area about which I've been particularly excited of late is metamathematics—where it's looking as if it may be possible to use our formalism to make what might be thought of as a "bulk theory of metamathematics".

Mathematics itself is about what we establish about mathematical systems. Metamathematics is about the infrastructure of how we get there—the structure of proofs, the network of theorems, and so on. And what I'm hoping is that we're going to be able to make an overall theory of how that has to work: a formal theory of the large-scale structure of metamathematics—that, among other things, can make statements about the general properties of "metamathematical space".

Like with physical space, however, there's not just pure underlying "geometry" to study. There's also actual "geography": in our human efforts to do mathematics over the last few millennia, where in metamathematical space have we gone, and "colonized"? There've been a few million mathematical theorems explicitly published in the history of human mathematics. What does the "empirical metamathematics" of them reveal? Some of it presumably reflects historical accidents, but some may instead reflect general features of metamathematics and metamathematical space.

A version of this document with immediately executable code is available at
writings.stephenwolfram.com/2020/09/the-empirical-metamathematics-of-euclid-and-beyond
Originally published September 28, 2020

*Email: s.wolfram@wolfram.com

I've wondered about empirical metamathematics for a long time, and tucked away on page 1176 at the end of the Notes for the section about "Implications for Mathematics and Its Foundations" in *A New Kind of Science* is something I wrote more than 20 years ago about it:

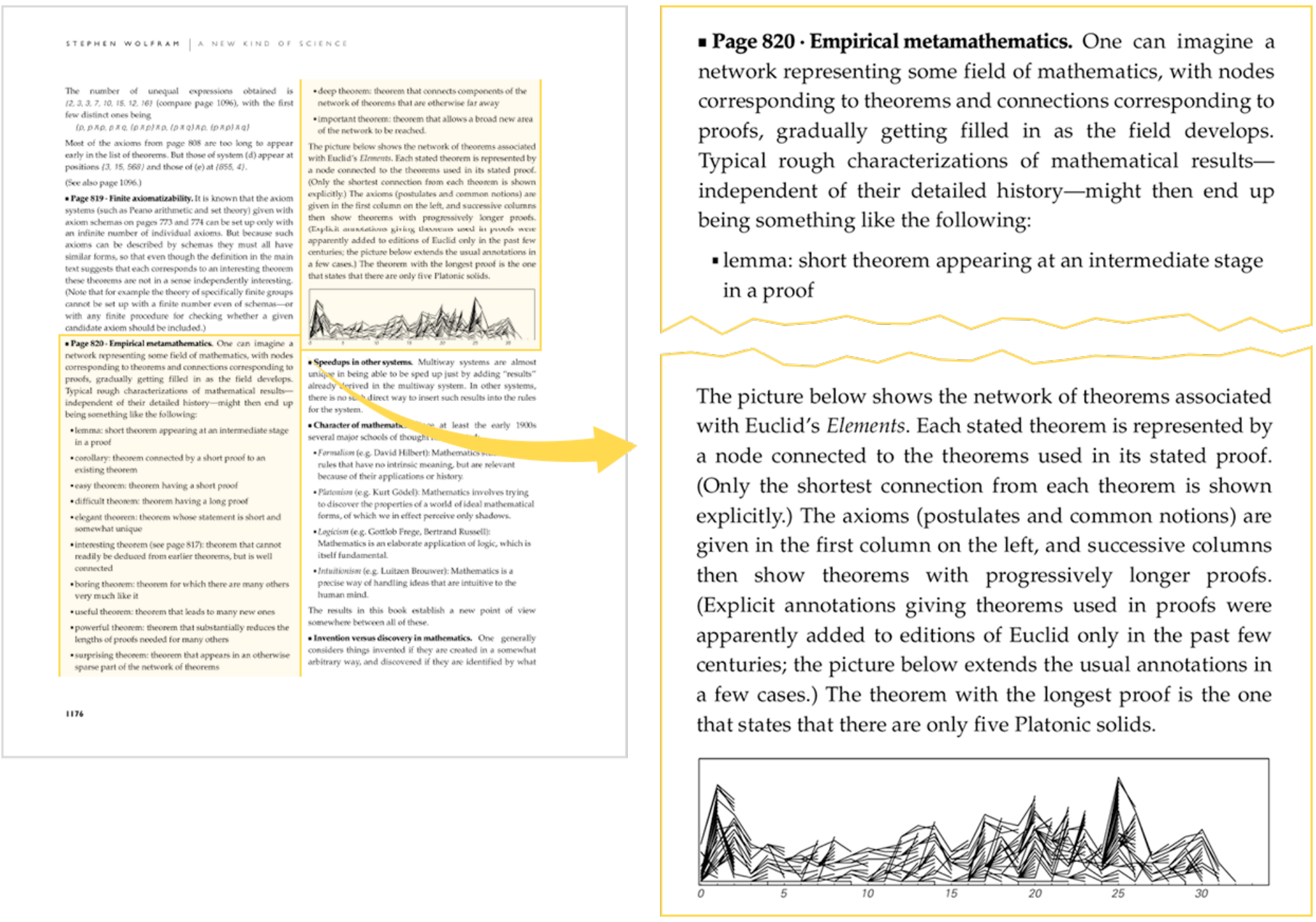

This note is mostly about what a descriptive theory of empirical metamathematics might be like—for example characterizing what one might mean by a powerful theorem, a deep theorem, a surprising theorem and so on. But at the end of the note is a graph: an actual piece of quantitative empirical metamathematics, based on the best-known structured piece of mathematics in history—Euclid's *Elements*.

The graph shows relationships between theorems in the *Elements*: a kind of causal graph of how different theorems make use of each other. As presented in *A New Kind of Science*, it's a small "footnote item" that doesn't look like much. But for more than 20 years, I've kept wondering what more there might be to learn from it. And now that I'm trying to make a general theory of metamathematics, it seemed like it was a good time to try to find out...

## The Most Famous Math Book in History

Euclid's *Elements* is an impressive achievement. Written in Greek around 300 BC (though presumably including many earlier results), the *Elements* in effect defined the way formal mathematics is done for more than two thousand years. The basic idea is to start from certain axioms that are assumed to be true, then—without any further "input from outside"—use purely deductive methods to establish a collection of theorems.

Euclid effectively had 10 axioms (5 “postulates” and 5 “common notions”), like “one can draw a straight line from any point to any other point”, or “things which equal the same thing are also equal to one another”. (One of his axioms was his fifth postulate—that parallel lines never meet—which might seem obvious, but which actually turns out not to be true for physical curved space in our universe.)

On the basis of his axioms, Euclid then gave 465 theorems. Many were about 2D and 3D geometry; some were about arithmetic and numbers. Among them were many famous results, like the Pythagorean theorem, the triangle inequality, the fact that there are five Platonic solids, the irrationality of $\sqrt{2}$ and the fact that there are an infinite number of primes. But certainly not all of them are famous—and some seem to us now pretty obscure. And in what has remained a (sometimes frustrating) tradition of pure mathematics for more than two thousand years, Euclid never gives any narrative about why he’s choosing the theorems he does, out of all the infinitely many possibilities.

We don’t have any original Euclids, but versions from a few centuries later exist. They’re written in Greek, with each theorem explained in words, usually by referring to a diagram. Mathematical notation didn’t really start getting invented until the 1400s or so (i.e. a millennium and a half later)—and even the notation for numbers in Euclid’s time was pretty unwieldy. But Euclid had basically modern-looking diagrams, and he even labeled points and angles with (Greek) letters—despite the fact that the idea of variables standing for numbers wouldn’t be invented until the end of the 1500s.

There’s a stylized—almost “legalistic”—way that Euclid states his theorems. And so far as we can tell, in the original version, all that was done was to state theorems; there was no explanation for why a theorem might be true—no proof offered. But it didn’t take long before people started filling in proofs, and there was soon a standard set of proofs, in which each particular theorem was built up from others—and ultimately from the axioms.

There’ve been more than a thousand editions of Euclid printed (probably more than any other book except the Bible), and reading Euclid was until quite recently part of any serious education. (At Eton—where I went to high school—it was only in the 1960s that learning “mathematics” began to mean much other than reading Euclid, in the original Greek of course.) Here’s an edition of Euclid from the 1800s that I happen to own, with the proof of every theorem giving little references to other theorems that are used:

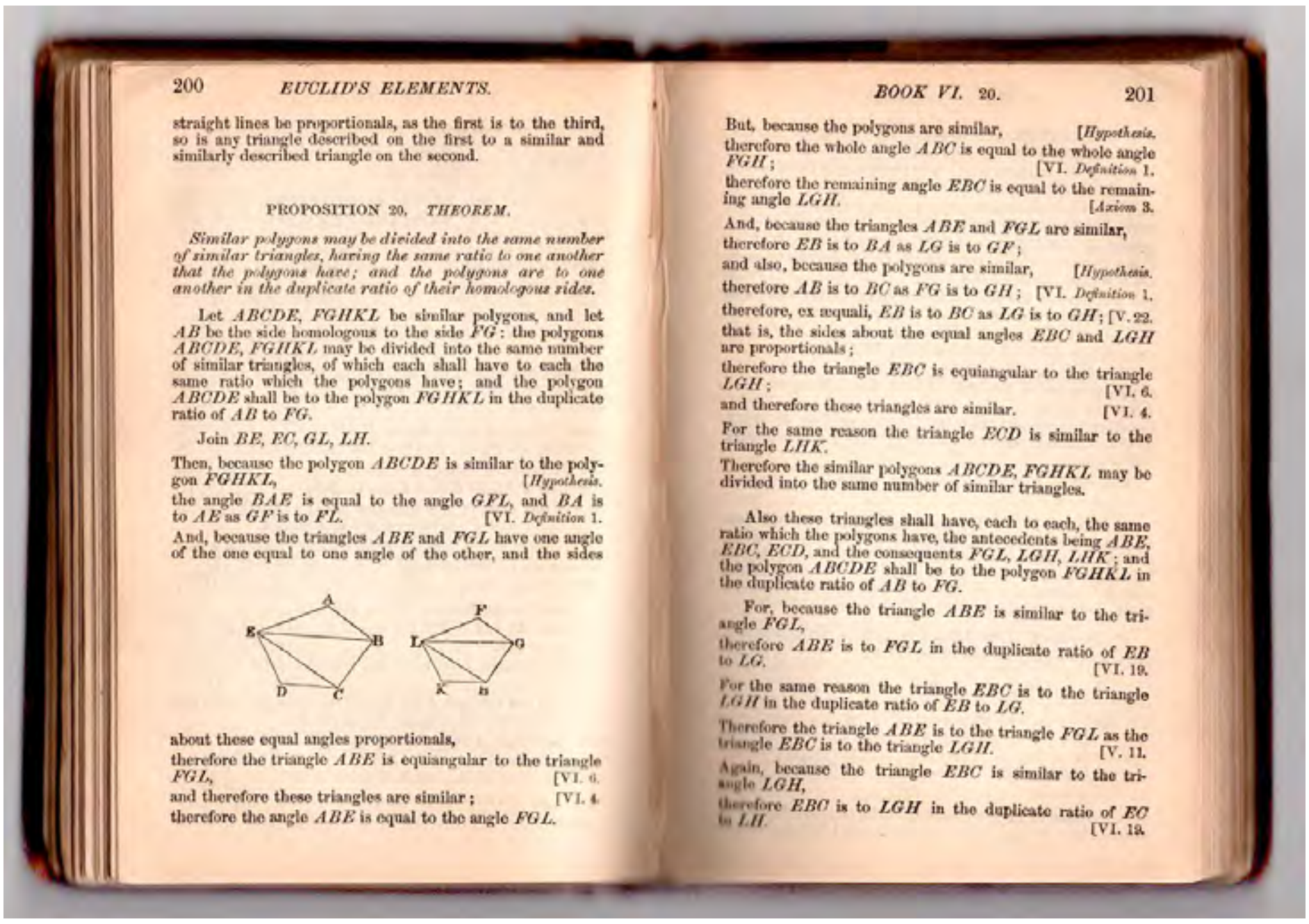

200 *EUCLID'S ELEMENTS.*

straight lines be proportionals, as the first is to the third, so is any triangle described on the first to a similar and similarly described triangle on the second.

PROPOSITION 20. *THEOREM.*

*Similar polygons may be divided into the same number of similar triangles, having the same ratio to one another that the polygons have; and the polygons are to one another in the duplicate ratio of their homologous sides.*

Let *ABCDE*, *FGHKL* be similar polygons, and let *AB* be the side homologous to the side *FG*: the polygons *ABCDE*, *FGHKL* may be divided into the same number of similar triangles, of which each shall have to each the same ratio which the polygons have; and the polygon *ABCDE* shall be to the polygon *FGHKL* in the duplicate ratio of *AB* to *FG*.

Join *BE*, *EC*, *GL*, *LH*.

Then, because the polygon *ABCDE* is similar to the polygon *FGHKL*, [*Hypothesis.*
the angle *BAE* is equal to the angle *GFL*, and *BA* is to *AE* as *GF* is to *FL*. [VI. *Definition* 1.
And, because the triangles *ABE* and *FGL* have one angle of the one equal to one angle of the other, and the sides about these equal angles proportionals,
therefore the triangle *ABE* is equiangular to the triangle *FGL*, [VI. 6.
and therefore these triangles are similar; [VI. 4.
therefore the angle *ABE* is equal to the angle *FGL*.

*BOOK VI.* 20. 201

But, because the polygons are similar, [*Hypothesis.*
therefore the whole angle *ABC* is equal to the whole angle *FGH*; [VI. *Definition* 1.
therefore the remaining angle *EBC* is equal to the remaining angle *LGH*. [*Axiom* 3.
And, because the triangles *ABE* and *FGL* are similar,
therefore *EB* is to *BA* as *LG* is to *GF*;
and also, because the polygons are similar, [*Hypothesis.*
therefore *AB* is to *BC* as *FG* is to *GH*; [VI. *Definition* 1.
therefore, ex æquali, *EB* is to *BC* as *LG* is to *GH*; [V. 22.
that is, the sides about the equal angles *EBC* and *LGH* are proportionals;
therefore the triangle *EBC* is equiangular to the triangle *LGH*; [VI. 6.
and therefore these triangles are similar. [VI. 4.
For the same reason the triangle *ECD* is similar to the triangle *LHK*.
Therefore the similar polygons *ABCDE*, *FGHKL* may be divided into the same number of similar triangles.

Also these triangles shall have, each to each, the same ratio which the polygons have, the antecedents being *ABE*, *EBC*, *ECD*, and the consequents *FGL*, *LGH*, *LHK*; and the polygon *ABCDE* shall be to the polygon *FGHKL* in the duplicate ratio of *AB* to *FG*.

For, because the triangle *ABE* is similar to the triangle *FGL*,
therefore *ABE* is to *FGL* in the duplicate ratio of *EB* to *LG*. [VI. 19.
For the same reason the triangle *EBC* is to the triangle *LGH* in the duplicate ratio of *EB* to *LG*.
Therefore the triangle *ABE* is to the triangle *FGL* as the triangle *EBC* is to the triangle *LGH*. [V. 11.
Again, because the triangle *EBC* is similar to the triangle *LGH*,
therefore *EBC* is to *LGH* in the duplicate ratio of *EC* to *LH*. [VI. 19.

But so what about the metamathematics of Euclid? Given all those theorems—and proofs—can we map out the structure of what Euclid did? That's what the graph in *A New Kind of Science* was about. A few years ago, we put the data for that graph into our Wolfram Data Repository—and I looked at it again, but nothing immediately seemed to jump out about it; it still just seemed like a complicated mess:

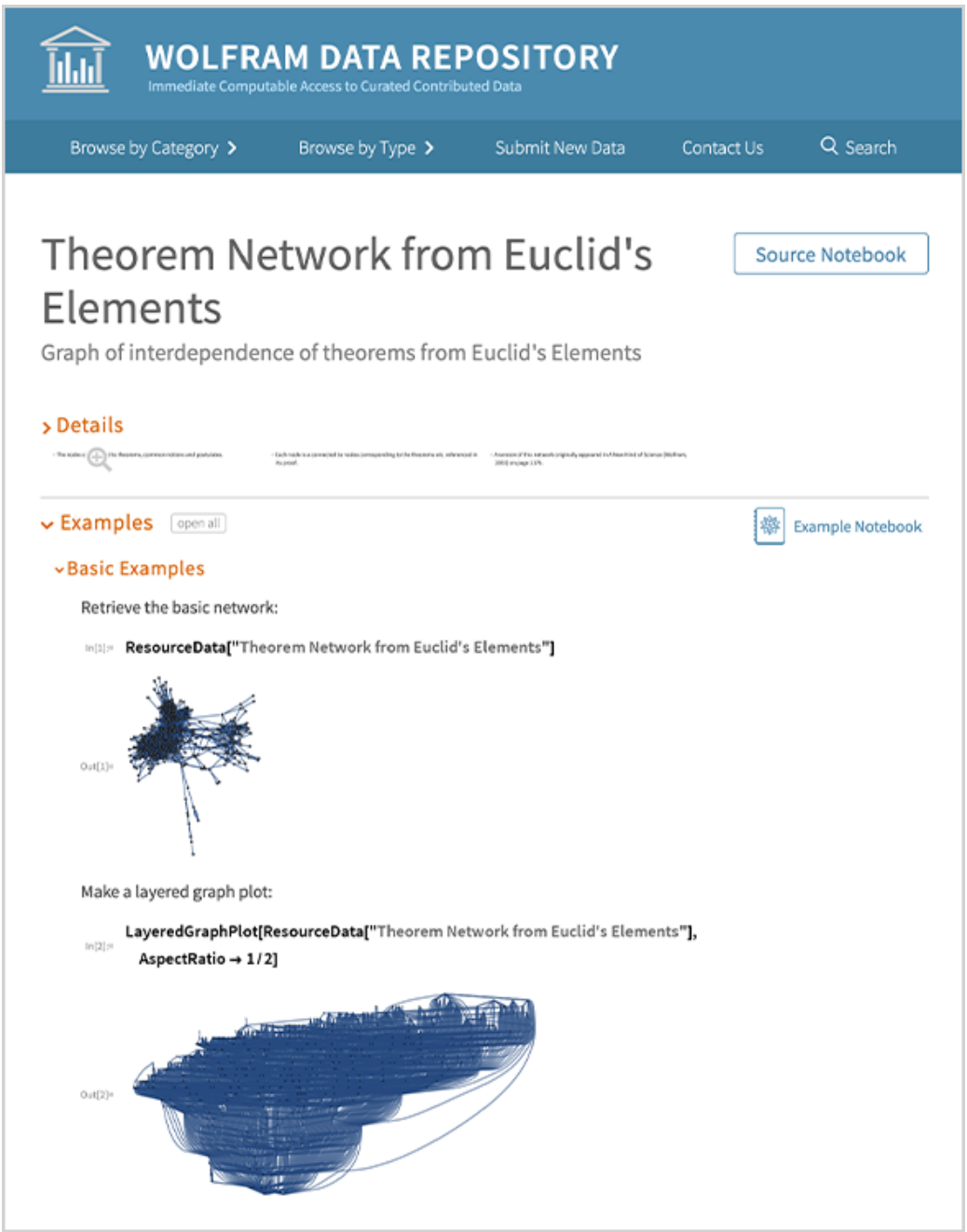

What else happened? One thing is that we added automated theorem proving to Mathematica and the Wolfram Language. Enter a potential theorem, and axioms from which to derive it, and FindEquationalProof will try to generate a proof. This works well for "structurally simple" mathematical systems (like basic logic), and indeed one can generate proofs with complex networks of lemmas that go significantly beyond what humans can do (or readily understand):

```
In[∘]:= FindEquationalProof[p · q == q · p, ∀{a,b,c} ((a · b) · c) · (a · ((a · c) · a)) == c]["ProofGraph"]
```

Out[∘]=

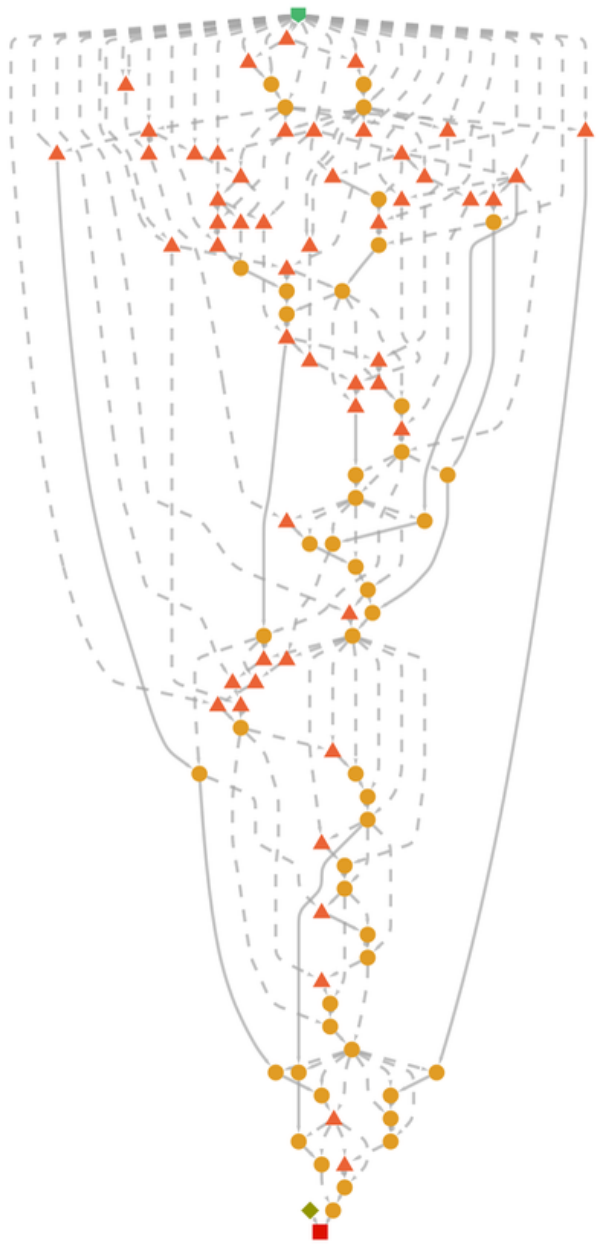

It's in principle possible to use these methods to prove theorems in Euclidean geometry too. But it's a different problem to make the proofs readily understandable to humans (like the step-by-step solutions of Wolfram|Alpha). So at least for now—even after 2000 years—the most effective source of information about the empirical metamathematics of proofs of Euclid's theorems is still basically going to be Euclid's *Elements*.

But when it comes to representing Euclid's theorems there's something new. The whole third-of-a-century story of the Wolfram Language has been about finding ways to represent more and more things in the world computationally. I had long wondered what it would take to represent Euclid-style geometry computationally. And in April I was excited to announce that we'd managed to do it:

With the new #WolfLang it'll be exciting to finally be able to say in computational language what Euclid said in Greek 2300 years ago...

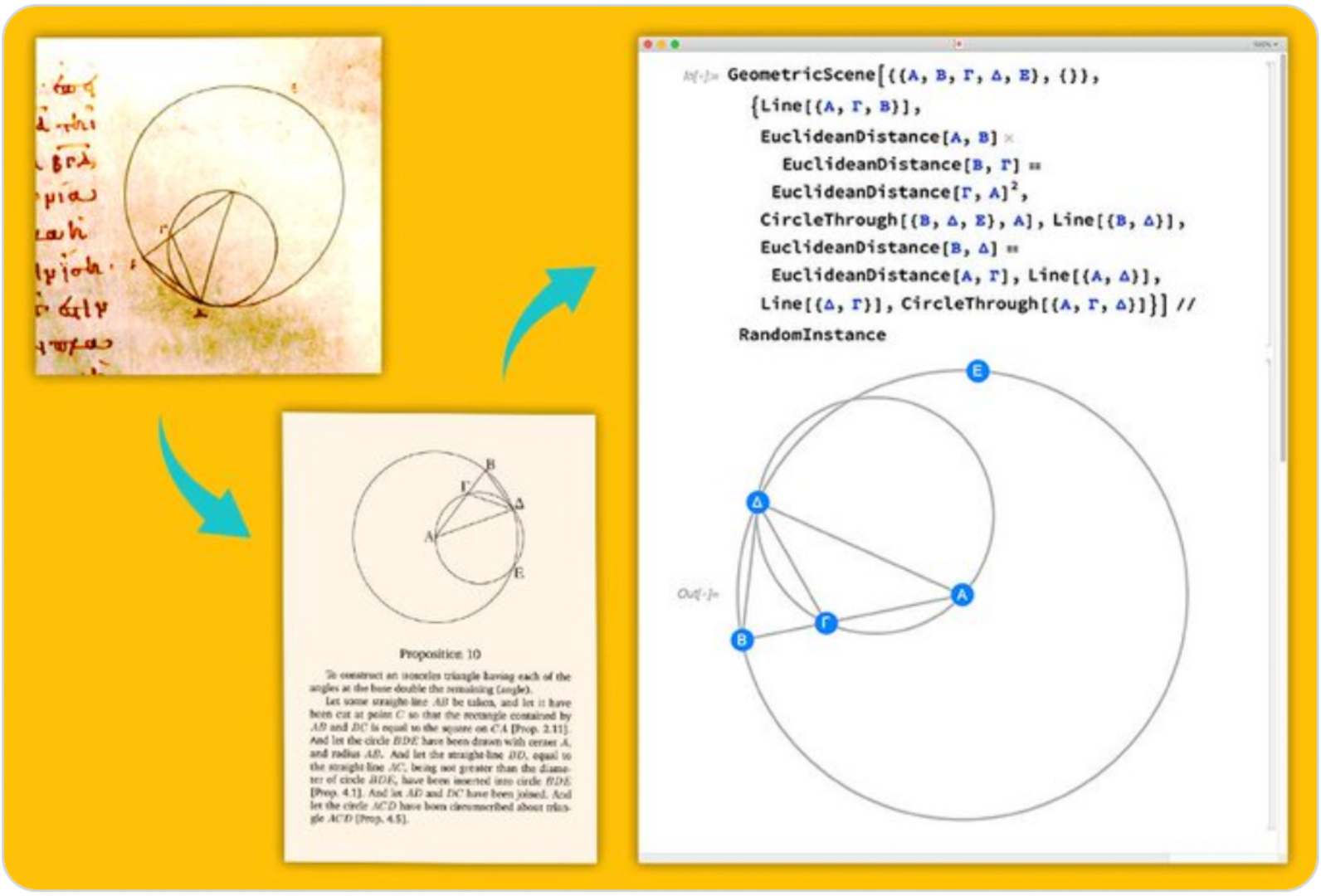

## Basic Statistics of Euclid

Euclid's *Elements* is divided into 13 "books", containing a total of 465 theorems (and 131 definitions):

| subjects | 2D geometry | | | | | | numbers | | | | 3D geometry | | |
|---|---|---|---|---|---|---|---|---|---|---|---|---|---|
| books | 1 | 2 | 3 | 4 | 5 | 6 | 7 | 8 | 9 | 10 | 11 | 12 | 13 |
| theorems | 48 | 14 | 37 | 16 | 25 | 33 | 39 | 27 | 36 | 115 | 39 | 18 | 18 |
| totals | (173) | | | | | | (217) | | | | (75) | | |
| definitions | 23 | 2 | 11 | 7 | 18 | 4 | 22 | 0 | 0 | 16 | 28 | 0 | 0 |
| totals | (65) | | | | | | (38) | | | | (28) | | |

Stating the theorems takes 9589 words (about 60k characters) of Greek (about 13,000 words in a standard English translation). (The 10 axioms take another 115 words in Greek or about 140 in English, and the definitions another 2369 words in Greek or about 3300 in English.)

A typical theorem (or "proposition")—in this case Book 1, Theorem 20—is stated as:

*παντὸς τριγώνου αἱ δύο πλευραὶ τῆς λοιπῆς μείζονές εἰσι πάντῃ μεταλαμβανόμεναι.*

In any triangle two sides taken together in any manner are greater than the remaining one.

(This is what we now call the triangle inequality. And of course, to make this statement we have to have defined what a triangle is, and Euclid does that earlier in Book 1.)

If we look at the statements of Euclid's theorems in Greek (or in English), there's a distribution of lengths (colored here by subjects, and reasonably fit by a Pascal distribution):

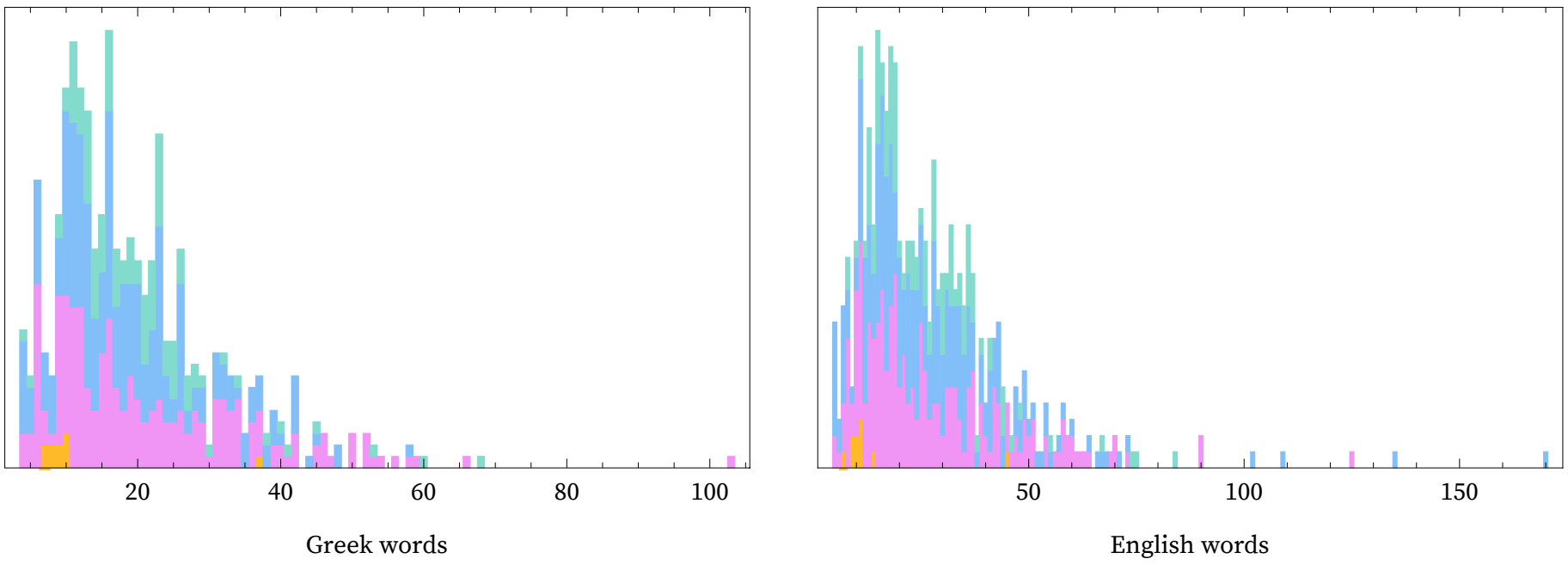

The "outlier" longest-to-state theorem (in both Greek and English) is the rather unremarkable 103-Greek-word 3.8

*ἐὰν κύκλου ληφθῇ τι σημεῖον ἐκτός, ἀπὸ δὲ τοῦ σημείου πρὸς τὸν κύκλον διαχθῶσιν εὐθεῖαί τινες, ὧν μία μὲν διὰ τοῦ κέντρου, αἱ δὲ λοιπαί, ὡς ἔτυχεν, τῶν μὲν πρὸς τὴν κοίλην περιφέρειαν προσπιπτουσῶν εὐθειῶν μεγίστη μέν ἐστιν ἡ διὰ τοῦ κέντρου, τῶν δὲ ἄλλων ἀεὶ ἡ ἔγγιον τῆς διὰ τοῦ κέντρου τῆς ἀπώτερον μείζων ἐστίν, τῶν δὲ πρὸς τὴν κυρτὴν περιφέρειαν προσπιπτουσῶν εὐθειῶν ἐλαχίστη μέν ἐστιν ἡ μεταξὺ τοῦ τε σημείου καὶ τῆς διαμέτρου, τῶν δὲ ἄλλων ἀεὶ ἡ ἔγγιον τῆς ἐλαχίστης τῆς ἀπώτερόν ἐστιν ἐλάττων, δύο δὲ μόνον ἴσαι ἀπὸ τοῦ σημείου προσπεσοῦνται πρὸς τὸν κύκλον ἐφ' ἑκάτερα τῆς ἐλαχίστης.*

If a point be taken outside a circle and from the point straight lines be drawn through to the circle, one of which is through the centre and the others are drawn at random, then, of the straight lines which fall on the concave circumference, that through the centre is greatest, while of the rest the nearer to that through the centre is always greater than the more remote, but, of the straight lines falling on the convex circumference, that between the point and the diameter is least, while of the rest the nearer to the least is always less than the more remote, and only two equal straight lines will fall on the circle from the point, one on each side of the least.

which can be illustrated as:

Out[◦]=

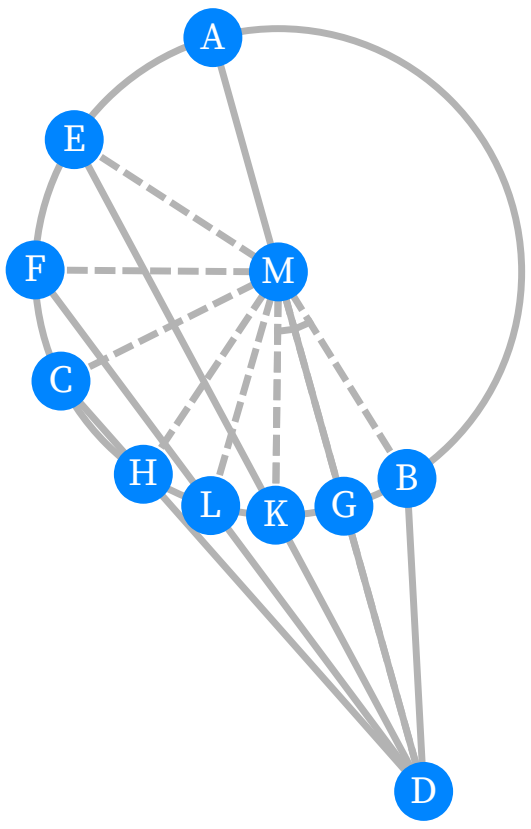

(The runner-up, at about two-thirds the length, is the also rather unremarkable 11.35.)

The nominally shortest-to-state theorems are in Book 10, Theorems 85 through 90, and all have just 4 Greek words:

*εὑρεῖν τὴν πρώτην ἀποτομήν.*

To find the first apotome.

⋮

*εὑρεῖν τὴν ἕκτην ἀποτομήν.*

To find the sixth apotome.

The shortness of these theorems is a bit of a cheat, since the successive “apotomes” (pronounced /əˈpɒtəmi/ like “hippopotamus”) actually have quite long definitions that are given elsewhere. And, yes, some emphasis in math has changed in the past 2000+ years; you don’t hear about apotomes these days. (An apotome is a number $x - y$ where $\frac{x}{y}$ isn’t rational, but $\frac{x^2}{y^2}$ is—as for $x=\sqrt{2}$, $y=1$. It’s difficult enough to describe even this without math notation. But then for a “first apotome” Euclid added the conditions that both $\frac{\sqrt{x^2-y^2}}{x}$ and $x$ must be rational—all described in words.)

At five words, we’ve got one more familiar theorem (3.30) and another somewhat obscure one (10.26):

*τὴν δοθεῖσαν περιφέρειαν δίχα τεμεῖν.*

To bisect a given circumference.

*μέσον μέσου οὐχ ὑπερέχει ῥητῷ.*

A medial area does not exceed a medial area by a rational area.

In our modern Wolfram Language representation, we’ve got a precise, symbolic way to state Euclid’s theorems. But Euclid had to rely on natural language (in his case, Greek). Some words he just assumed people would know the meanings of. But others he defined. Famously, he started at the beginning of Book 1 with his Definition 1—and in a sense changing how we think about this is what launched our whole Physics Project:

*σημεῖόν ἐστιν, οὗ μέρος οὐθέν.*

A point is that which has no part.

There is at least an implicit network of dependencies among Euclid’s definitions. Having started by defining points and lines, he moves on to defining things like triangles, and equilaterality, until eventually, for example, by Book 11 Definition 27 he’s saying things like “An icosahedron is a solid figure contained by twenty equal and equilateral triangles”.

Of course, Euclid didn’t ultimately have to set up definitions; he could just have repeated the content of each definition every time he wanted to refer to that concept. But like words in natural language—or functions in our computational language—definitions are an important form of compression for making statements. And, yes, you have to pick the right definitions to make the things you want to say easy to say. And, yes, your definitions will likely play at least some role in determining what kinds of things you choose to talk about. (Apotomes, anyone?)

## The Interdependence of Theorems

All the theorems Euclid states represent less than 10,000 words of Greek. But the standard proofs of them are perhaps 150,000 words of Greek. (They're undoubtedly not minimal proofs—but the fact that the same ones are being quoted after more than 2000 years presumably tells us at least something.)

Euclid is very systematic. Every theorem throughout the course of his *Elements* is proved in terms of earlier theorems (and ultimately in terms of his 10 axioms). Thus, for example, the proof of 1.14 (i.e. Book 1, Theorem 14) uses 1.13 as well as the axioms P2 (i.e. Postulate 2), P4, CN1 (i.e. Common Notion 1) and CN3. By the time one's got to 12.18 the proof is written only in terms of other theorems (in this case 12.17, 12.2, 5.14 and 5.16) and not directly in terms of axioms.

The total number of theorems (or axioms) directly referenced in a given proof varies from 0 (for axioms) to 21 (for 12.17, which is about inscribing polyhedra in spheres); the average is 4.3:

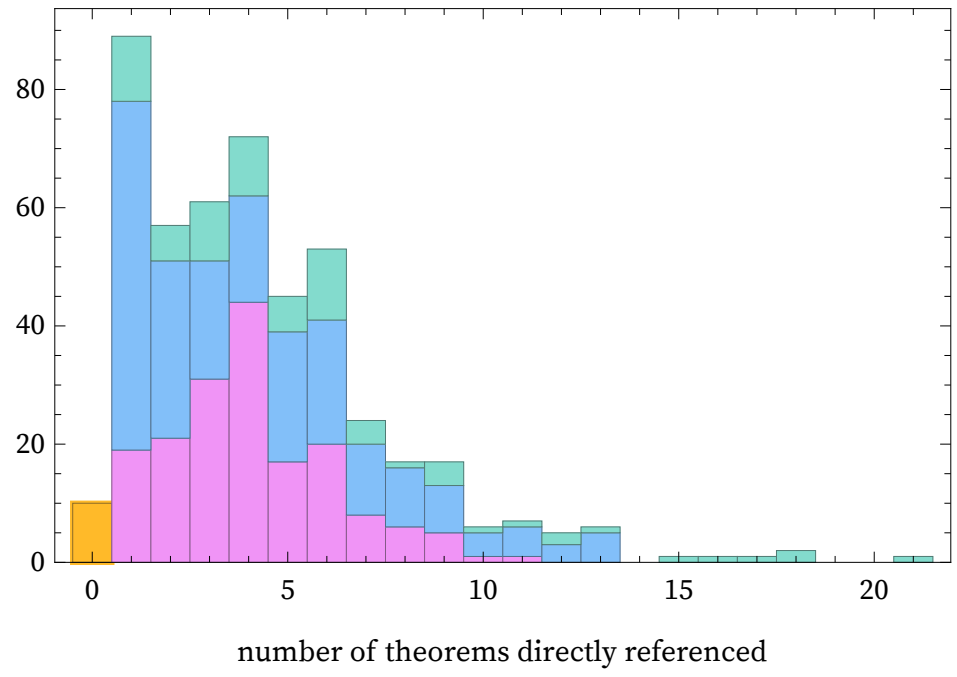

If we put Euclid's axioms and theorems in order, we can represent which axioms or theorems occur in a given proof by an arrangement of dots across the page. For example, for 1.12 through 1.17 we have:

| | CN1 | CN2 | CN3 | CN4 | CN5 | P1 | P2 | P3 | P4 | P5 | 1.1 | 1.2 | 1.3 | 1.4 | 1.5 | 1.6 | 1.7 | 1.8 | 1.9 | 1.10 | 1.11 | 1.12 | 1.13 | 1.14 | 1.15 | 1.16 |
|---|---|---|---|---|---|---|---|---|---|---|---|---|---|---|---|---|---|---|---|---|---|---|---|---|---|---|
| 1.12 | | | | | | • | | • | | | | | | | | | | • | | • | | | | | | |
| 1.13 | • | • | | | | | | | | | | | | | | | | | | | • | | | | | |
| 1.14 | • | | • | | | | • | | • | | | | | | | | | | | | | | • | | | |
| 1.15 | • | | • | | | | | | • | | | | | | | | | | | | | | • | | | |
| 1.16 | | | | | • | • | • | | | | | | • | • | | | | | | • | | | | | • | |
| 1.17 | | | | | | | • | | | | | | | | | | | | | | | | • | | | |

Doing this for all the theorems we get:

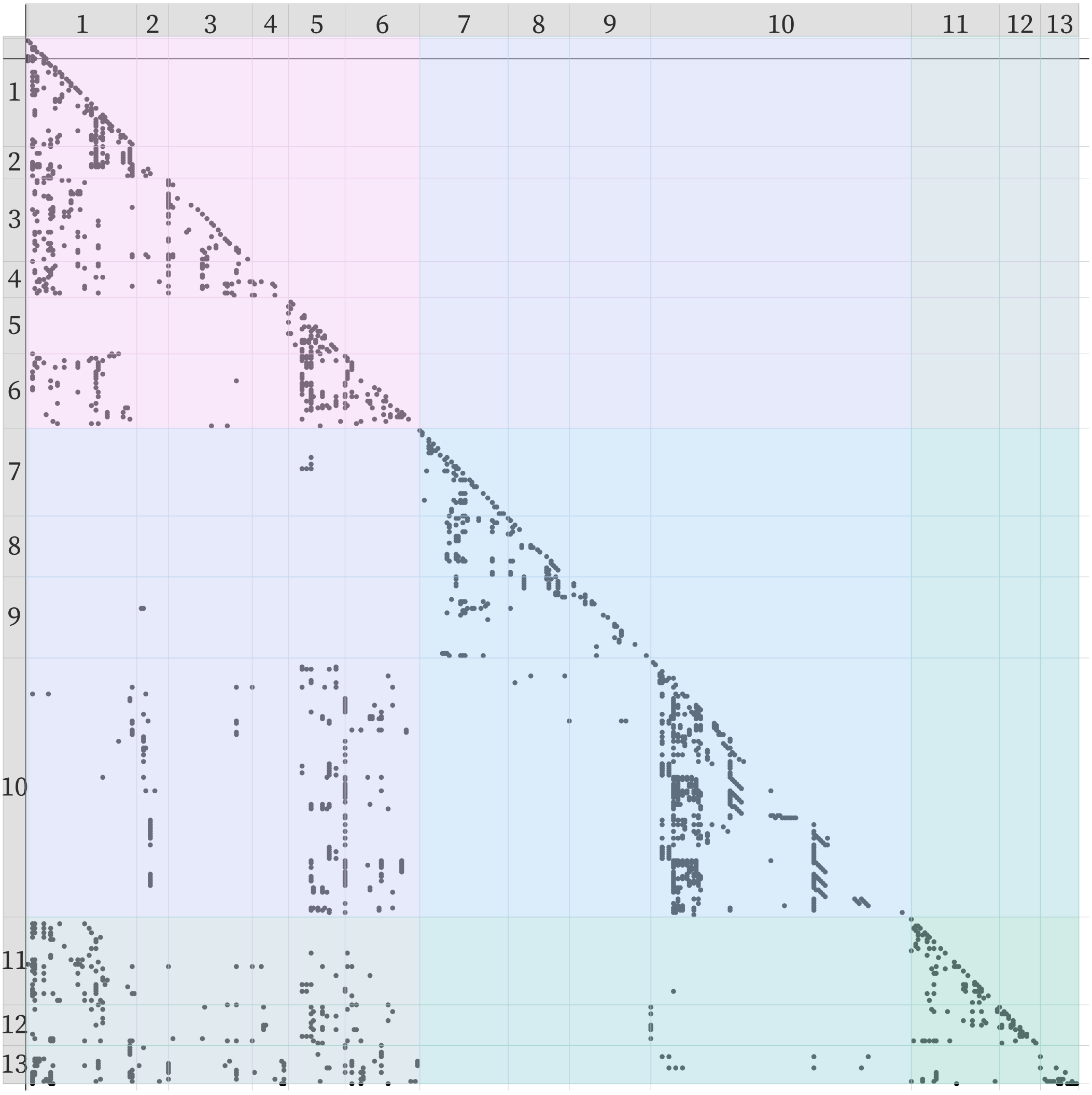

We can see there's lots of structure here. For example, there are clearly "popular" theorems near the beginning of Book 6 and Book 10, to which lots of at least "nearby" theorems refer. There are also "gaps": ranges of theorems that no theorems in a given book refer to.

At a coarse level, something we can do is to look at cross-referencing within and between books:

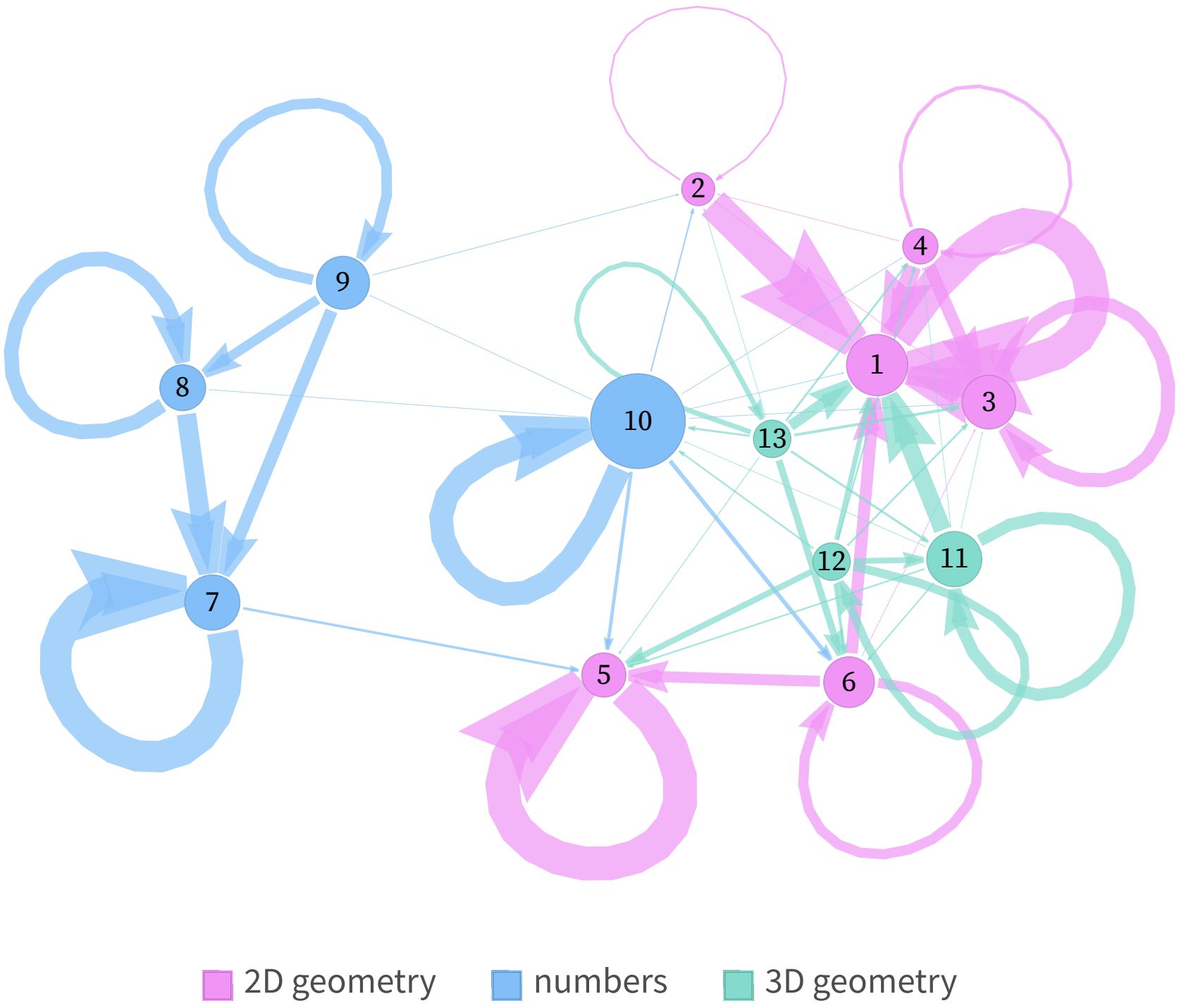

The size of each node represents the number of theorems in each book. The thickness of each arrow represents the fraction of references in the proofs of those theorems going to different books. The self-loops are from theorems in a given book that refer to theorems in the same book. Needless to say, the self-loop is large for Book 1, since it doesn't have any previous book to refer to. Book 7 again has a large self-loop, because it's the first book about numbers, and doesn't refer much to the earlier books (which are about 2D geometry).

It's interesting to see that Books 7, 8 and 9—which are about numbers rather than geometry—"keep to themselves", even though Book 10, which is also about numbers, is more central. It's also interesting to see the interplay between the books on 2D and 3D geometry over on the right-hand side of the graph.

But, OK, what about individual theorems? What is their network of dependencies?

Here's 1.5, whose proof is given in terms of 1.3 and 1.4, as well as the axioms P1, P2 and CN3:

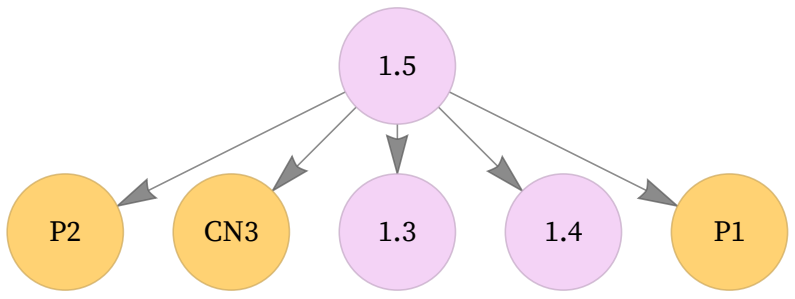

But now we can continue this, and show what 1.3 and 1.4 depend on—all the way down to the axioms:

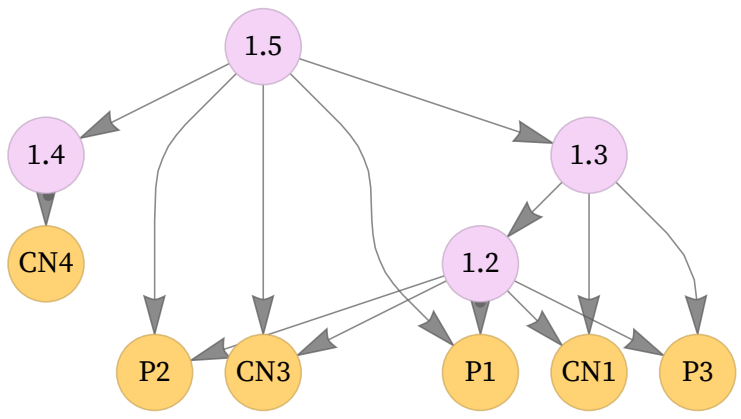

Later theorems depend on much more. Here are the direct dependencies for 12.18:

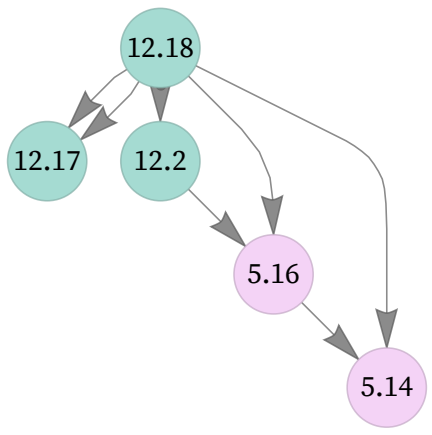

Here's what happens if one goes another step:

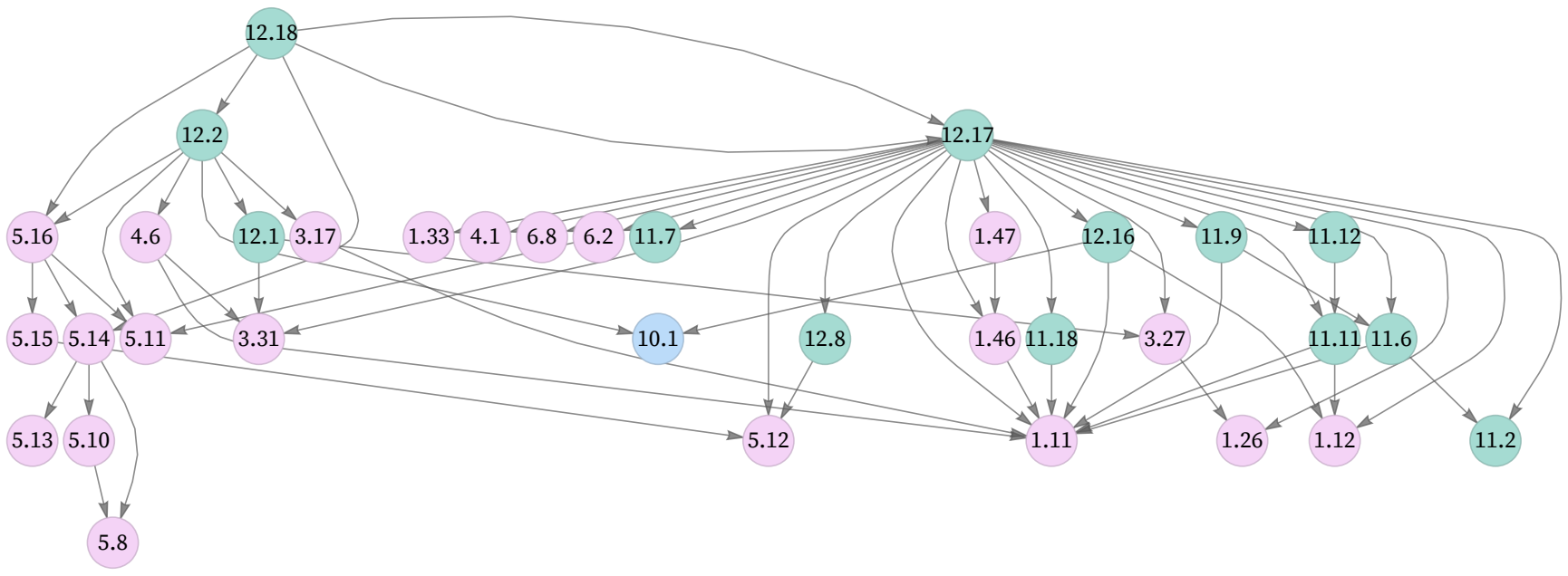

Here's 3 steps:

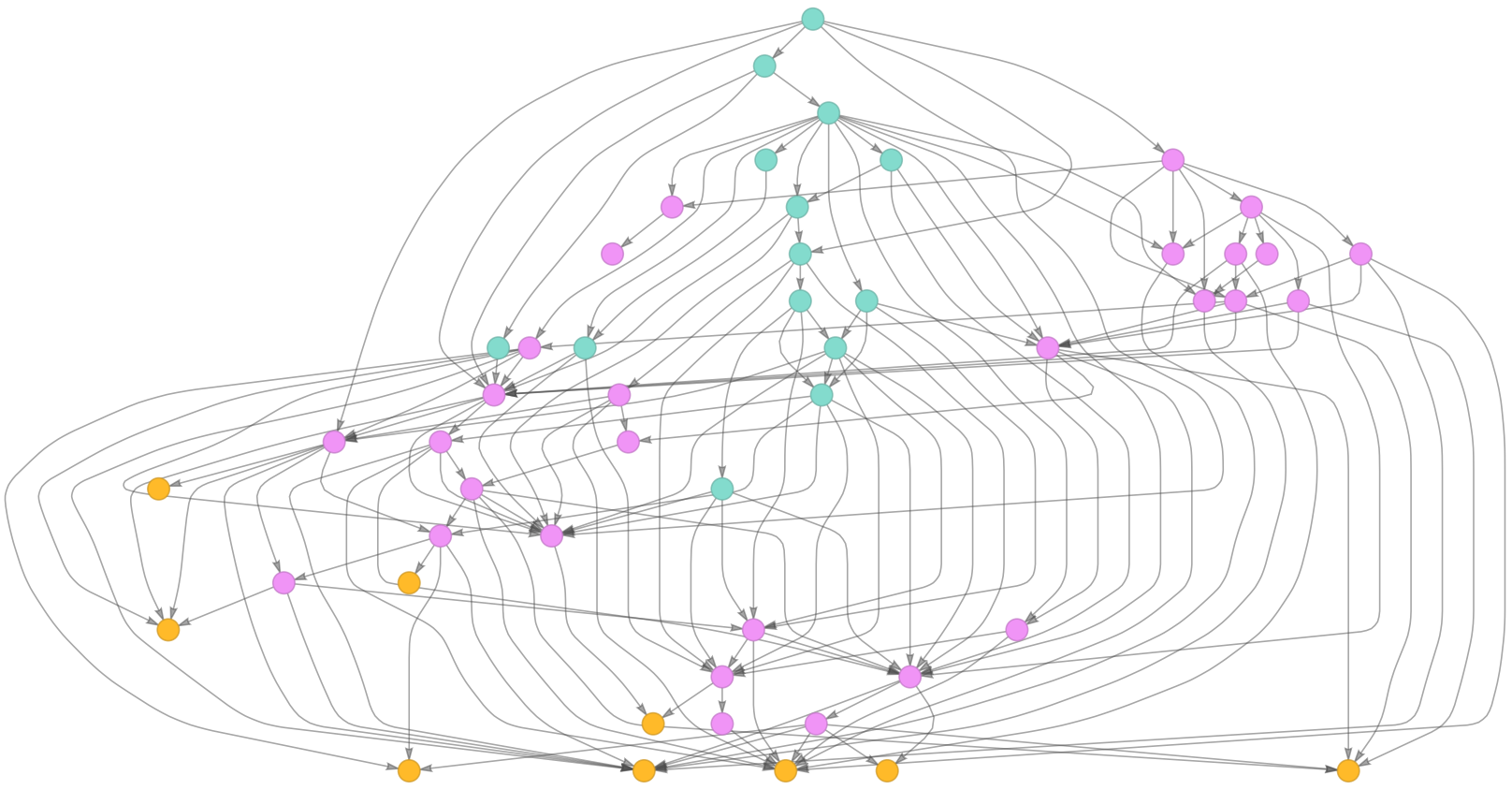

And here's what happens if one goes all the way down to the axioms (which in this case takes 5 steps):

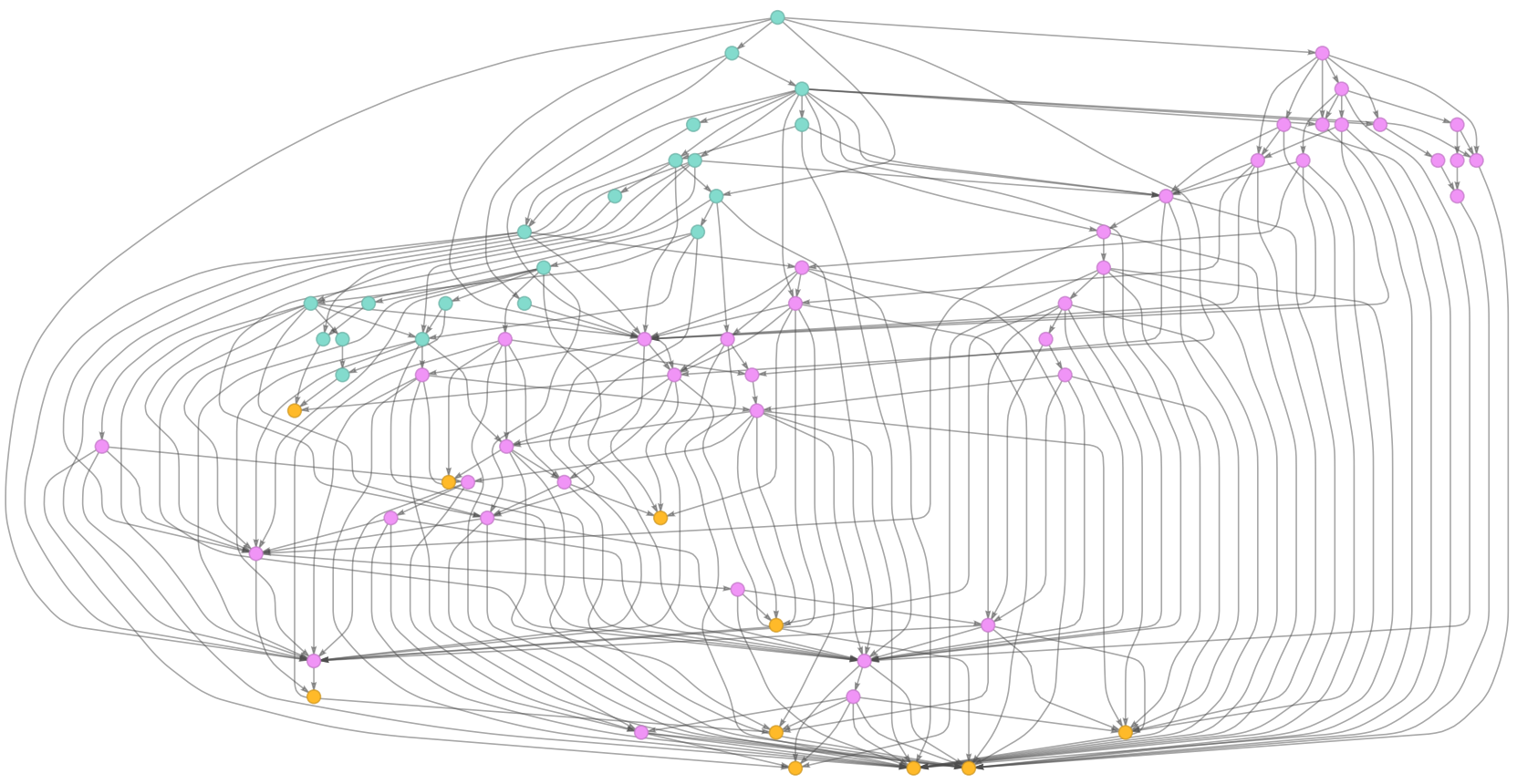

Things look a little simpler if we consider the transitive reduction of this graph. We're no longer faithfully representing what's in the text of Euclid, but we're still capturing the core dependency information. If theorem A in Euclid refers to B, and B refers to C, then even if in Euclid A refers to C we won't mention that. And, yes, graph theoretically A→C is just the transitive closure of A→B and B→C. But it could still be that the pedagogical structure of the proof of theorem A makes it desirable to refer to theorem B, even if in principle one could rederive theorem B from theorem C.

Here's the original 1-step graph for 12.18, along with its transitive reduction:

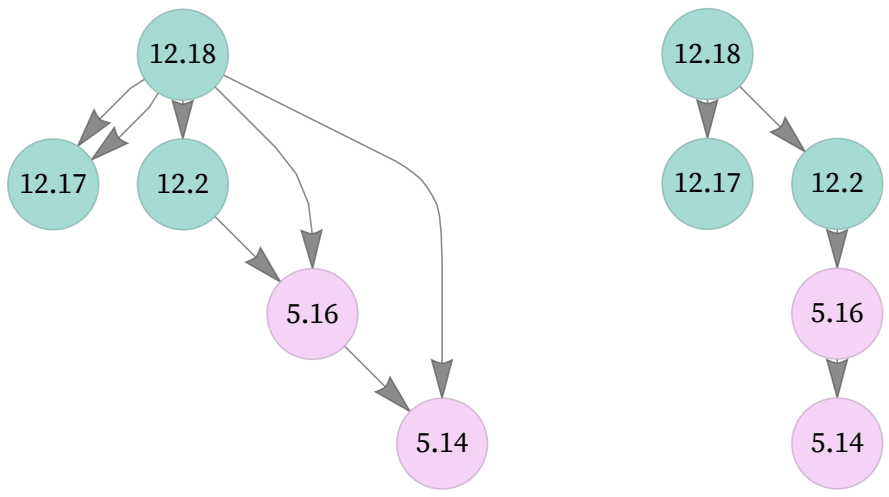

And here, by the way, is also the "fully pedantic" transitive closure, including all indirect connections, whether they're mentioned by Euclid or not:

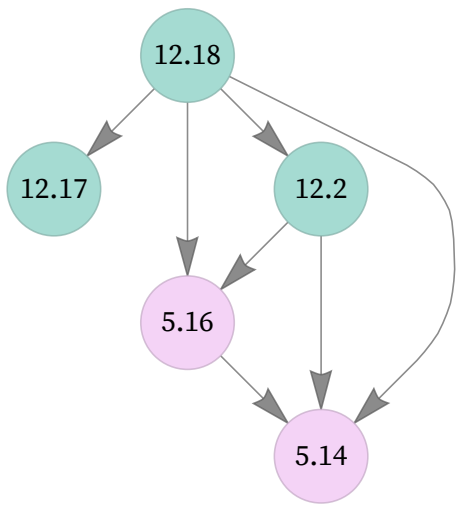

And now here's the transitive reduction of the full 12.8 dependency graph, all the way down to the axioms:

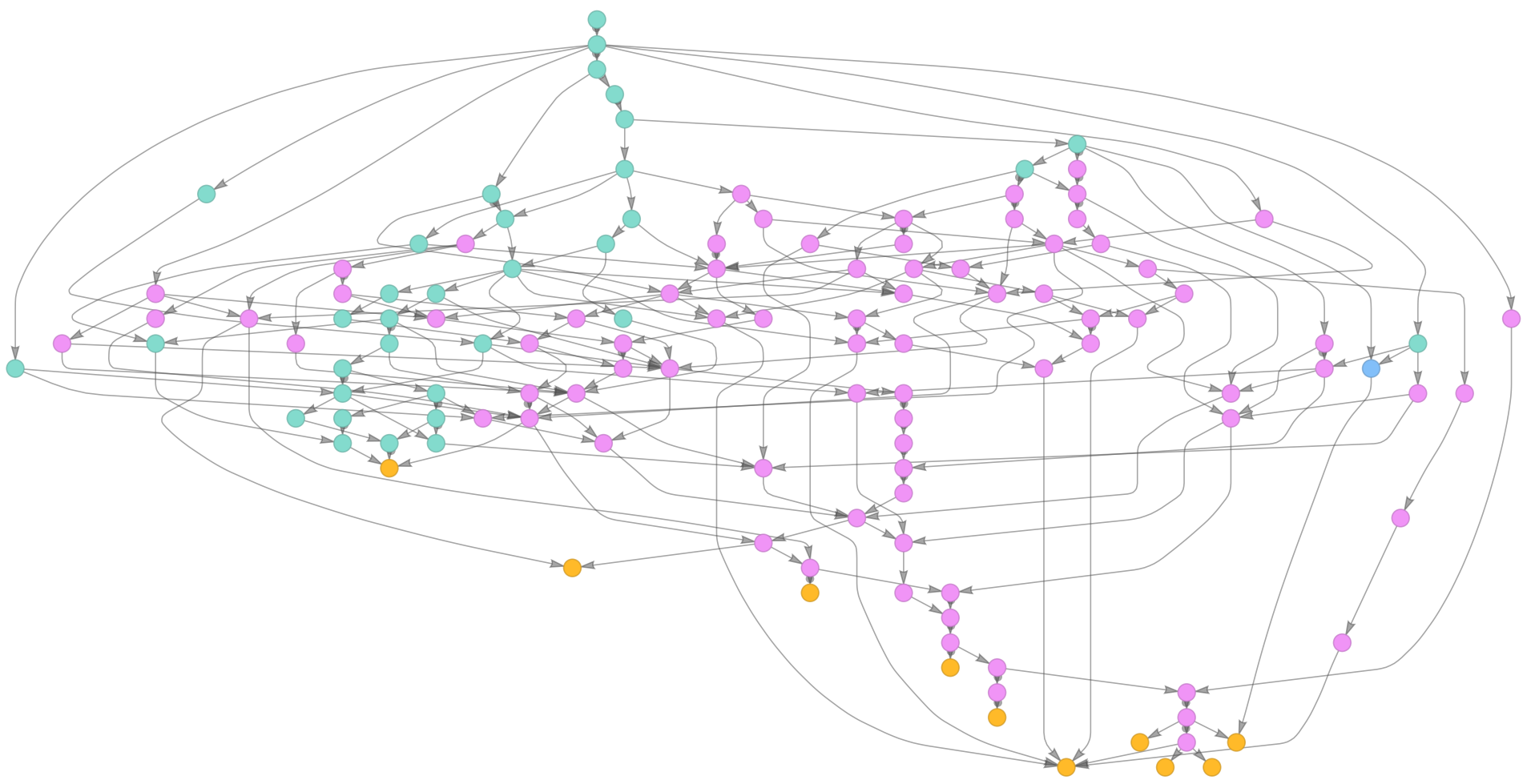

And what all these graphs show is that even to prove one theorem, one's making use of lots of other theorems. To make this quantitative, we can plot the total number of theorems that appear anywhere in the "full proof" of a given theorem, ultimately working all the way down to the axioms:

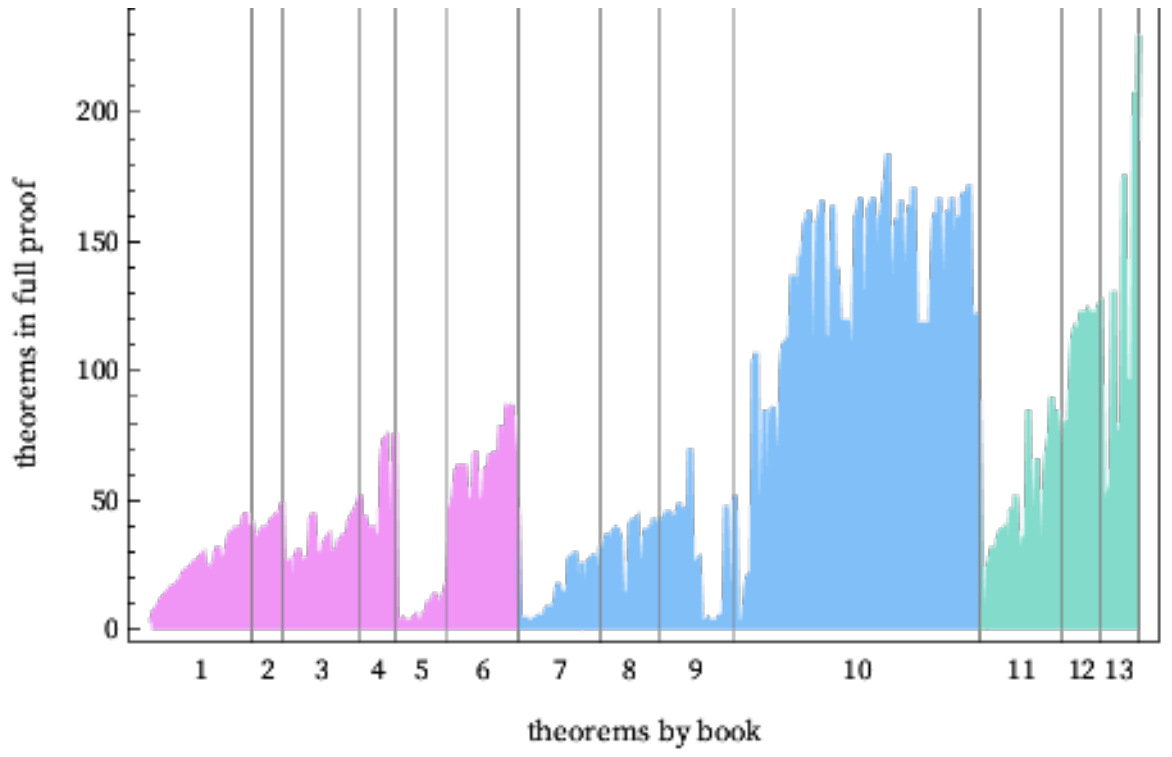

At the beginnings of many of the books, there tend to be theorems that are proved more directly from the axioms, so they don't depend on as much. But as one progresses through the books, one's relying on more and more theorems—sometimes, as we saw above, in the same book, and sometimes in earlier books.

From the picture above, we can see that Euclid in a sense builds up to a "climax" at the end—with his very last theorem (13.18) depending on more theorems than anything else. We'll be discussing "Euclid's last theorem" some more below...

# The Graph of All Theorems

OK, so what is the full interdependence graph for all the theorems in Euclid? It's convenient to go the opposite way than in our previous graphs—and put the axioms at the top, and show how theorems below are derived from them. Here's the graph one gets by doing that:

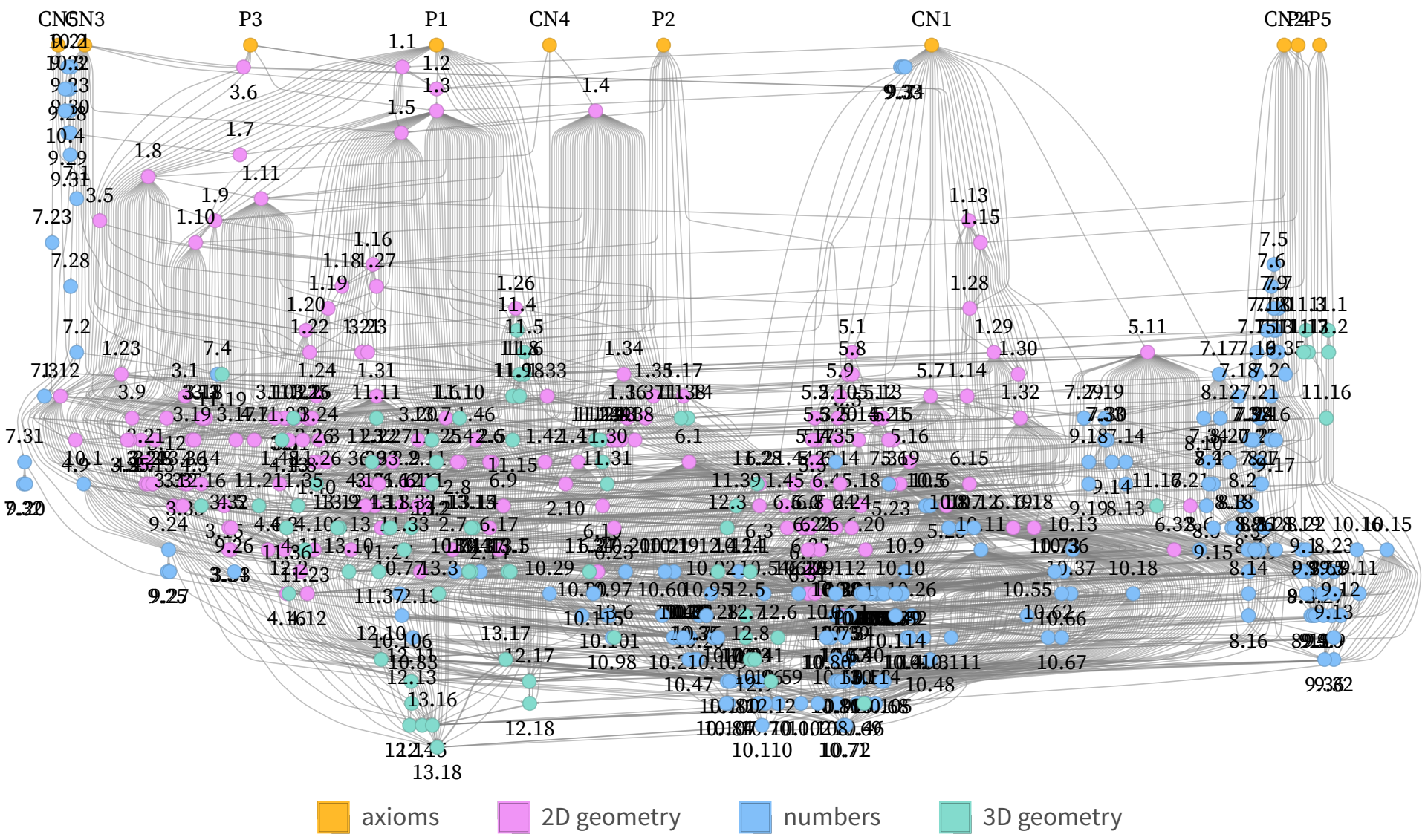

One can considerably simplify this by looking just at the transitive reduction graph (the full graph has 2054 connections; this reduction has 974, while if one went "fully pedantic" with transitive closure, one would have 25,377 connections):

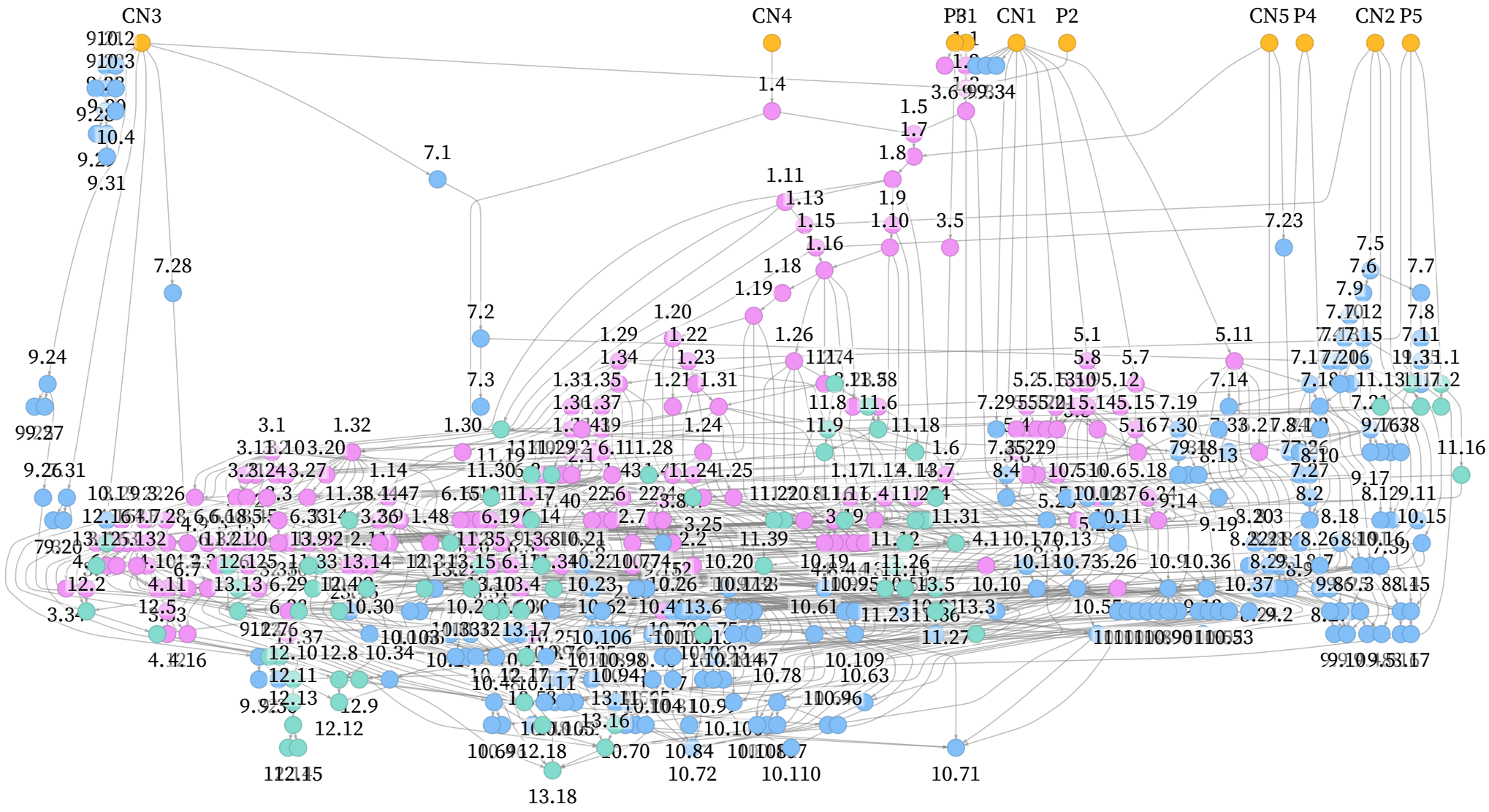

What can we see from this? Probably the most obvious thing is that the graphs start fairly sparse, then become much denser. And what this effectively means is that one starts off by proving certain "preliminaries", and then after one's done that, it unlocks a mass of other theorems. Or, put another way, if we were exploring this metamathematical space starting from the axioms, progress might seem slow at first. But after proving a bunch of preliminary theorems, we'd be able to dramatically speed up.

Here's another view of this, plotting how many subsequent theorems depend on each different theorem:

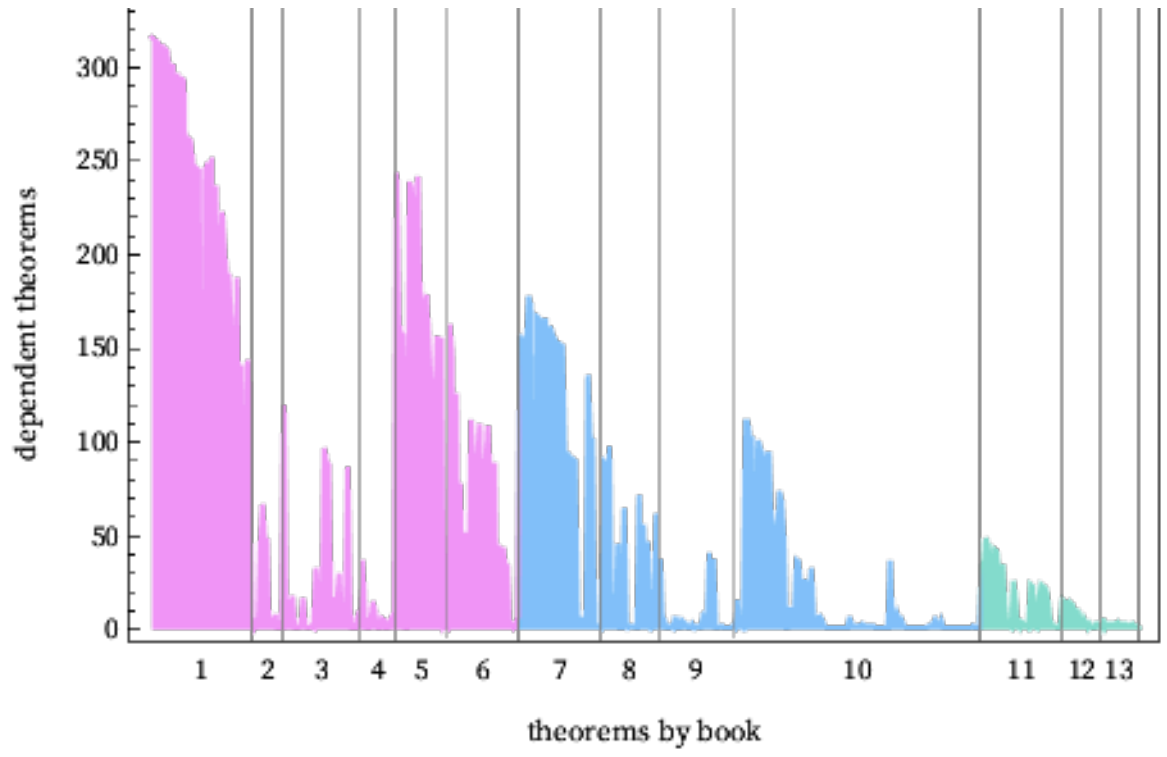

In a sense, this is complementary to the plot we made above, that showed how many theorems a given theorem depends on. (From a graph-theoretical point of view they're very directly complementary: this plot involves `VertexInComponent`; the previous one involved `VertexOutComponent`.)

And what the plot shows is that there are a bunch of early theorems (particularly in Book 1) that have lots of subsequent theorems depending on them—so that they're effectively foundational to much of what follows. The plot also shows that in most of the books the early theorems are the most "foundational", in the sense that the most subsequent theorems depend on them.

By the way, we can also look at the overall form of the basic dependency graph, not layering it starting from the axioms:

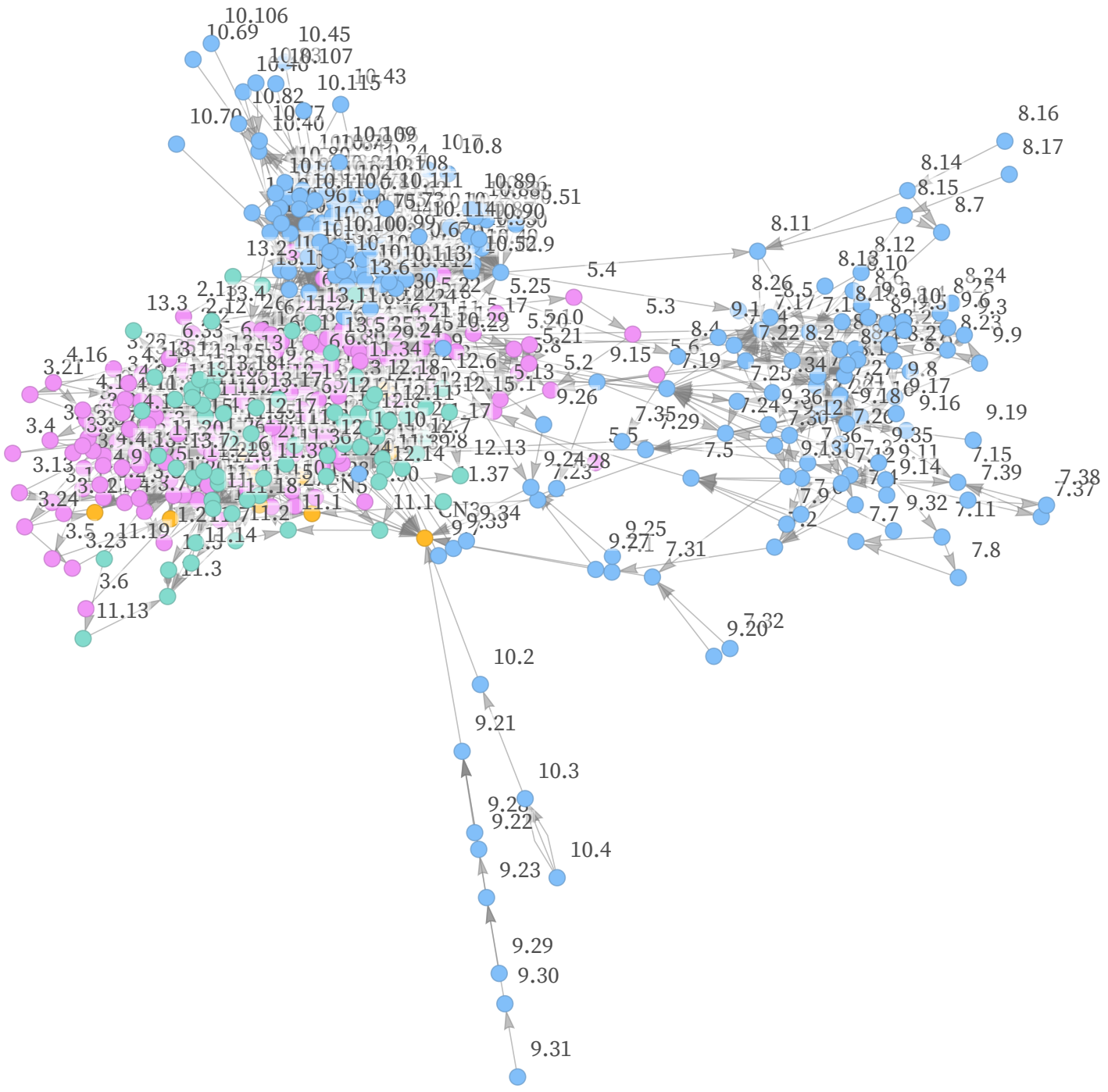

The transitive reduction is slightly easier to interpret:

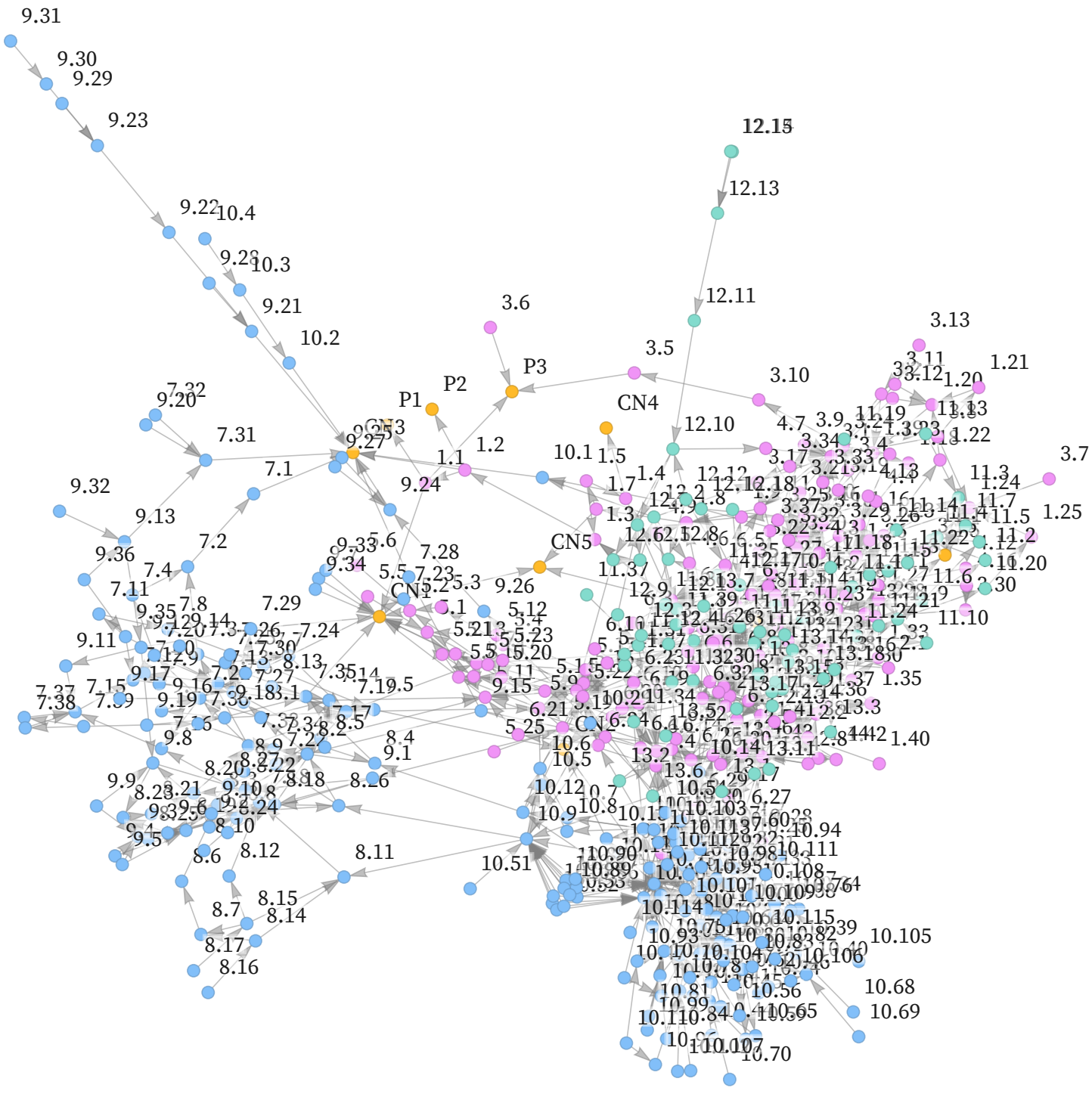

And the main notable feature is the presence of "prongs" associated, for example, with Book 9 theorems about the properties of even and odd numbers.

## The Causal Graph Analogy

Knowing about the Wolfram Physics Project, there's an obvious analog of theorem dependency graphs: they're like causal graphs. You start from a certain set of "initial events" (the "big bang"), corresponding to the axioms. Then each subsequent theorem is like an event, and the theorem dependency graph is tracing out the causal connections between these events.

Just like the causal graph, the theorem dependency graph defines a partial ordering: you can't write down the proof of a given theorem until the theorems that will appear in it have been proved. Like in the causal graph, one can define light cones: there's a certain set of "future" theorems that can be affected by any given theorem. Here is the "future light cone" of Book 1, Theorem 5:

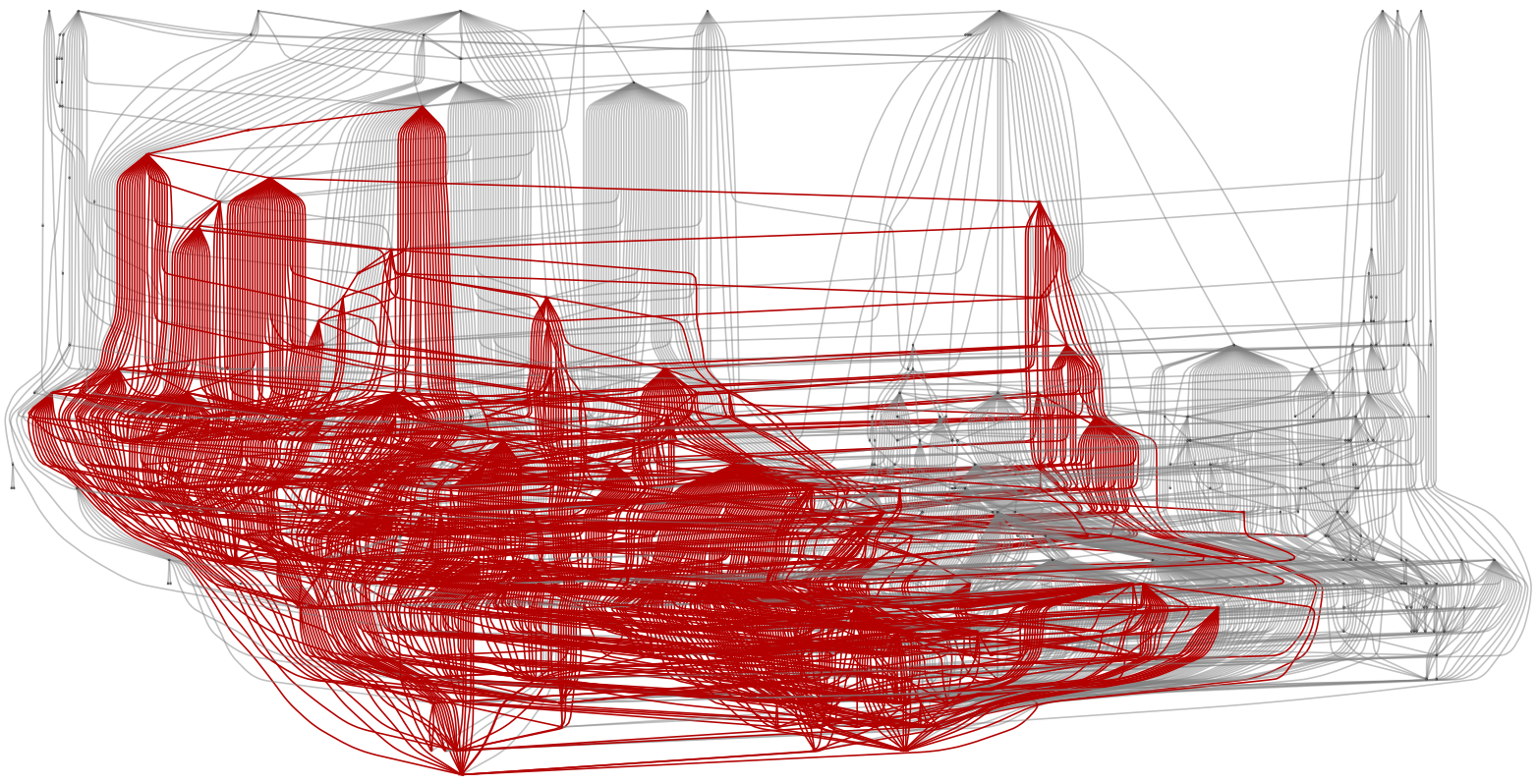

And here is the corresponding transitive reduction graph:

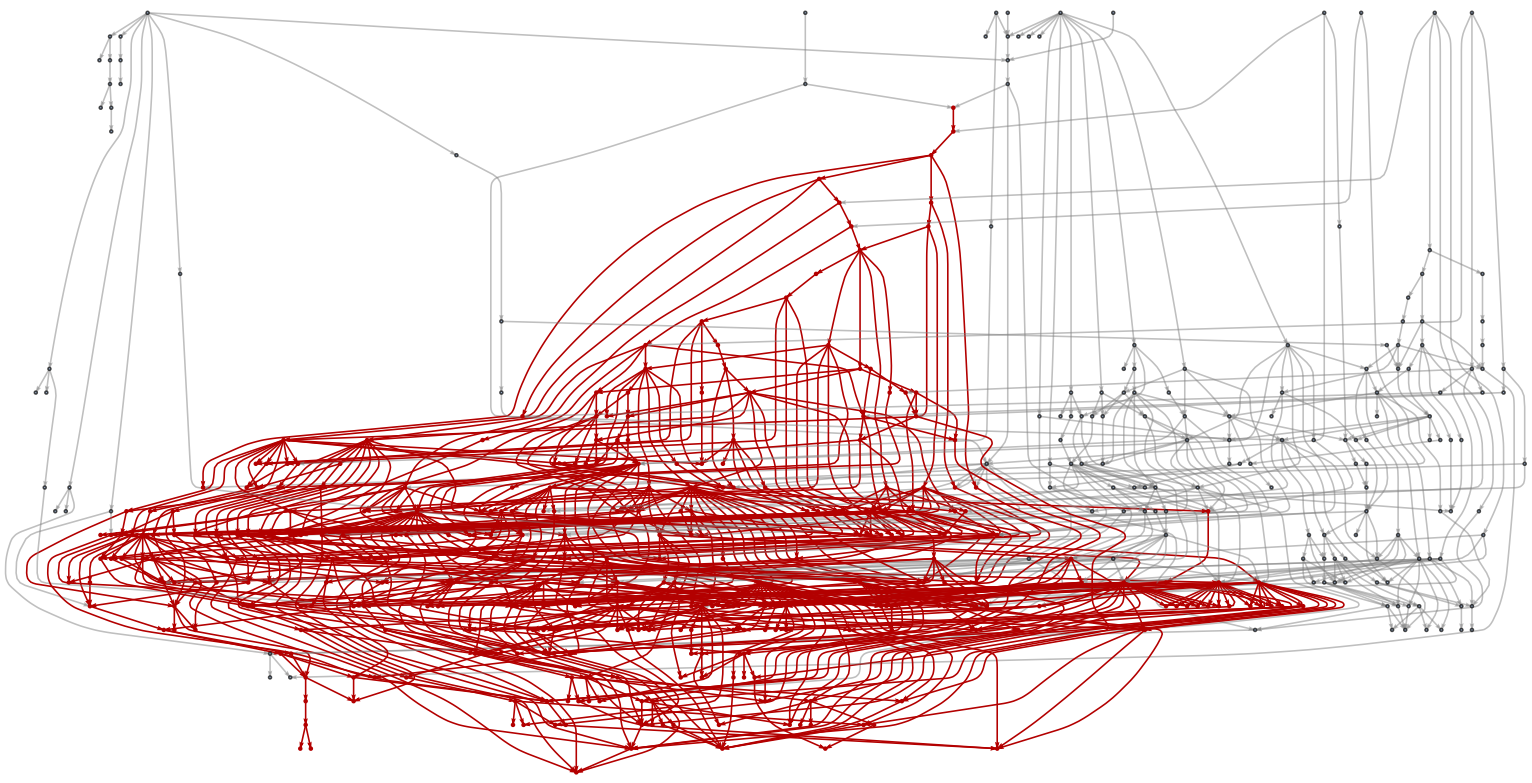

But now let's think about the notion of time in the theorem dependency graph. Imagine you were rederiving the theorems in Euclid in a series of "time steps". What would you have to do at each time step? The theorem dependency graph tells you what you will have to have done in order to derive a particular theorem. But just like for spacetime causal graphs, there are many different foliations one can use to define consistent time steps.

Here's an obvious one, effectively corresponding to a "cosmological rest frame" in which at each step one "does as much as one consistently can at that step":

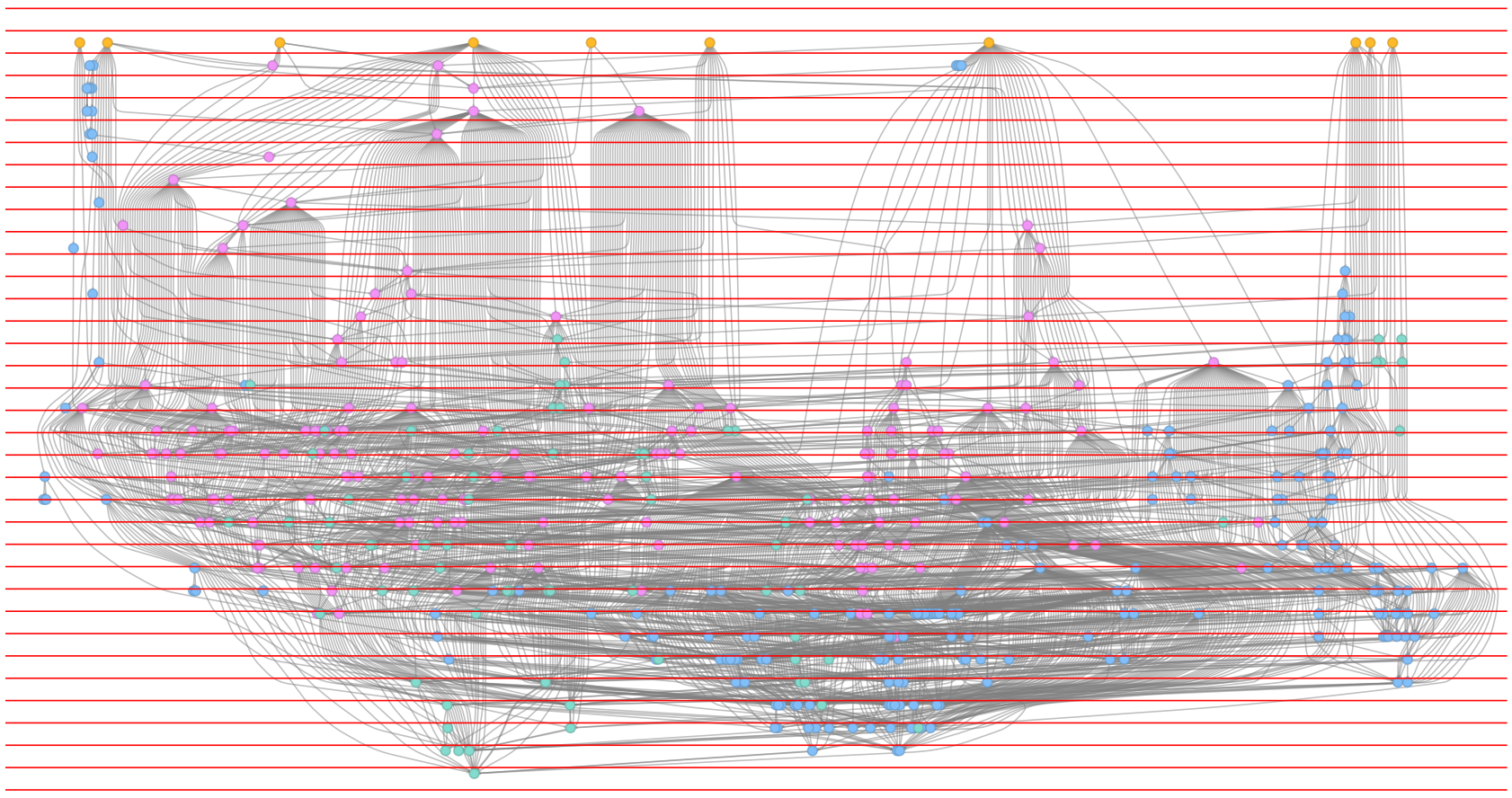

And here are the number of theorems that appear on each slice (in effect each theorem appears on the slice determined by its longest path to any axiom):

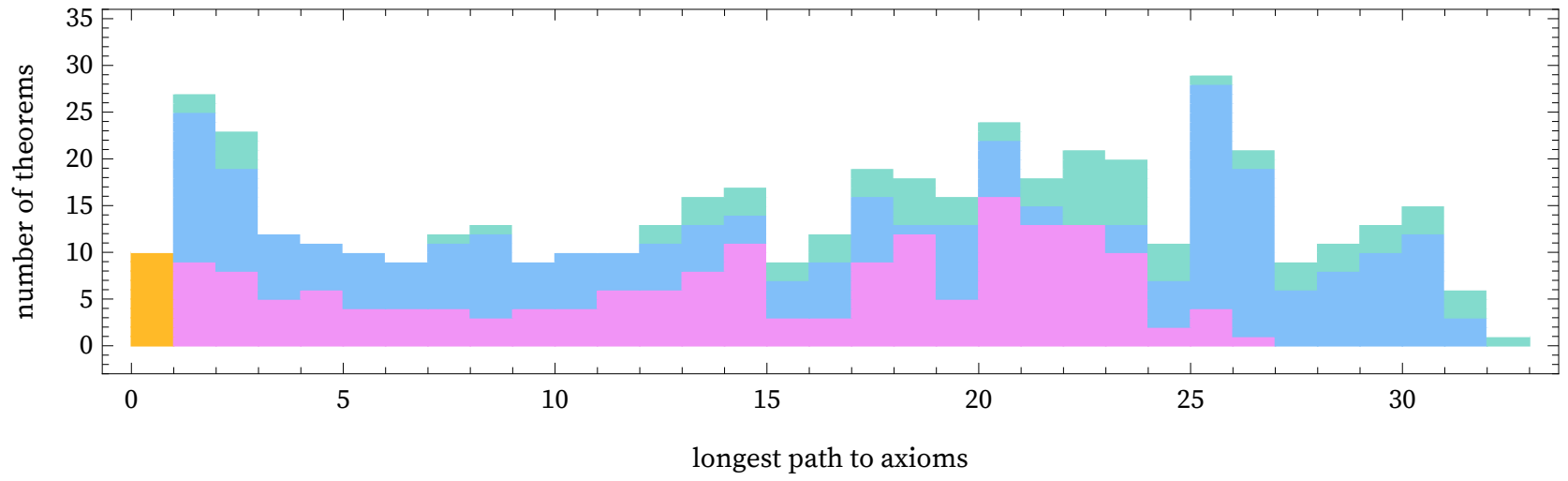

But there are many other foliations that are possible, in which one for example concentrates first on a particular group of theorems, only doing others when one "needs to".

Each choice of foliation can be thought of as corresponding to a different reference frame—and a different choice of how one explores the analog of spacetime in Euclid. But, OK, if the foliations define successive moments in time—or successive "simultaneity surfaces"—what is the analog of space? In effect, the "structure of space" is defined by the way that theorems are laid out on the slices defined by the foliations. And a convenient way to probe this is to look at branchial graphs, in which pairs of theorems on a given slice are connected by an edge if they have an immediate common ancestor on the slice before.

So here are the branchial graphs for all successive slices of Euclid in the "cosmological rest frame":

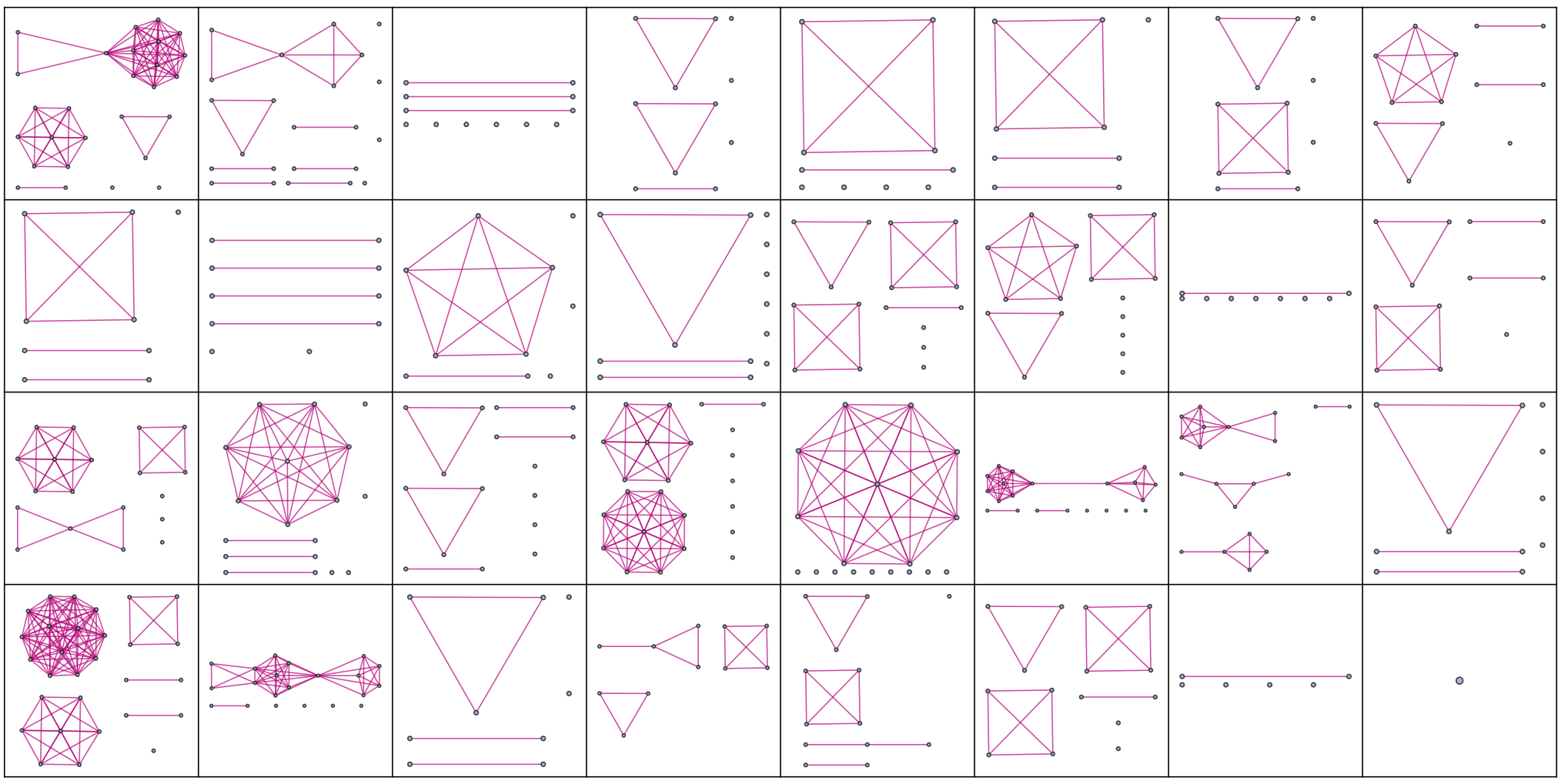

And here are the branchial graphs specifically from slices 23 and 26:

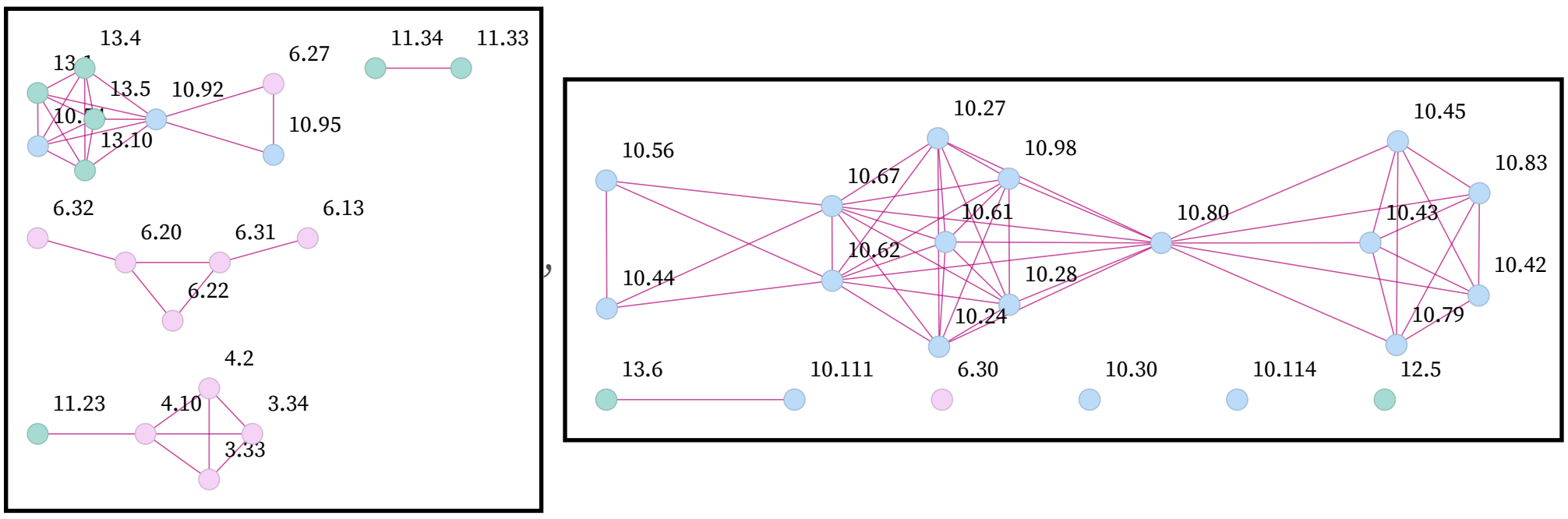

How should we interpret these graphs? Just like in quantum mechanics, they effectively define a map of “entanglements”, but now these are “entanglements” not between quantum states but between theorems. But potentially we can also interpret these graphs as showing how theorems are laid out in a kind of “instantaneous metamathematical space”—or, in effect, we can use the graphs to define “distances between theorems”.

We can generalize our ordinary branchial graphs by connecting theorems that have common ancestors not just one slice back, but also up to $\delta t$ slices back. Here are the results for slice 26 (in the cosmological rest frame):

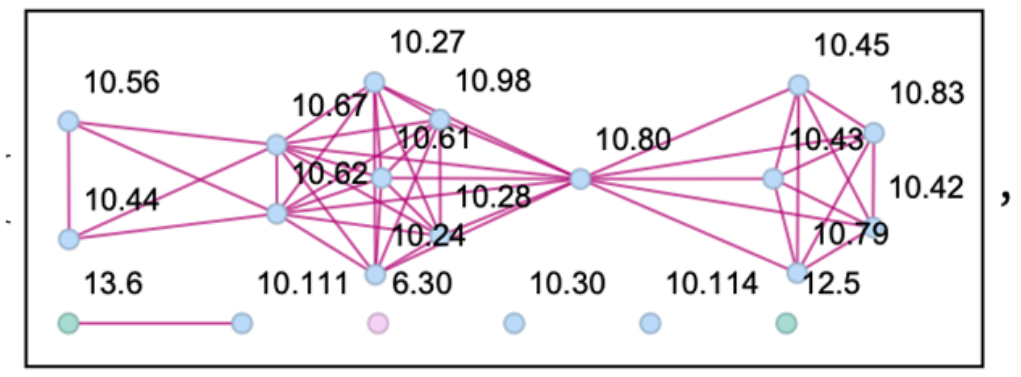

$\delta t = 1$

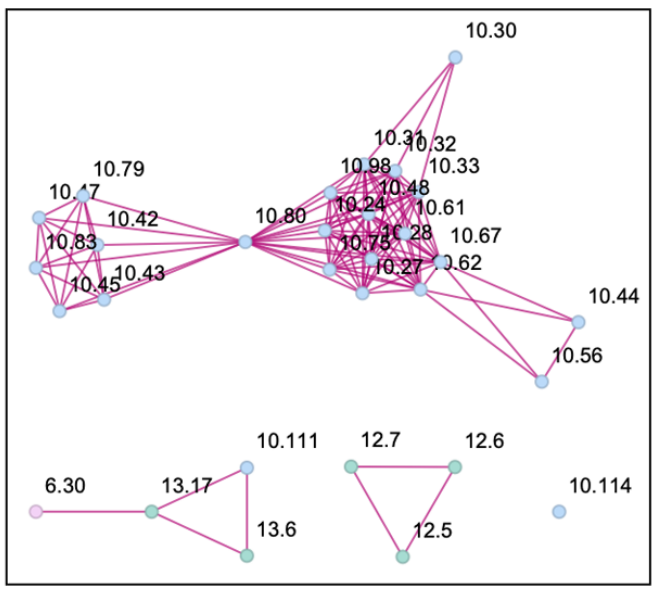

δt = 2

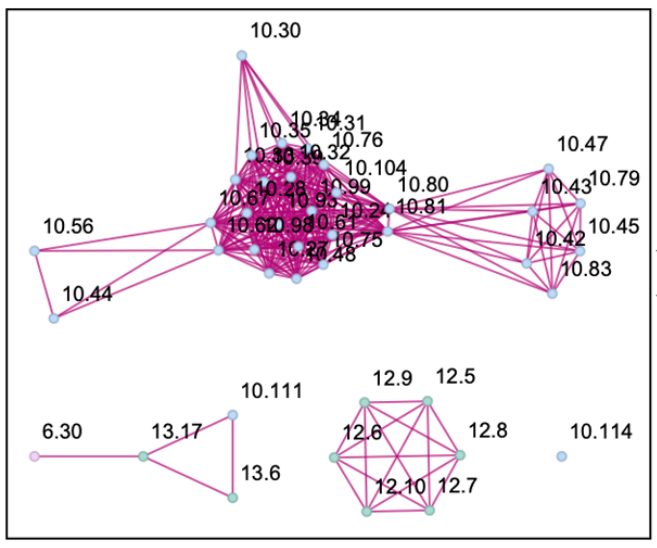

δt = 3

If we went all the way back to the axioms (the analog of the “big bang”) then we’d just get a complete graph, connecting all the theorems on slice 26. But here we’re seeing in effect “fuzzier and fuzzier” versions of how the theorems that exist at slice 26 can be thought of as being “metamathematically laid out”. The disconnected components in these branchial graphs represent theorems that have no recent shared history—so that in some sense they’re “causally disconnected”.

In thinking about “theorem search”, it’s interesting to try to imagine measures of “distance between theorems”—and in effect branchial distance captures some of this. And even for Euclid there are presumably things to learn about the “layout” of theorems, and what should count as “close to” what.

There are only 465 theorems in Euclid’s *Elements*. But what if there were many more? What might the “metamathematical space” they define be like? Just as for the hypergraphs—or, for that matter, the multiway graphs—in our models of physics we can ask questions about the limiting emergent geometry of this space. And—ironically enough—one thing we can immediately say is that it seems to be far from Euclidean!

But does it for example have some definite effective dimension? There isn’t enough data to say much about the branchial slices we just saw. But we can say a bit more about the complete theorem dependency graph—which is the analog of the multiway graph in our physics models. For example, starting with the axioms (the analog of the “big bang”) we can ask how many theorems are reached in successive steps. The result (counting the axioms) is:

{10, 81, 325, 444, 470, 475, 475, 475, 475, 475, 475}

If we were dealing with something that approximated a $d$-dimensional manifold, we’d expect these numbers to be of order $r^d$. Computing their logarithmic differences to fit for $d$ gives

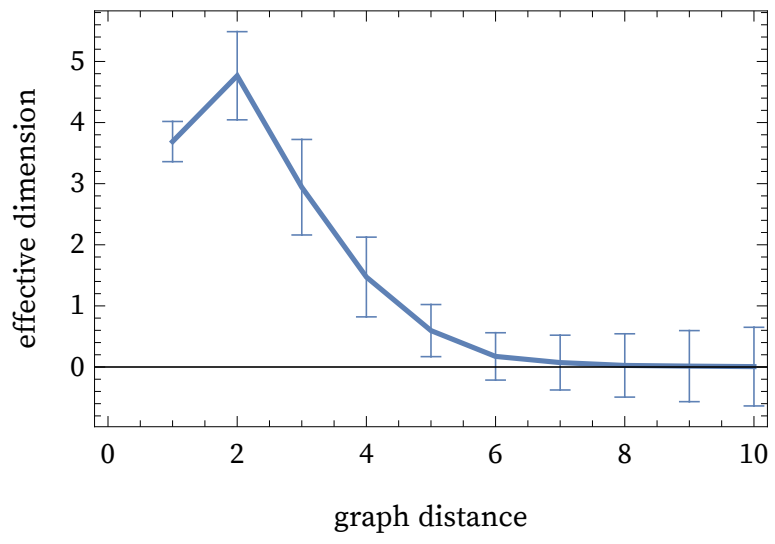

if one starts from the axioms, and

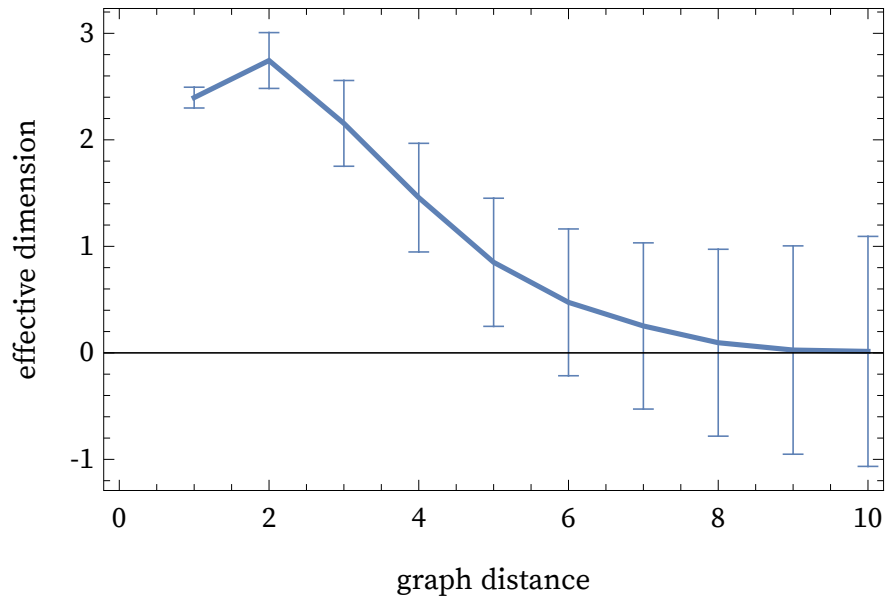

if one starts from all possible theorems in the network.

One gets somewhat different results if one deals not with the actual theorem dependency graph in Euclid, but instead with its transitive reduction—removing all "unnecessary" direct connections. Now the number of theorems reached on successive steps is:

{10, 44, 122, 267, 390, 445, 465, 470, 471, 474, 475}

The "dimension estimate" based on theorems reached starting from the axioms is

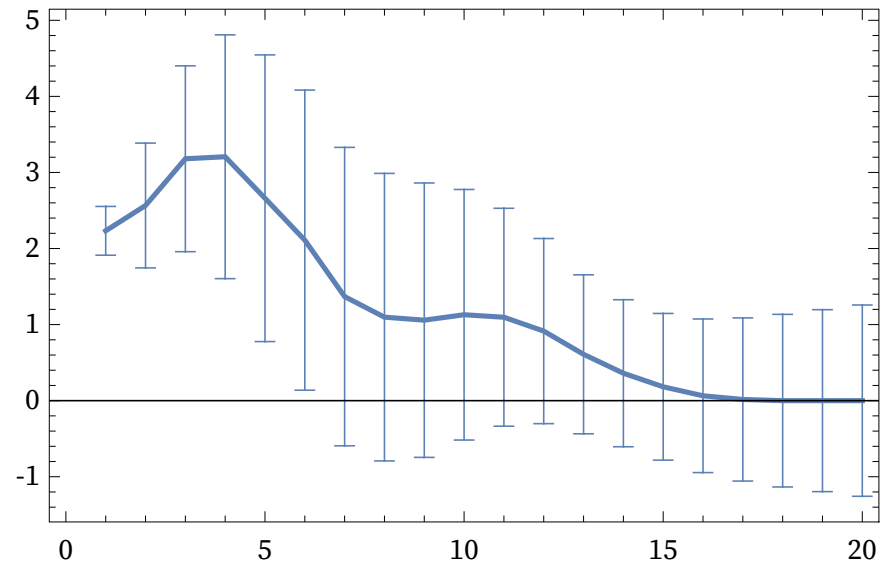

while the corresponding result starting from all theorems is:

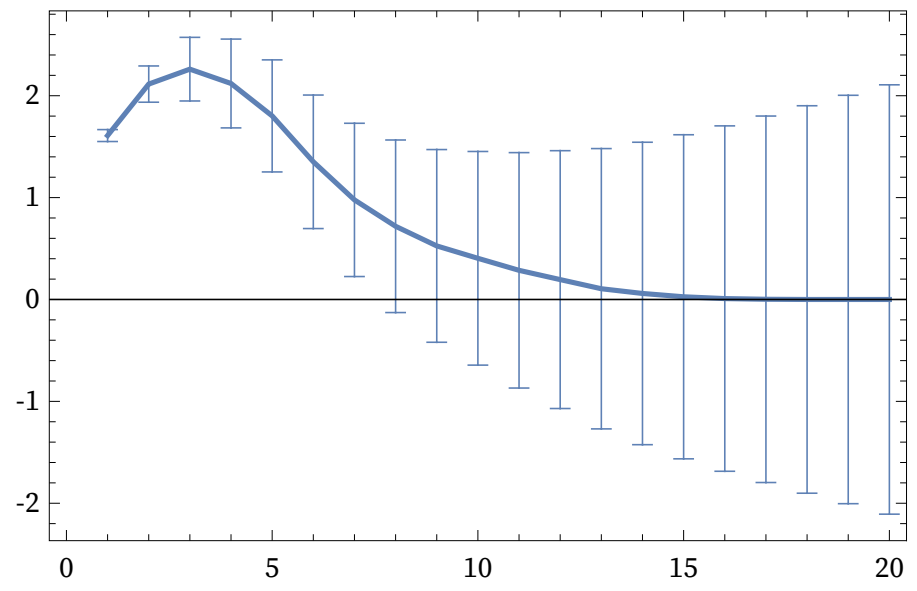

Euclid’s *Elements* represents far too little data to make a definite statement, but perhaps there’s a hint of 2-dimensional structure, with positive curvature.

# The Most Difficult Theorem in Euclid

One way to assess the “difficulty” of a theorem is to look at what results have to have already been built up in order to prove the theorem. And by this measure, the most difficult theorem in Euclid’s *Elements* is the very last theorem in the last book—what one might call “Euclid’s last theorem”, the climax of the *Elements*—Book 13, Theorem 18, which amounts to the statement that there are five Platonic solids, or more specifically:

*τὰς πλευρὰς τῶν πέντε σχημάτων ἐκθέσθαι καὶ συγκρῖναι πρὸς ἀλλήλας.*

This theorem uses all 10 axioms, and 219 of the 464 previous theorems. Here’s its graph of dependencies:

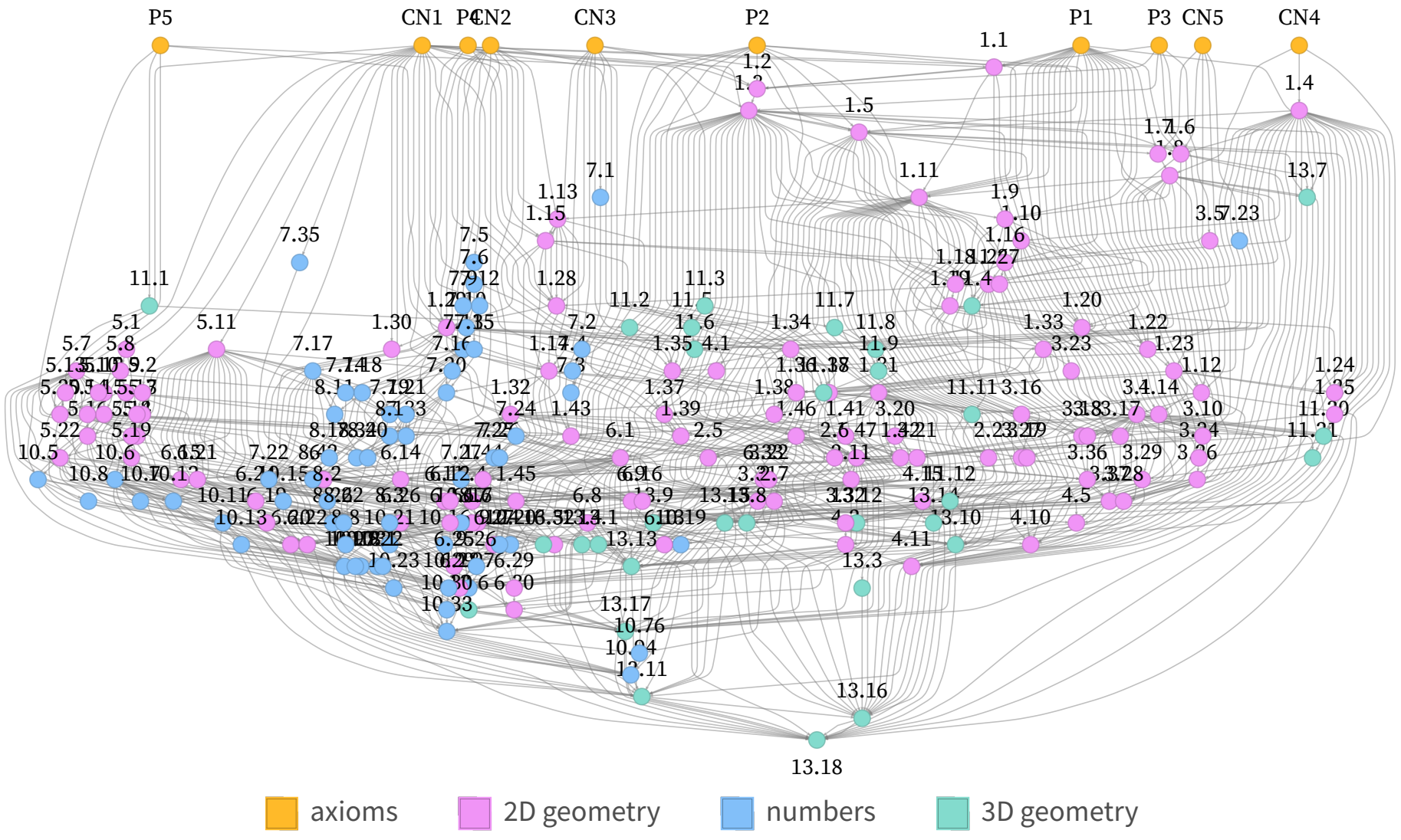

And here is the transitive reduction of this—notably with different subject areas being more obviously separated:

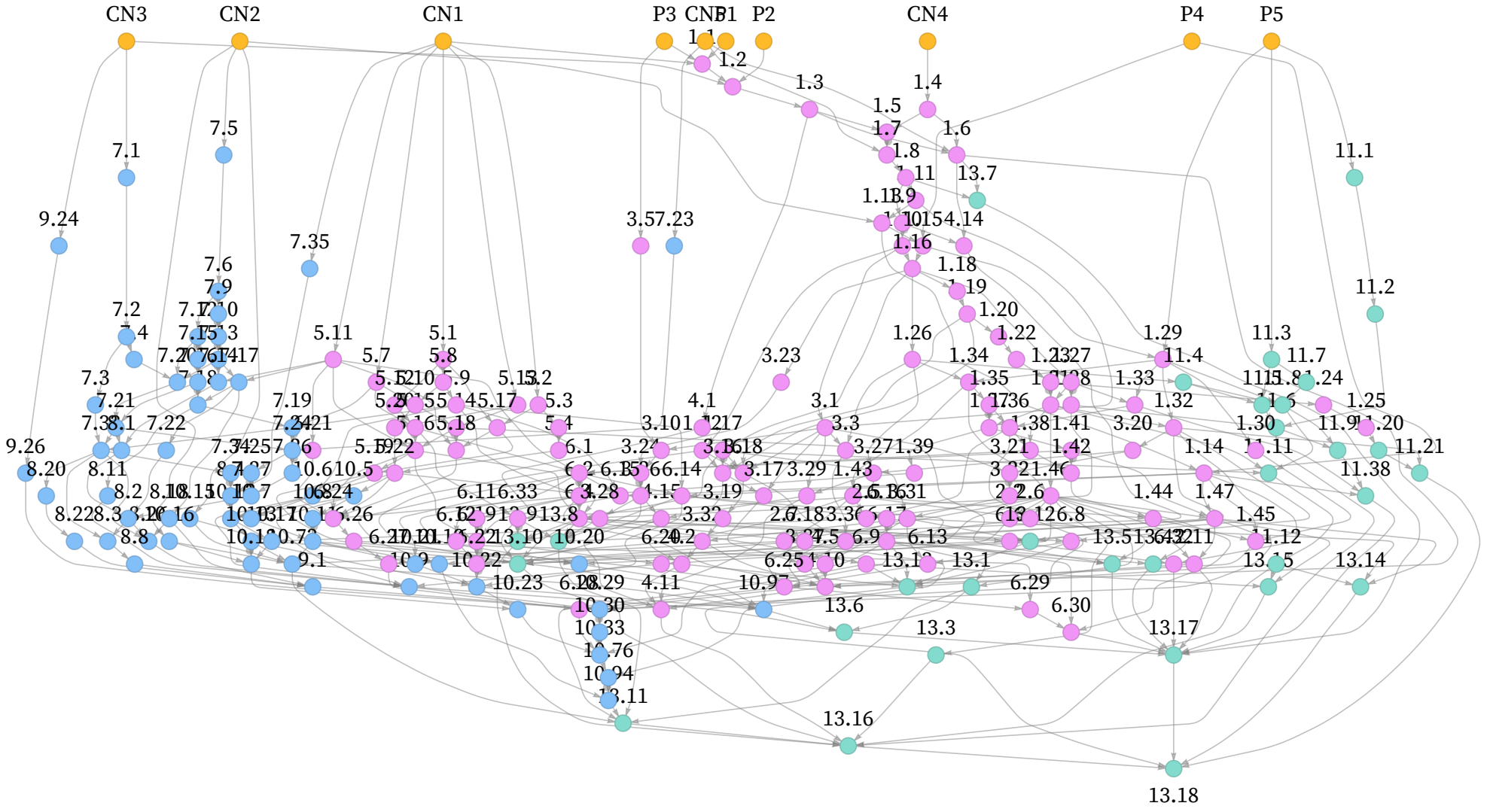

This shows how 13.18 and its prerequisites (its “past light cone”) sit inside the whole theorem dependency graph:

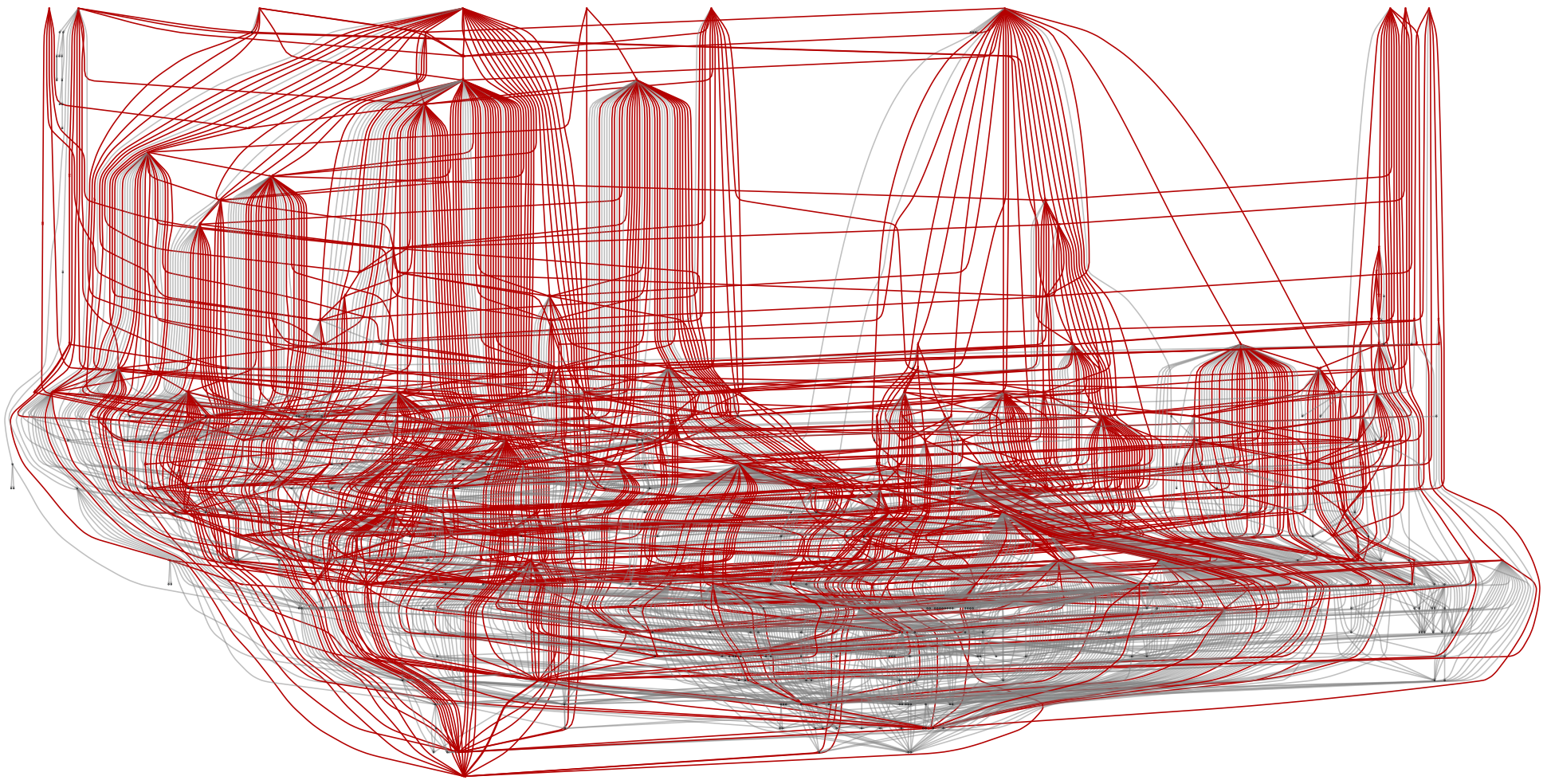

If we started from the axioms, the longest chains of theorems we'd have to prove to get to 13.18 is:

| | |
|---|---|
| CN1 | → 1.1 → 1.2 → 1.3 → 1.5 → 1.7 → 1.8 → 1.11 → 1.13 → 1.15 → 1.16 → 1.18 → 1.19 → 1.20 → 1.22 → 1.23 → 1.31 → 1.37 → 1.41 → 6.1 → 6.2 → 6.11 → 6.19 → 6.20 → 10.9 → 10.29 → 10.30 → 10.33 → 10.76 → 10.94 → 13.11 → 13.16 → 13.18 |
| CN2 | → 1.13 → 1.15 → 1.16 → 1.18 → 1.19 → 1.20 → 1.22 → 1.23 → 1.31 → 1.37 → 1.41 → 6.1 → 6.2 → 6.11 → 6.19 → 6.20 → 10.9 → 10.29 → 10.30 → 10.33 → 10.76 → 10.94 → 13.11 → 13.16 → 13.18 |
| CN3 | → 1.2 → 1.3 → 1.5 → 1.7 → 1.8 → 1.11 → 1.13 → 1.15 → 1.16 → 1.18 → 1.19 → 1.20 → 1.22 → 1.23 → 1.31 → 1.37 → 1.41 → 6.1 → 6.2 → 6.11 → 6.19 → 6.20 → 10.9 → 10.29 → 10.30 → 10.33 → 10.76 → 10.94 → 13.11 → 13.16 → 13.18 |
| CN4 | → 1.4 → 1.5 → 1.7 → 1.8 → 1.11 → 1.13 → 1.15 → 1.16 → 1.18 → 1.19 → 1.20 → 1.22 → 1.23 → 1.31 → 1.37 → 1.41 → 6.1 → 6.2 → 6.11 → 6.19 → 6.20 → 10.9 → 10.29 → 10.30 → 10.33 → 10.76 → 10.94 → 13.11 → 13.16 → 13.18 |
| CN5 | → 1.7 → 1.8 → 1.11 → 1.13 → 1.15 → 1.16 → 1.18 → 1.19 → 1.20 → 1.22 → 1.23 → 1.31 → 1.37 → 1.41 → 6.1 → 6.2 → 6.11 → 6.19 → 6.20 → 10.9 → 10.29 → 10.30 → 10.33 → 10.76 → 10.94 → 13.11 → 13.16 → 13.18 |
| P1 | → 1.1 → 1.2 → 1.3 → 1.5 → 1.7 → 1.8 → 1.11 → 1.13 → 1.15 → 1.16 → 1.18 → 1.19 → 1.20 → 1.22 → 1.23 → 1.31 → 1.37 → 1.41 → 6.1 → 6.2 → 6.11 → 6.19 → 6.20 → 10.9 → 10.29 → 10.30 → 10.33 → 10.76 → 10.94 → 13.11 → 13.16 → 13.18 |
| P2 | → 1.2 → 1.3 → 1.5 → 1.7 → 1.8 → 1.11 → 1.13 → 1.15 → 1.16 → 1.18 → 1.19 → 1.20 → 1.22 → 1.23 → 1.31 → 1.37 → 1.41 → 6.1 → 6.2 → 6.11 → 6.19 → 6.20 → 10.9 → 10.29 → 10.30 → 10.33 → 10.76 → 10.94 → 13.11 → 13.16 → 13.18 |
| P3 | → 1.1 → 1.2 → 1.3 → 1.5 → 1.7 → 1.8 → 1.11 → 1.13 → 1.15 → 1.16 → 1.18 → 1.19 → 1.20 → 1.22 → 1.23 → 1.31 → 1.37 → 1.41 → 6.1 → 6.2 → 6.11 → 6.19 → 6.20 → 10.9 → 10.29 → 10.30 → 10.33 → 10.76 → 10.94 → 13.11 → 13.16 → 13.18 |
| P4 | → 1.15 → 1.16 → 1.18 → 1.19 → 1.20 → 1.22 → 1.23 → 1.31 → 1.37 → 1.41 → 6.1 → 6.2 → 6.11 → 6.19 → 6.20 → 10.9 → 10.29 → 10.30 → 10.33 → 10.76 → 10.94 → 13.11 → 13.16 → 13.18 |
| P5 | → 1.29 → 1.34 → 1.35 → 1.36 → 1.38 → 6.1 → 6.2 → 6.11 → 6.19 → 6.20 → 10.9 → 10.29 → 10.30 → 10.33 → 10.76 → 10.94 → 13.11 → 13.16 → 13.18 |

Or in other words, from CN1 and from P1 and P3 we'd have to go 33 steps to reach 13.18. If we actually look at the paths, however, we see that after different segments at the beginning, they all merge at Book 6, Theorem 1, and then are the same for the last 14 steps:

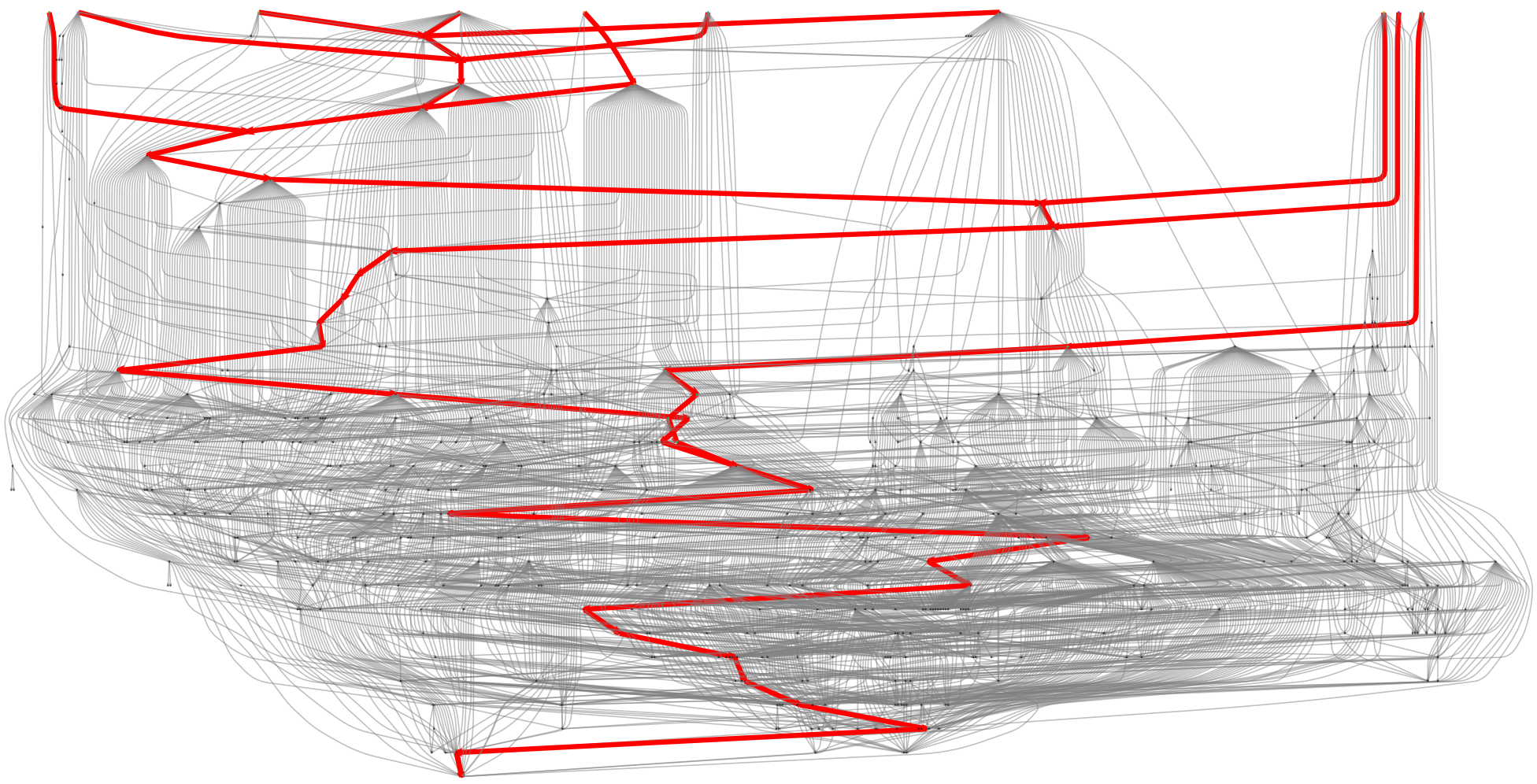

(Theorem 6.1 is the statement that both triangles and parallelograms that have the same base and same height have the same area, i.e. one can skew a triangle or parallelogram without changing its area.)

How much more difficult than other theorems is 13.18? Here's a histogram of maximum path lengths for all theorems (ignoring cases to be discussed later where a particular theorem does not use a given axiom at all):

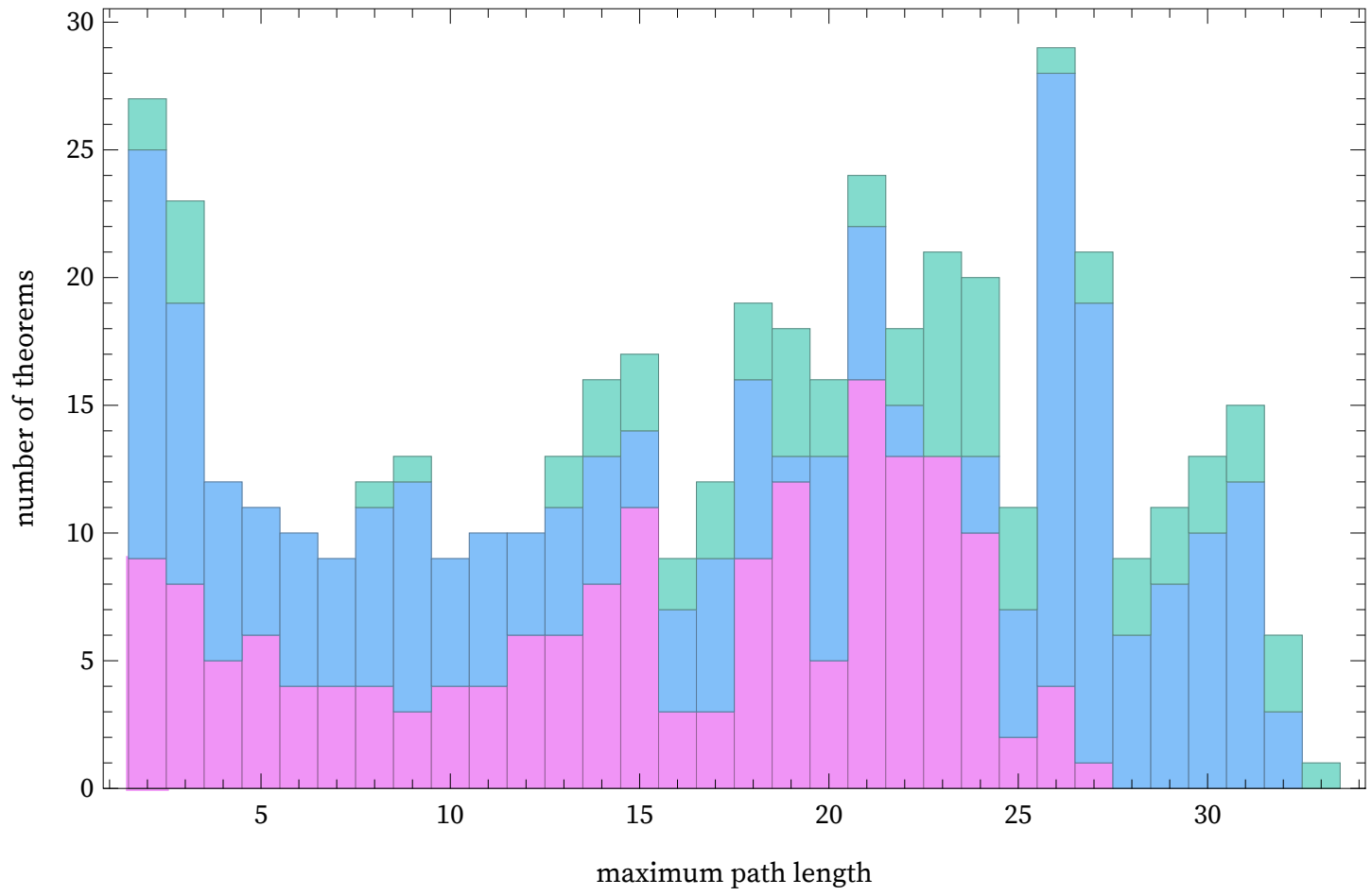

And here's how the maximum path length varies through the sequence of all 465 theorems:

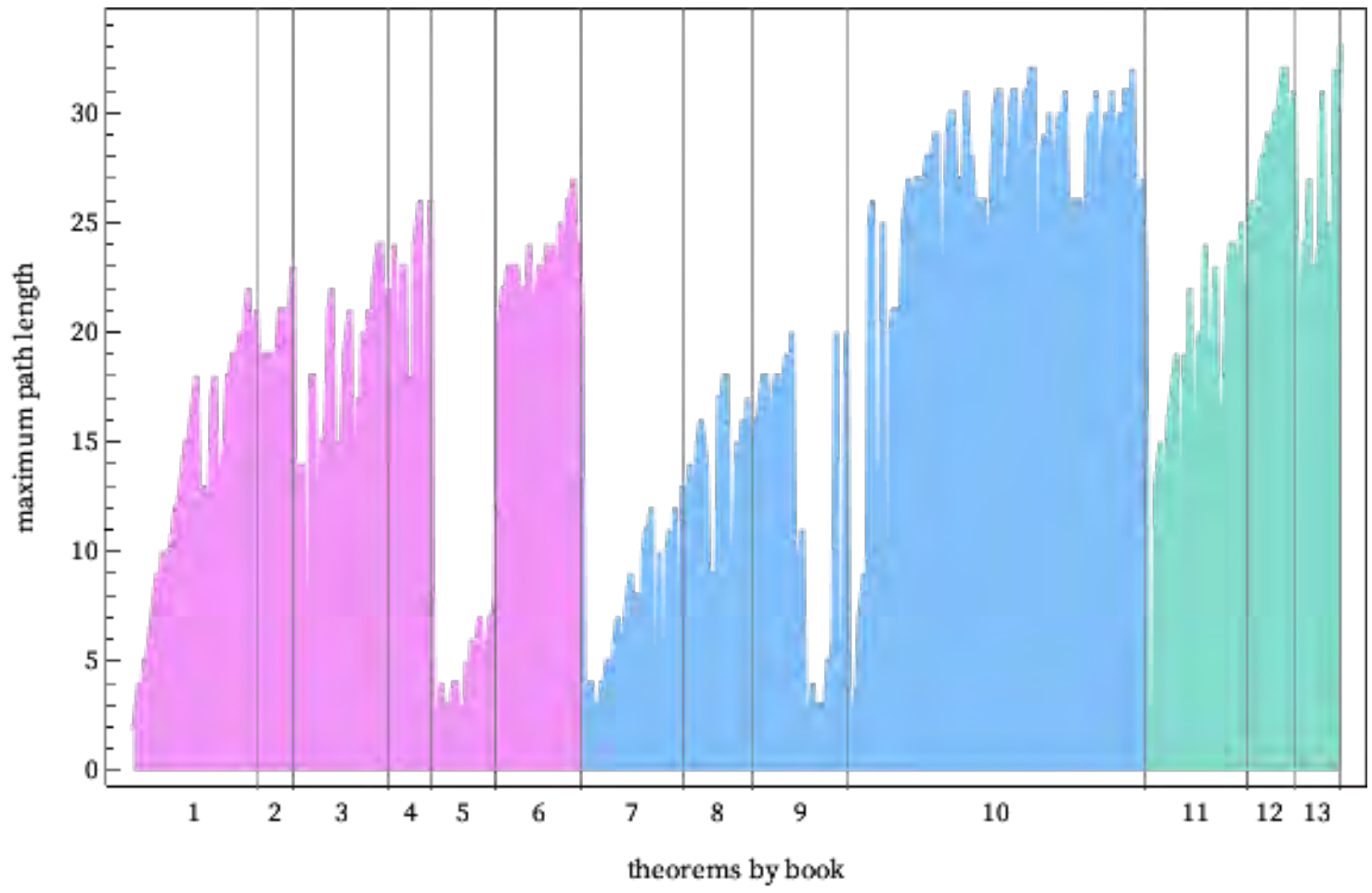

In the causal graph interpretation, and using the "flat foliation" (i.e. the "cosmological rest frame") what this basically shows is at what "time slice" a given theorem first emerges from Euclid's proofs. Or, in other words, if one imagines exploring the "metamathematical space of Euclid" by going "one level of theorems at a time", the order in which one will encounter theorems is:

| |
|---|
| CN1 CN2 CN3 CN4 CN5 P1 P2 P3 P4 P5 |
| 1.1 1.4 3.5 3.6 5.1 5.2 5.7 5.11 5.13 7.1 7.5 7.23 7.28 7.29 7.31 7.35 9.7 9.21 9.24 9.33 9.34 10.1 10.2 10.15 10.16 11.1 11.3 |
| 1.2 5.3 5.5 5.6 5.8 5.12 5.17 6.21 7.2 7.6 7.7 7.32 9.20 9.22 9.25 9.26 9.27 9.28 10.3 11.2 11.7 11.13 11.16 |
| 1.3 5.4 5.9 5.10 5.15 7.3 7.4 7.8 7.9 7.12 9.23 10.4 |
| 1.5 1.6 4.1 5.14 5.20 5.21 7.10 7.11 7.15 9.29 9.30 |
| 1.7 5.16 5.18 5.22 7.13 7.16 7.37 7.38 9.11 9.31 |
| 1.8 5.19 5.23 5.24 7.14 7.17 7.20 9.35 10.5 |
| 1.9 1.11 4.9 5.25 7.18 7.21 7.22 8.13 8.18 10.6 10.8 13.7 |
| 1.10 1.13 4.14 7.19 8.1 8.10 8.11 8.12 8.19 9.16 9.17 10.7 11.19 |
| 1.12 1.14 1.15 3.1 7.24 7.30 7.33 9.18 10.11 |
| 1.16 1.29 3.9 3.10 7.25 7.26 7.34 8.20 9.14 9.19 |
| 1.17 1.18 1.26 1.27 1.30 3.23 7.27 7.36 8.4 8.22 |
| 1.19 1.28 1.33 1.34 3.3 3.24 7.39 8.2 8.5 8.24 10.12 11.4 11.14 |
| 1.20 1.35 1.43 3.2 3.4 3.16 3.18 3.26 8.3 8.9 8.26 8.27 10.13 11.5 11.8 11.28 |
| 1.21 1.22 1.36 3.11 3.12 3.17 3.19 3.28 4.4 4.7 4.13 8.6 8.8 8.21 11.6 11.18 12.16 |
| 1.23 3.13 3.30 8.7 8.23 9.1 9.2 11.9 11.29 |
| 1.24 1.31 3.25 8.14 8.15 8.25 9.3 9.6 9.8 11.10 11.30 11.38 |
| 1.25 1.32 1.37 1.38 1.46 2.1 3.7 3.8 4.8 8.16 8.17 9.4 9.5 9.9 9.10 9.12 11.11 11.22 11.24 |
| 1.39 1.40 1.41 2.2 2.3 2.4 2.5 2.6 2.7 2.8 3.20 4.3 9.13 11.12 11.15 11.20 11.25 13.3 |
| 1.42 1.47 3.21 3.27 6.1 9.15 9.32 9.36 10.17 10.36 10.73 10.74 10.77 11.21 11.26 11.31 |
| 1.44 1.48 2.9 2.10 2.11 2.12 2.13 3.14 3.22 3.29 3.35 3.36 6.2 6.14 6.15 6.33 10.18 10.19 10.20 10.21 10.37 10.106 11.39 13.2 |
| 1.45 3.15 3.31 3.37 4.15 6.3 6.4 6.9 6.10 6.11 6.12 6.16 6.24 10.55 10.115 11.17 11.35 12.3 |
| 2.14 3.32 4.5 4.6 6.5 6.6 6.7 6.8 6.17 6.18 6.19 6.23 6.26 11.27 11.32 11.36 13.8 13.9 13.12 13.14 13.15 |
| 3.33 3.34 4.2 4.10 6.13 6.20 6.22 6.27 6.31 6.32 10.54 10.92 10.95 11.23 11.33 11.34 13.1 13.4 13.5 13.10 |
| 4.11 6.25 10.9 10.14 10.22 10.91 10.109 11.37 12.1 12.4 13.13 |
| 4.12 4.16 6.28 6.29 10.10 10.23 10.25 10.26 10.29 10.38 10.49 10.50 10.51 10.52 10.53 10.60 10.66 10.85 10.86 10.87 10.88 10.89 10.90 10.97 10.101 10.103 10.112 10.113 12.2 |
| 6.30 10.24 10.27 10.28 10.30 10.42 10.43 10.44 10.45 10.56 10.61 10.62 10.67 10.79 10.80 10.83 10.98 10.111 10.114 12.5 13.6 |
| 10.31 10.32 10.33 10.47 10.48 10.75 12.6 12.7 13.17 |
| 10.34 10.35 10.39 10.76 10.81 10.93 10.99 10.104 12.8 12.9 12.10 |
| 10.40 10.41 10.57 10.63 10.68 10.78 10.82 10.94 10.100 10.105 12.11 12.12 12.17 |
| 10.46 10.58 10.59 10.64 10.65 10.69 10.70 10.84 10.96 10.102 10.107 10.108 12.13 12.18 13.11 |
| 10.71 10.72 10.110 12.14 12.15 13.16 |
| 13.18 |

A question one might ask is whether “short-to-state” theorems are somehow “easier to prove” than longer-to-state ones. This shows the maximum path length to prove theorems as a function of the length of their statement in Euclid’s Greek. Remarkably little correlation is seen.

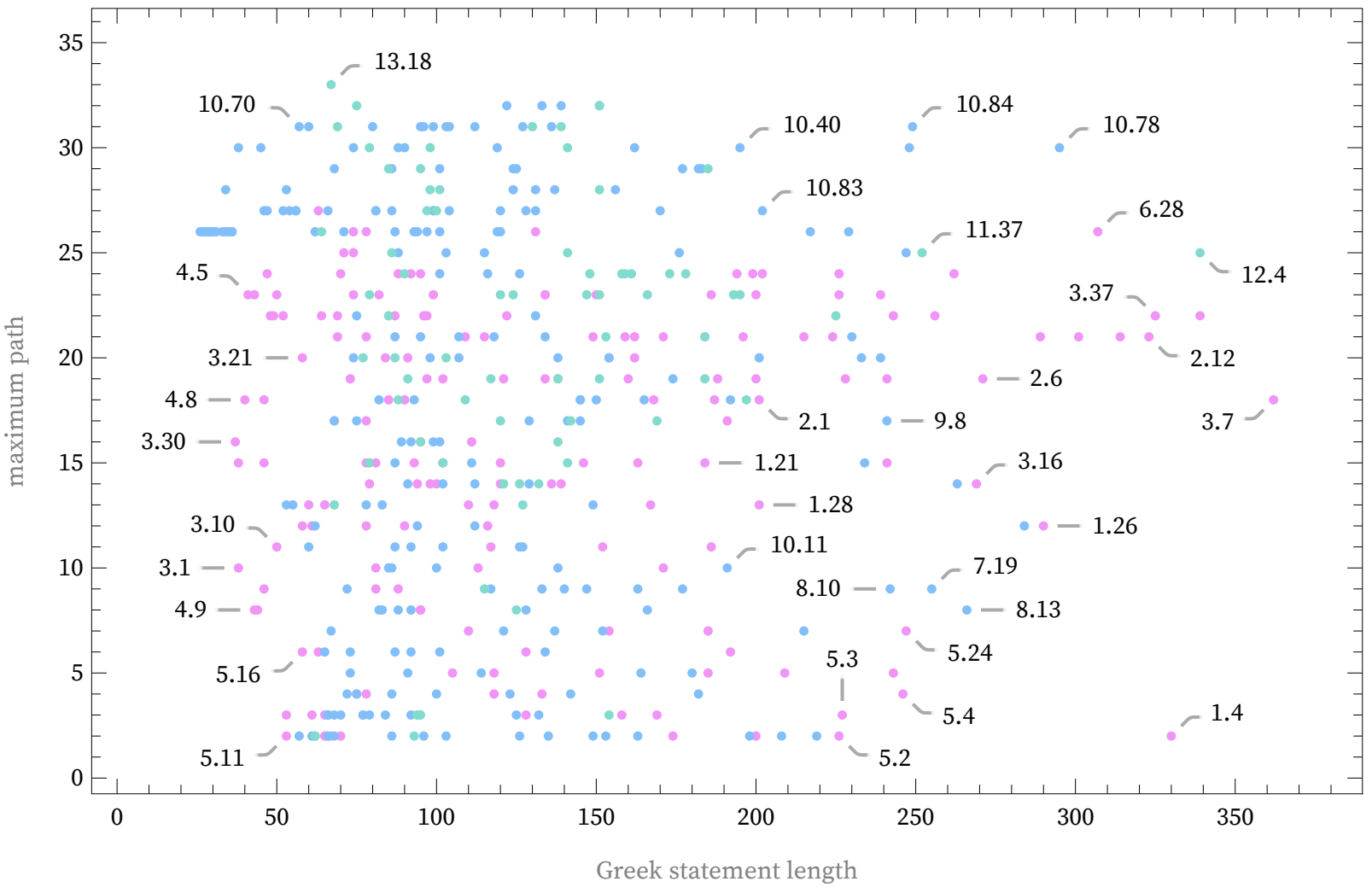

This plot shows instead the number of “prerequisite theorems” as a function of statement length:

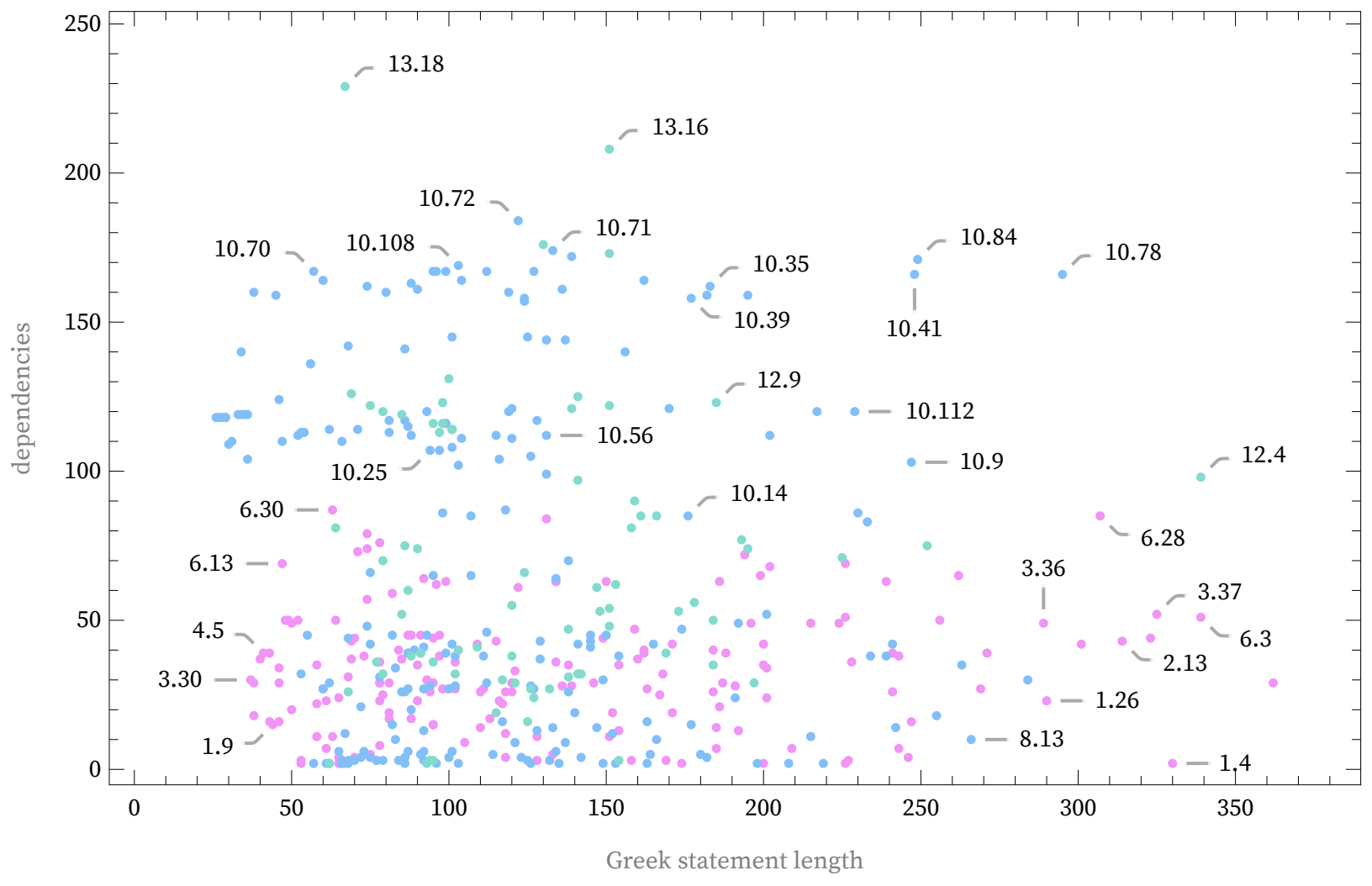

Once again there is poor correlation.

# The Most Popular Theorems in Euclid

How often do particular theorems get used in the proofs of other theorems? The "most popular" theorems in terms of being directly quoted in the proofs of other theorems are:

| | | | | | | | | | |
|---|---|---|---|---|---|---|---|---|---|
| 10.11 | *59* | CN1 | *29* | 5.19 | *22* | 5.22 | *16* | 5.9 | *13* |
| 6.1 | *52* | 10.73 | *27* | 10.15 | *21* | 3.31 | *16* | CN3 | *13* |
| 5.11 | *46* | 10.22 | *26* | 10.36 | *21* | 10.9 | *16* | 10.18 | *13* |
| 1.3 | *46* | 5.7 | *26* | 1.32 | *21* | 2.7 | *16* | 7.20 | *12* |
| 10.6 | *42* | P1 | *25* | 1.46 | *20* | 1.12 | *16* | P2 | *11* |
| 1.4 | *41* | 1.34 | *23* | 1.5 | *20* | 10.20 | *16* | 2.4 | *11* |
| 1.31 | *39* | 1.47 | *23* | 1.29 | *18* | 7.17 | *15* | CN2 | *11* |
| 10.13 | *32* | 5.16 | *23* | 6.17 | *18* | 6.4 | *14* | 3.16 | *11* |
| 1.11 | *32* | 10.21 | *22* | 10.23 | *18* | 1.10 | *14* | 10.17 | *10* |
| 3.1 | *29* | 1.8 | *22* | 1.23 | *18* | 10.12 | *13* | 6.2 | *10* |

Notably, all but one of 10.11's direct mentions are in other theorems in Book 10. Theorem 6.1 (which we already encountered above) is used in 4 books.

By the way, there is some subtlety here, because 26 theorems reference a particular theorem more than once in their proofs: for example, 10.4 references 10.3 three times, while 13.18 references both 13.17 and 13.16 twice:

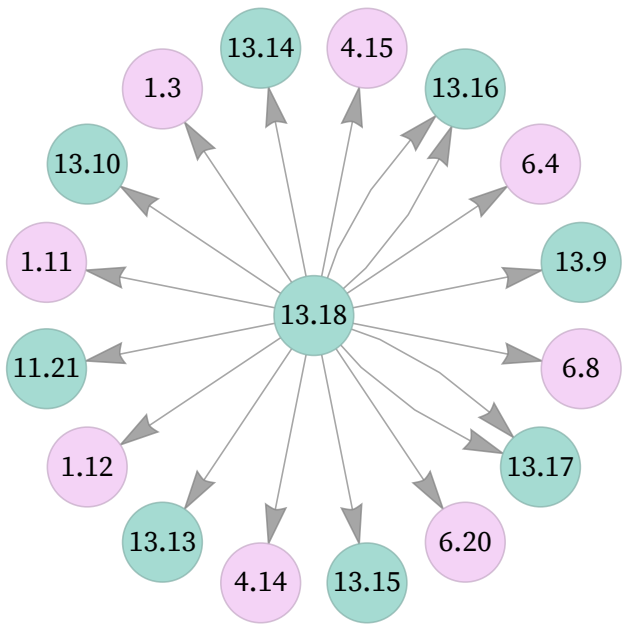

But looking simply at the distribution of the number of direct uses (here on a log scale), we see that the vast majority of theorems are very rarely used—with just a few being quite widely used:

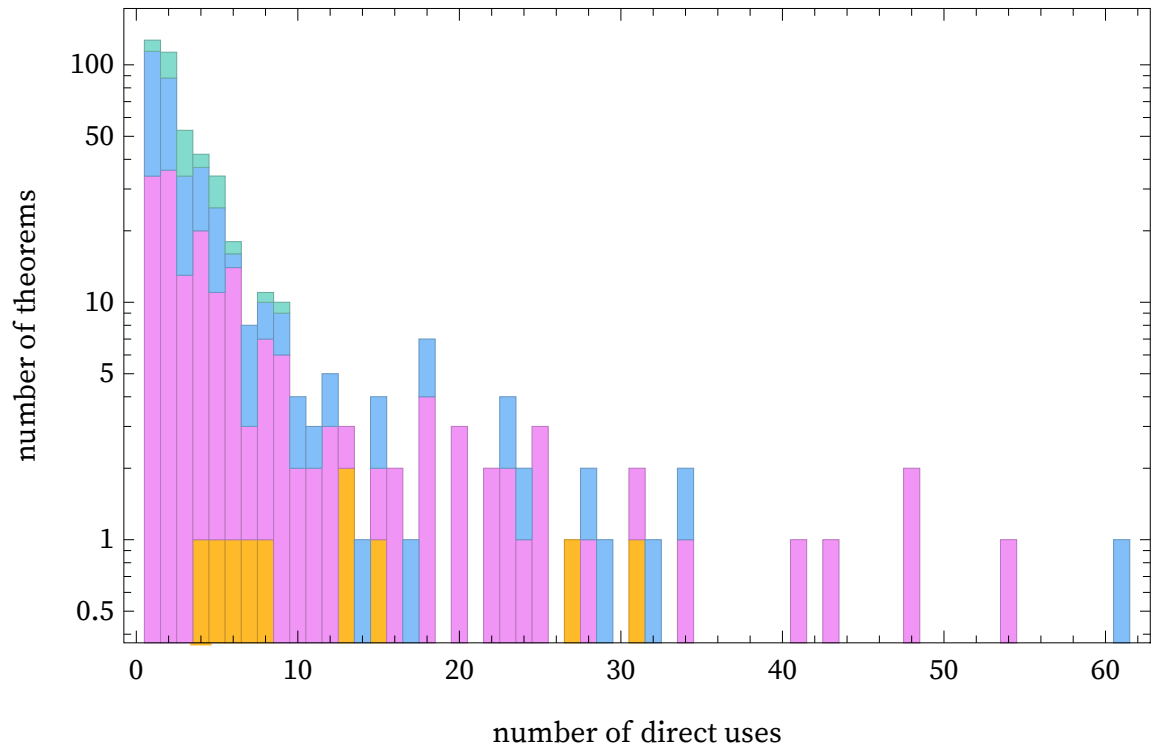

Indicating the number of direct uses by size, here are the “directly popular” theorems:

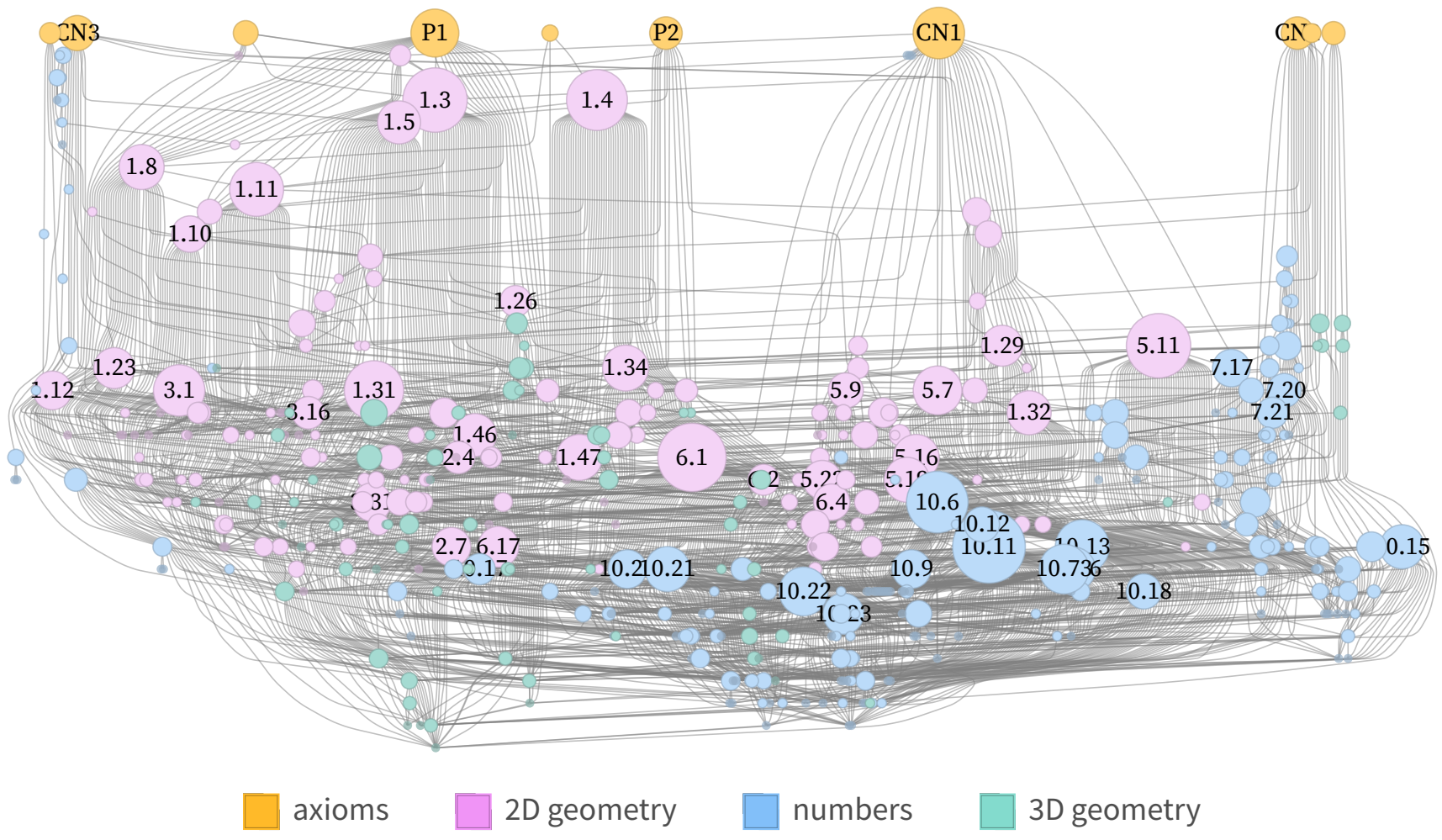

If we ask also about indirect uses, the results are as follows:

| | | | | | | | | | |
|---|---|---|---|---|---|---|---|---|---|
| CN1 | *411* | 1.3 | *313* | 1.15 | *294* | 1.23 | *244* | 1.37 | *189* |
| CN3 | *390* | 1.4 | *312* | 1.16 | *291* | 5.1 | *243* | 1.41 | *187* |
| CN2 | *377* | 1.5 | *310* | P5 | *262* | 5.11 | *241* | 1.36 | *187* |
| CN5 | *344* | 1.7 | *309* | 1.18 | *261* | 5.7 | *238* | 5.12 | *179* |
| P3 | *318* | 1.8 | *308* | 1.19 | *260* | 1.31 | *236* | 5.15 | *178* |
| P1 | *316* | 1.11 | *301* | 1.29 | *251* | 5.8 | *235* | 7.5 | *177* |
| P2 | *315* | 1.9 | *298* | 1.20 | *250* | 1.34 | *222* | 7.6 | *173* |
| 1.1 | *315* | P4 | *296* | 1.27 | *248* | 5.9 | *221* | 1.38 | *169* |
| 1.2 | *314* | 1.13 | *296* | 1.22 | *245* | 1.35 | *202* | 5.13 | *166* |
| CN4 | *313* | 1.10 | *296* | 1.26 | *244* | 1.33 | *192* | 5.10 | *166* |

Not too surprisingly, the axioms and early theorems are the most popular. But overall, the distribution of total number of uses is somewhat broader than the distribution of direct uses:

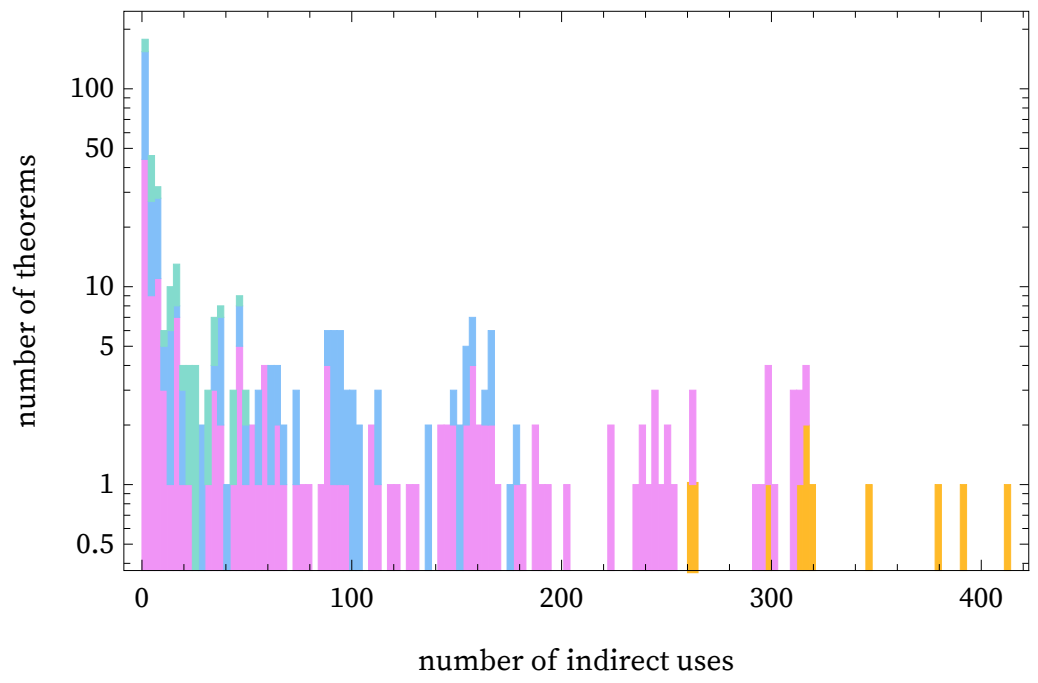

This shows all theorems, with their sizes in the graph essentially determined by the sizes of their "future light cone" in the theorem dependency graph:

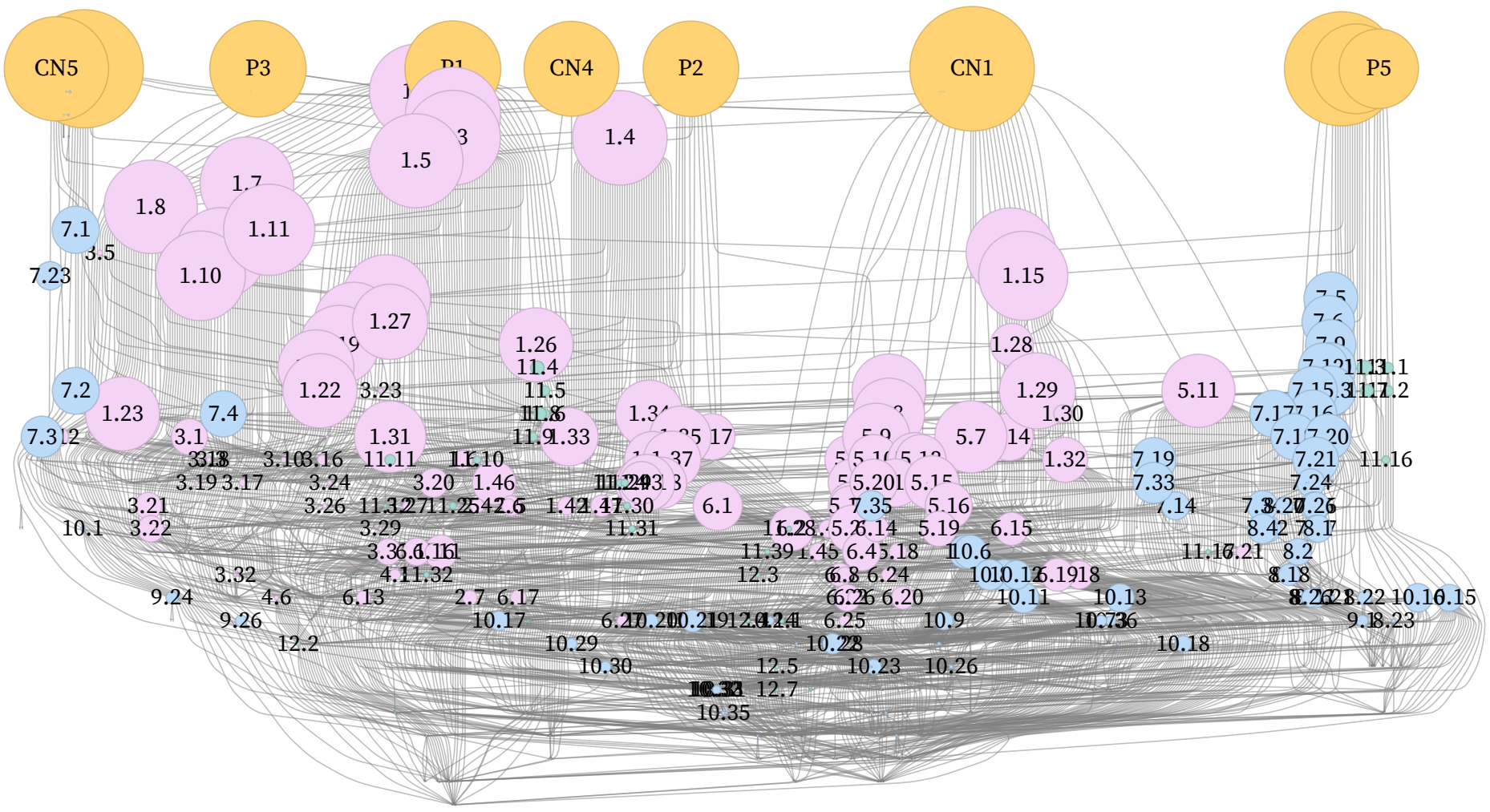

In addition to asking about direct and indirect uses, one can also assess the "centrality" of a given theorem by various graph-theoretical measures. One example is betweenness centrality (the fraction of shortest paths that pass through a given node):

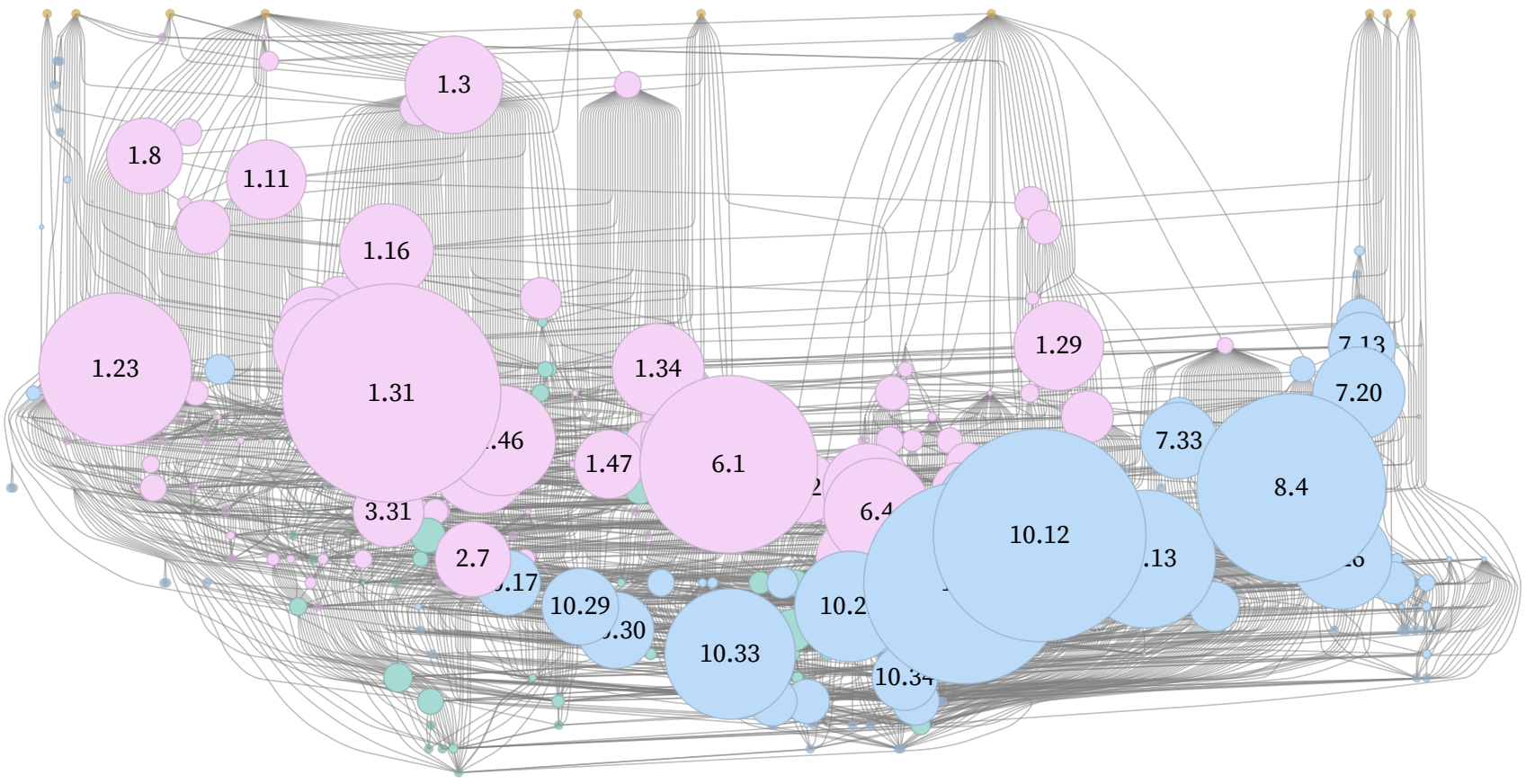

The theorems with top betweenness centralities are 1.31 (construction of parallel lines), 10.12 (transitivity of commensurability), 10.9 (commensurabilty in squares), 8.4 (continued ratios in lowest terms), etc.

For closeness centrality (average inverse distance to all other nodes) one gets:

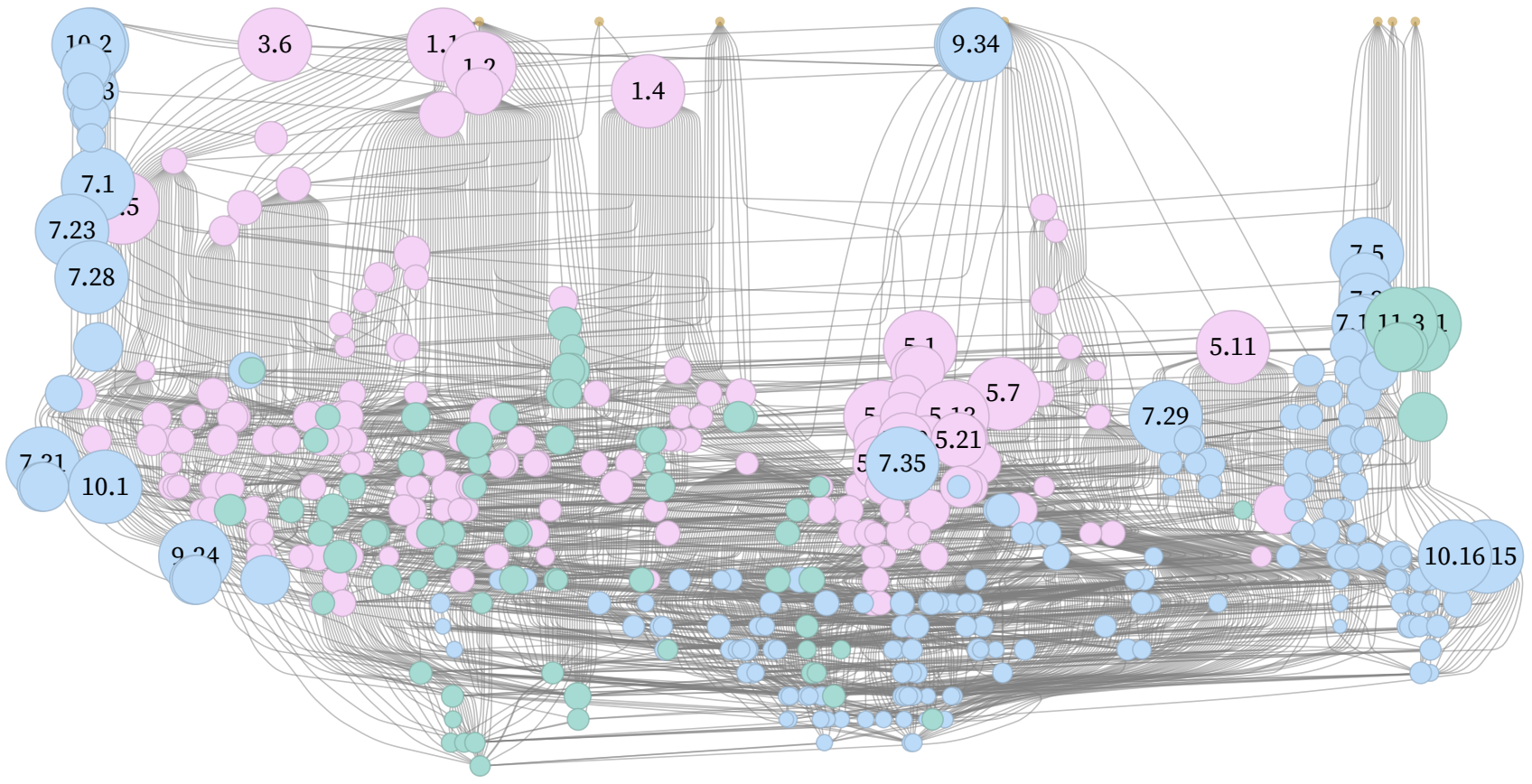

## What Really Depends on What?

Euclid's *Elements* starts with 10 axioms, from which all the theorems it contains are derived. But what theorems really depend on what axioms? This shows how many of the 465 theorems depend on each of the Common Notions and Postulates according to the proofs given in Euclid:

| CN1 | CN2 | CN3 | CN4 | CN5 | P1 | P2 | P3 | P4 | P5 |
|---|---|---|---|---|---|---|---|---|---|
| 411 | 377 | 390 | 313 | 344 | 316 | 315 | 318 | 296 | 262 |

The famous fifth postulate (that parallel lines do not cross) has the fewest theorems depending on it. (And actually, for many centuries there was a suspicion that no theorems really depended on it—so people tried to find proofs that didn't use it, although ultimately it became clear it actually was needed.)

Interestingly, at least according to Euclid, more than half (255 out of 465) of the theorems actually depend on all 10 axioms, though one sees definite variation through the course of the *Elements*:

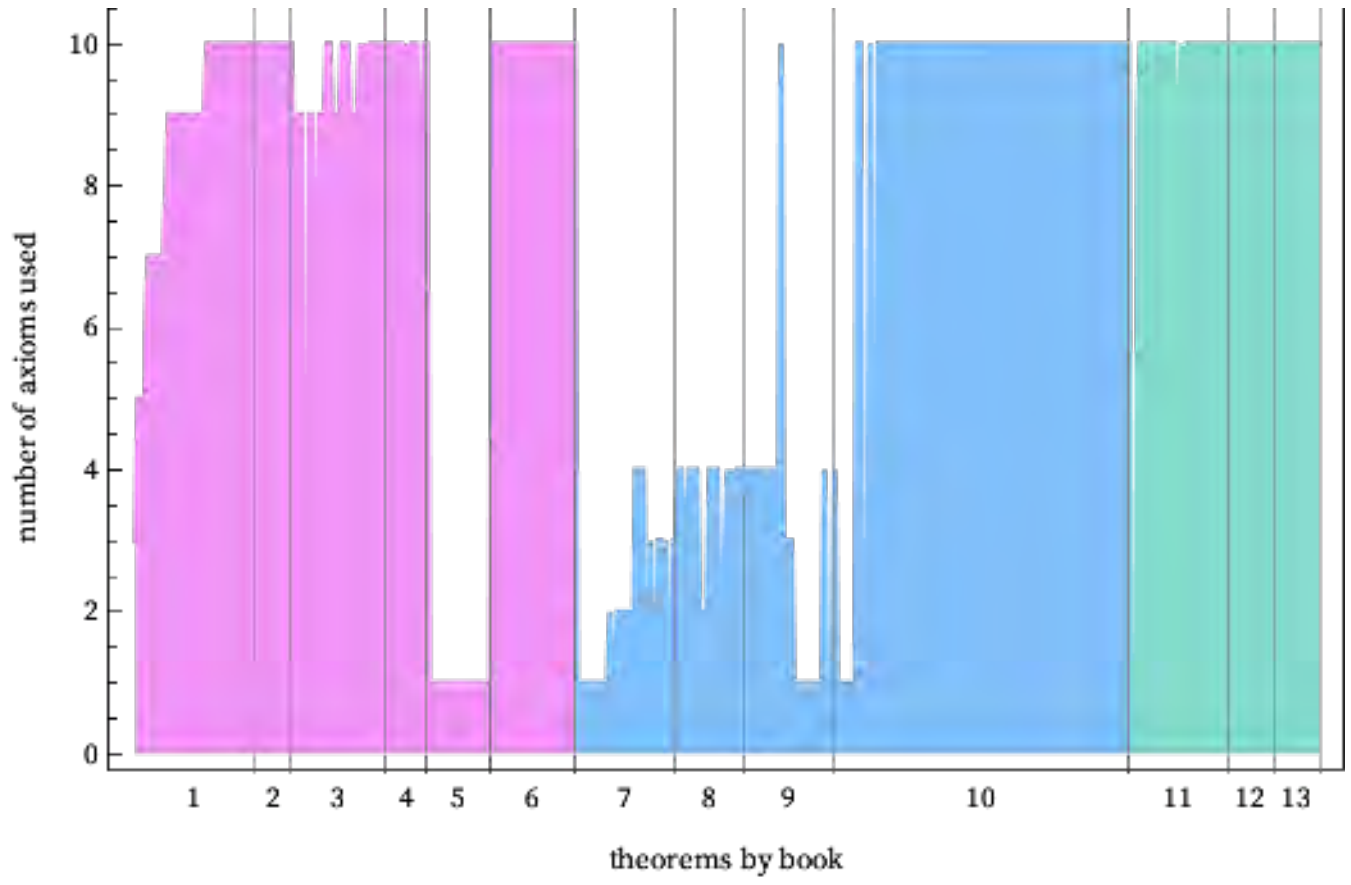

The number of theorems depending on different numbers of axioms is:

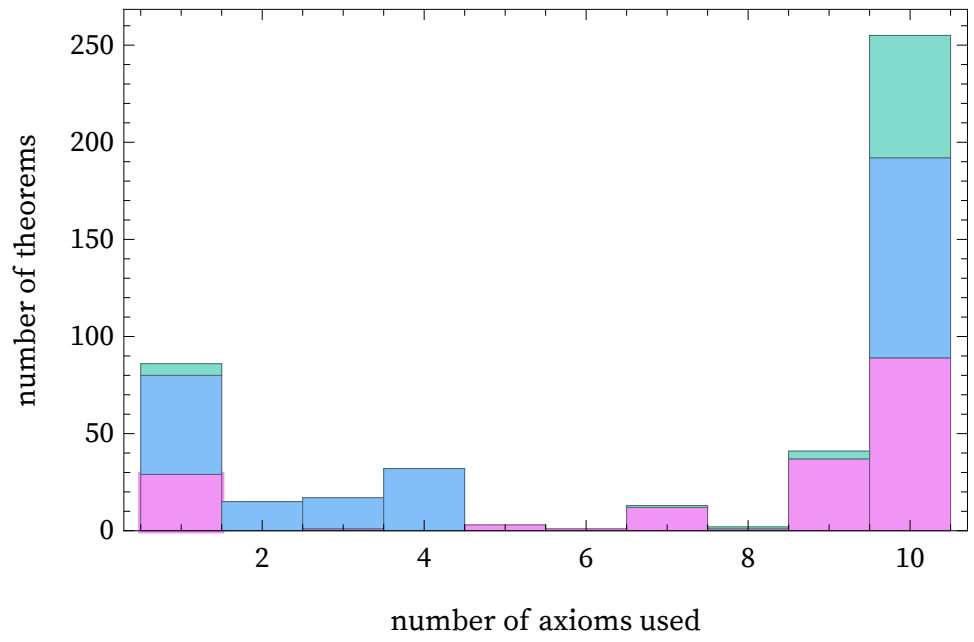

Scattered through the *Elements* there are 86 theorems that depend only on one axiom, most often CN1 (which is transitivity of equality):

| CN1 | CN3 | CN2 | P5 | P3 | CN5 | CN4 |
|---|---|---|---|---|---|---|
| 36 | 23 | 17 | 6 | 2 | 1 | 1 |

In most cases, the dependence is quite direct, but there are cases in which it is actually quite elaborate, such as:

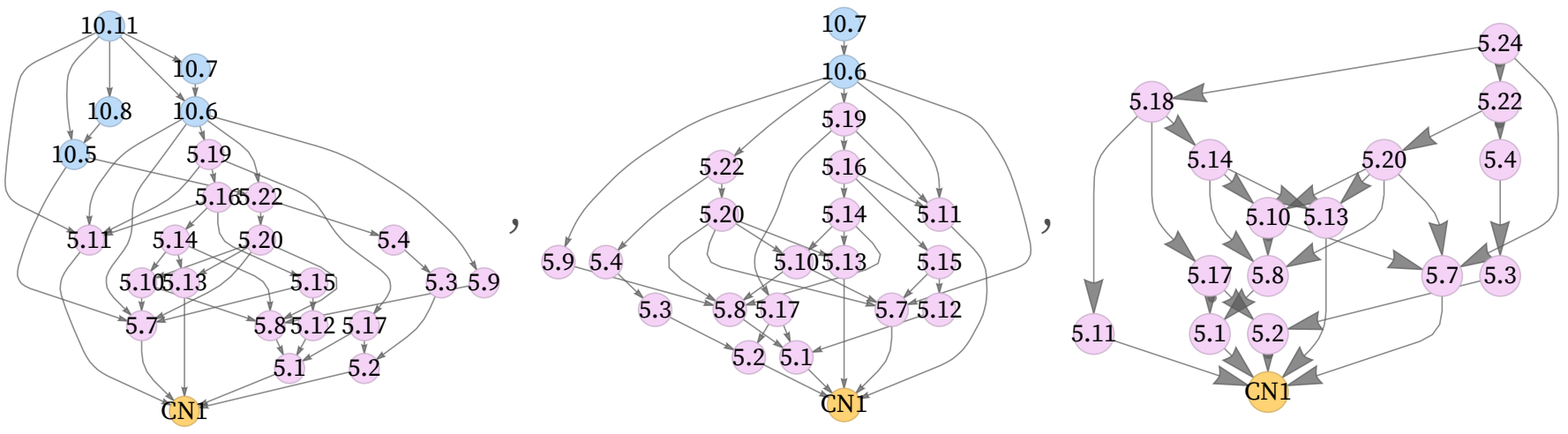

These get slightly simpler after transitive reduction:

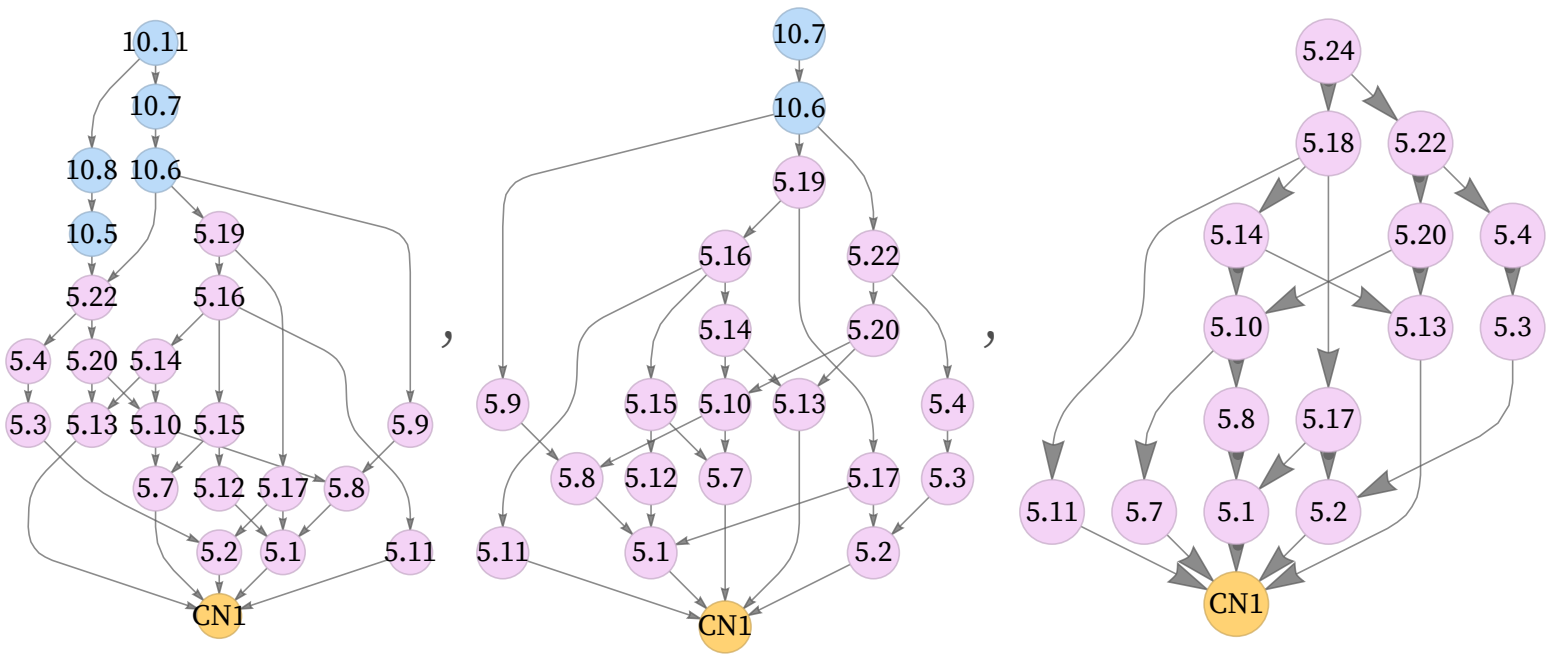

We can now also ask the opposite question of how many theorems don't depend on any given axiom (and, yes, this immediately follows from what we listed above):

| CN1 | CN2 | CN3 | CN4 | CN5 | P1 | P2 | P3 | P4 | P5 |
|---|---|---|---|---|---|---|---|---|---|
| 63 | 97 | 84 | 161 | 130 | 158 | 159 | 156 | 178 | 212 |

And in general we can ask what subsets of the axioms different theorems depend on. Interestingly, of the 1024 possible such subsets, only 19 actually occur, suggesting some considerable correlation between the axioms. Here is a representation of the partial ordering of the subsets that occur, indicating in each case for how many theorems that subset of dependencies occurs:

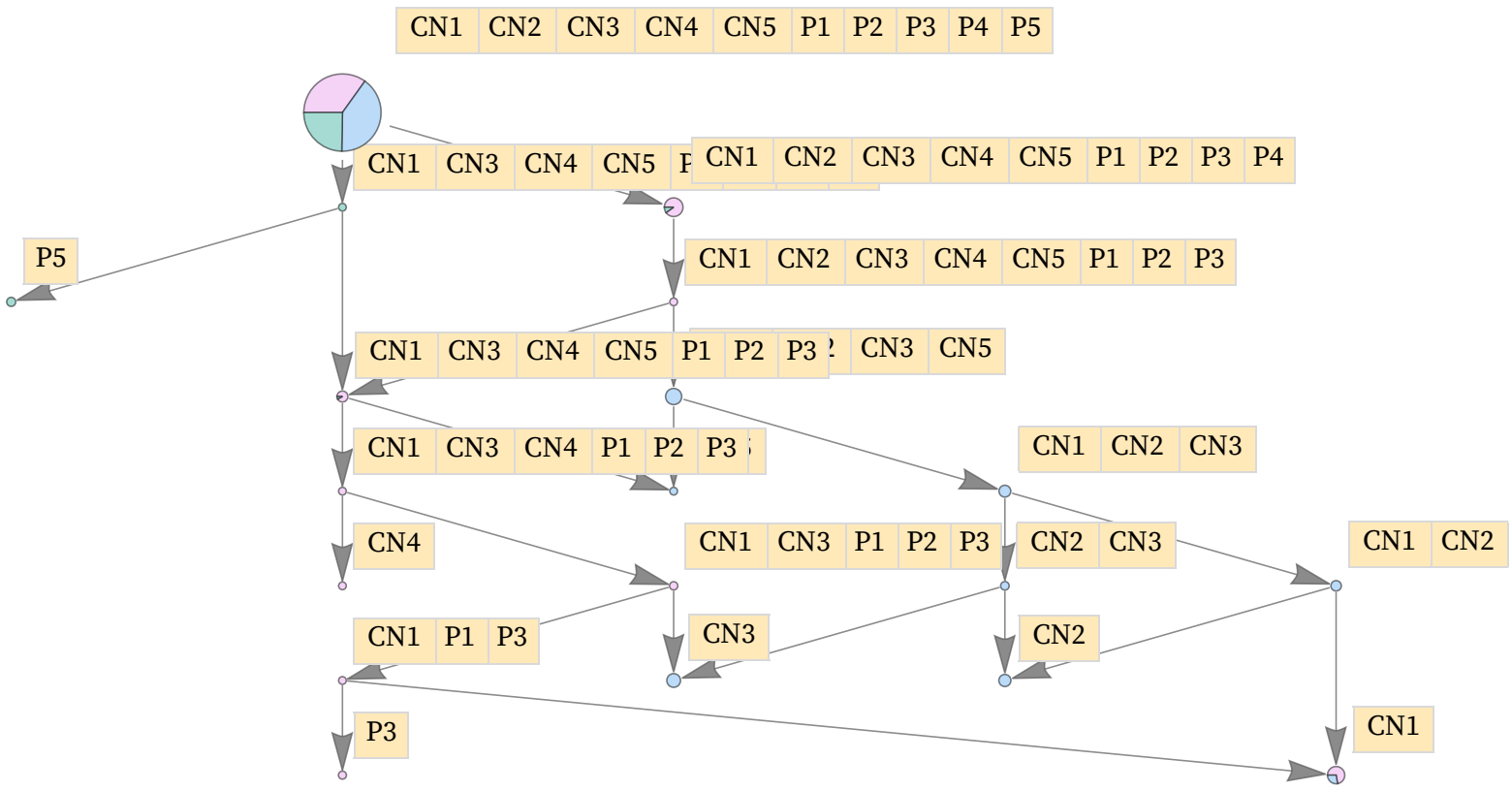

# The Machine Code of Euclid: All the Way Down to Axioms

Any theorem in Euclid can ultimately be proved just by using Euclid's axioms enough times. In other words, the proofs Euclid gave were stated in terms of "intermediate theorems"—but we can always in principle just "compile things down" so we just get a sequence of axioms. And here for example is how that works for Book 1, Theorem 5:

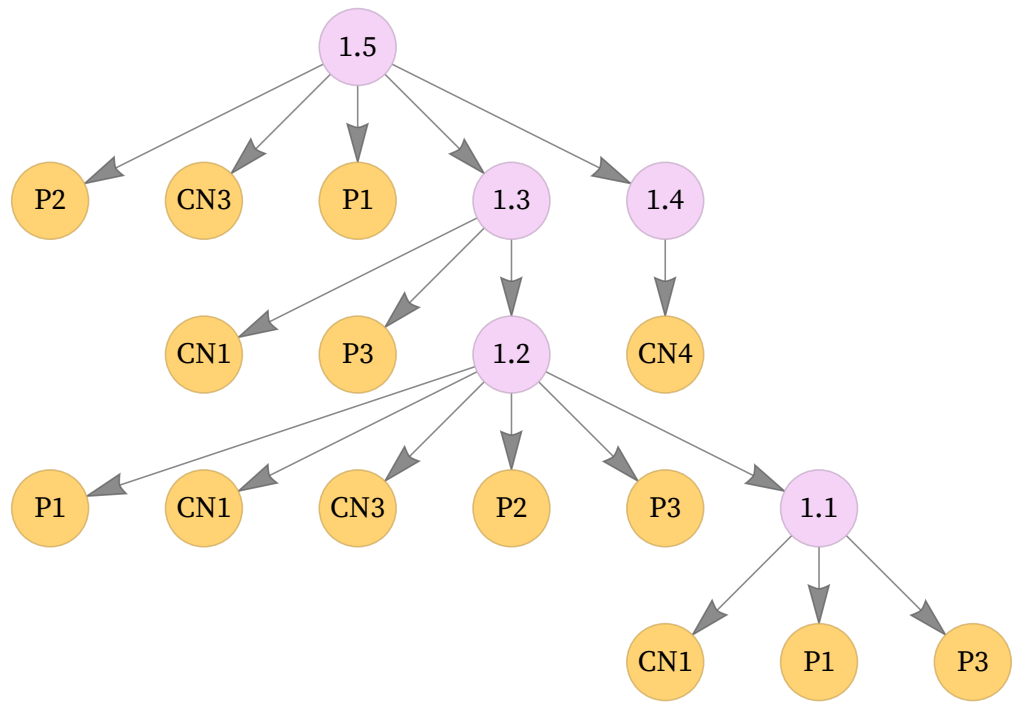

Of course it's much more efficient to "share the work" by using intermediate theorems:

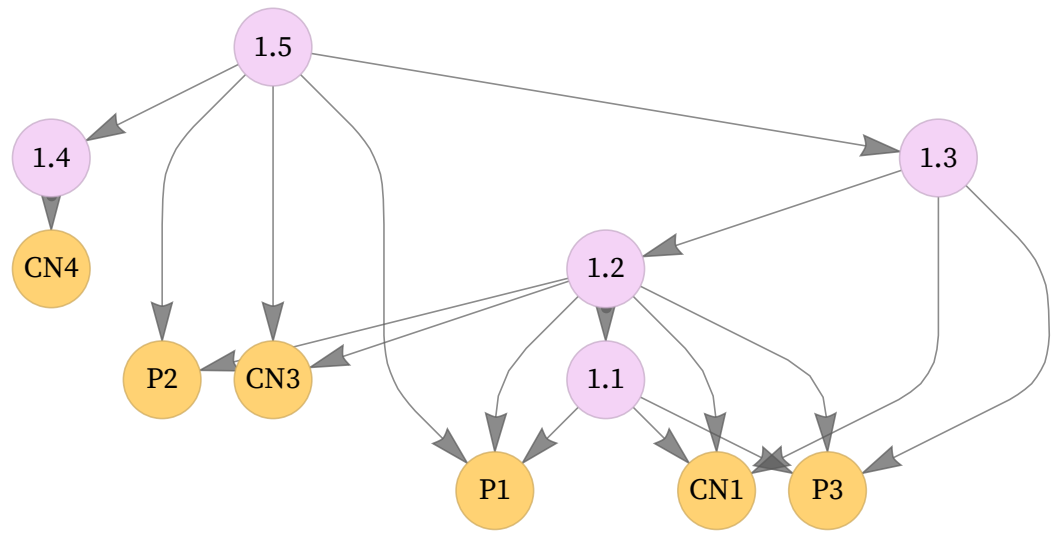

This doesn't change the "depth"—i.e. the length of any given path to get to the axioms. But it reduces the number of independent paths that have to be followed, because every time one reaches the same theorem (or axiom) one just "uses what one already knows about it".

But to get a sense of the "axiomatic machine code" of Euclid we can just "compile" the proof of every theorem down to its underlying sequence of axioms. And for example if we do this for 3.18 the final sequence of axioms we get has length 835,416. These are broken down among the various axioms according to:

| CN1 | CN2 | CN3 | CN4 | CN5 | P1 | P2 | P3 | P4 | P5 |
|---|---|---|---|---|---|---|---|---|---|
| 203 000 | 21 240 | 80 232 | 54 446 | 24 060 | 190 249 | 80 254 | 170 675 | 8091 | 3169 |

Here is a plot of the lengths of axiom sequences for all the theorems, shown on a log scale:

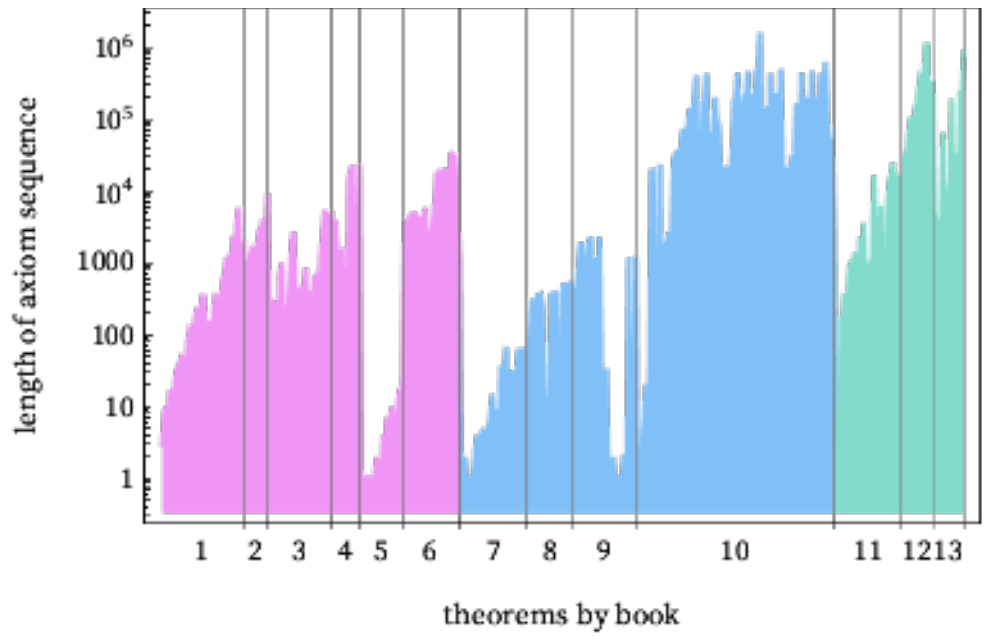

Interestingly, 3.18 isn't the theorem with the longest axiom sequence; it's in 4th place, and the top 10 are (in gray are the results with intermediate theorems allowed):

| 10.72 | 12.14 | 12.15 | 13.18 | 10.110 | 10.84 | 10.102 | 10.65 | 10.96 | 10.59 |
|---|---|---|---|---|---|---|---|---|---|
| 1 584 084 | 1 048 323 | 1 048 316 | 835 416 | 616 779 | 493 513 | 461 107 | 458 109 | 438 610 | 432 526 |
| *184* | *122* | *122* | *229* | *172* | *171* | *167* | *167* | *167* | *167* |

(10.72 is about addition of incommensurable medial areas, and is never referenced anywhere; 12.14 says the volumes of cones and cylinders with equal bases are proportional to their heights; 12.15 says the heights and bases of cones and cylinders with equal volumes are inversely proportional; etc.)

Here's the distribution of the lengths of axiom sequences across all theorems:

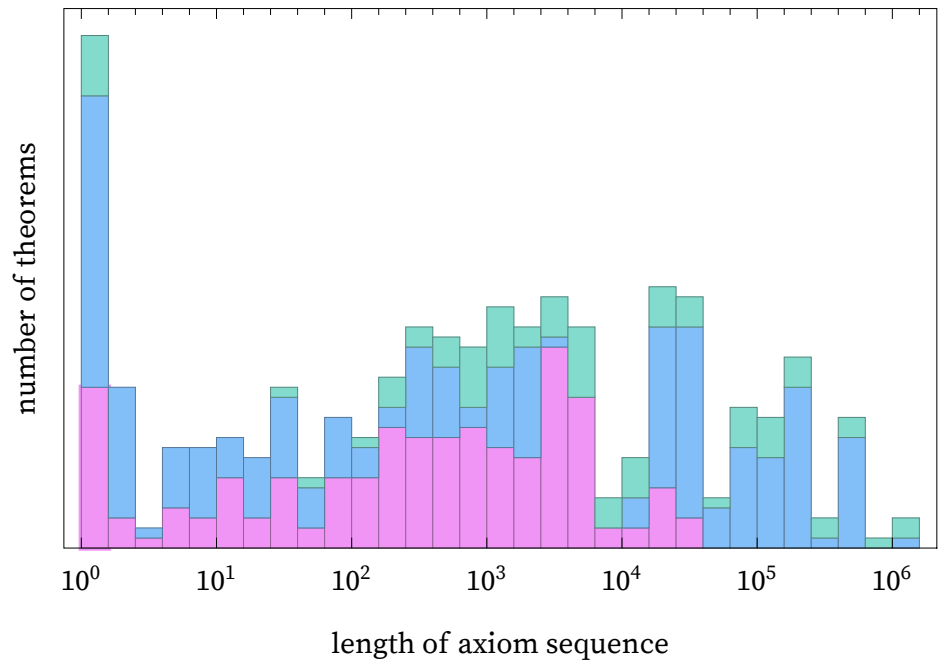

We can get some sense of the dramatic value of "remembering intermediate theorems" by comparing the total number of "intermediate steps" obtained with and without merging different instances of the same theorem:

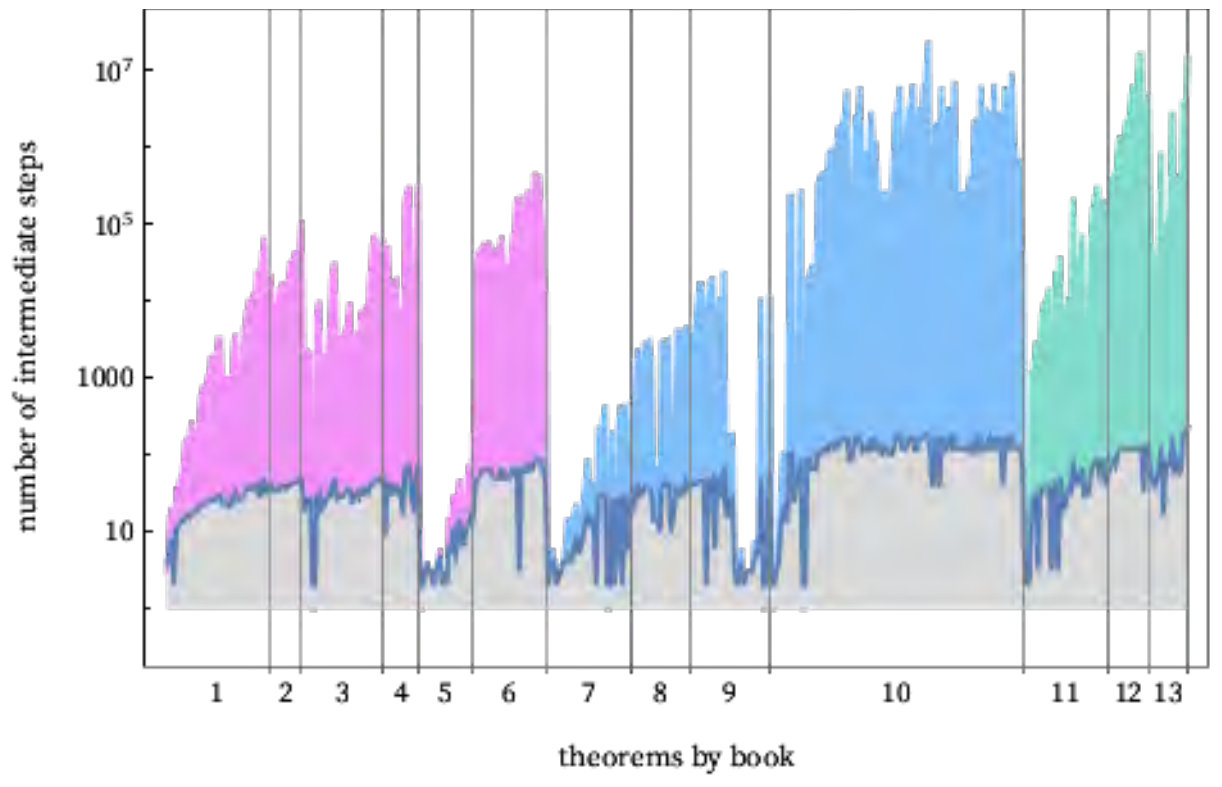

For example, for 8.13, 229 steps are needed when intermediate theorems are remembered, while 14,412,576 steps are needed otherwise. (For 10.72, it's 184 vs. 23,921,481 steps.)

## Superaxioms, or What Are the Most Powerful Theorems?

Euclid's 10 axioms are ultimately all we need in order to prove all the 465 theorems in the *Elements*. But what if we supplement these axioms with some of the theorems? Are there small sets of theorems we can add that will make the proofs of many theorems much shorter? To get a full understanding of this, we'd have to redo all the proofs. But we can get some sense of it just from the theorem dependency graph.

Consider the graph representing the proof of 1.12, with 1.7 highlighted:

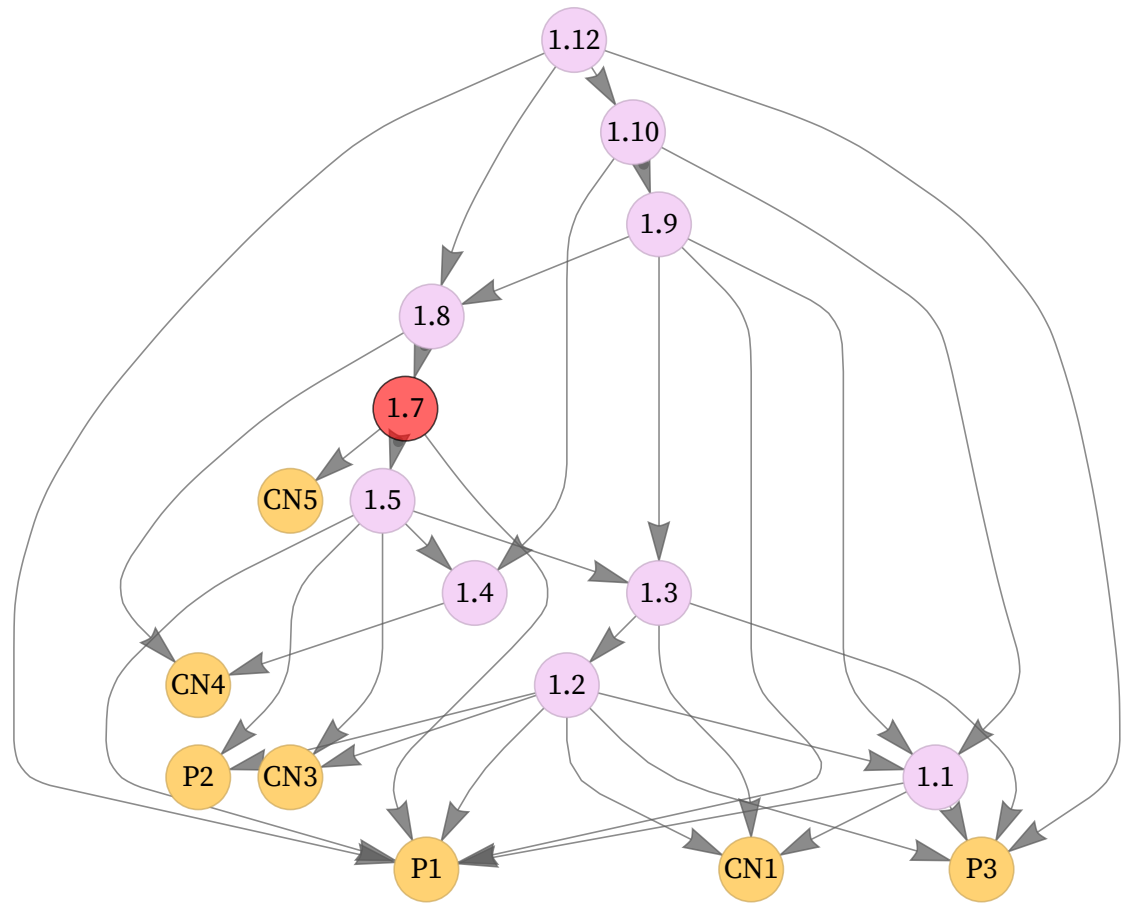

Now imagine adding 1.7 as a "superaxiom". Doing this, we can get a smaller proof graph for 1.12—with 4 nodes (and 14 connections) fewer:

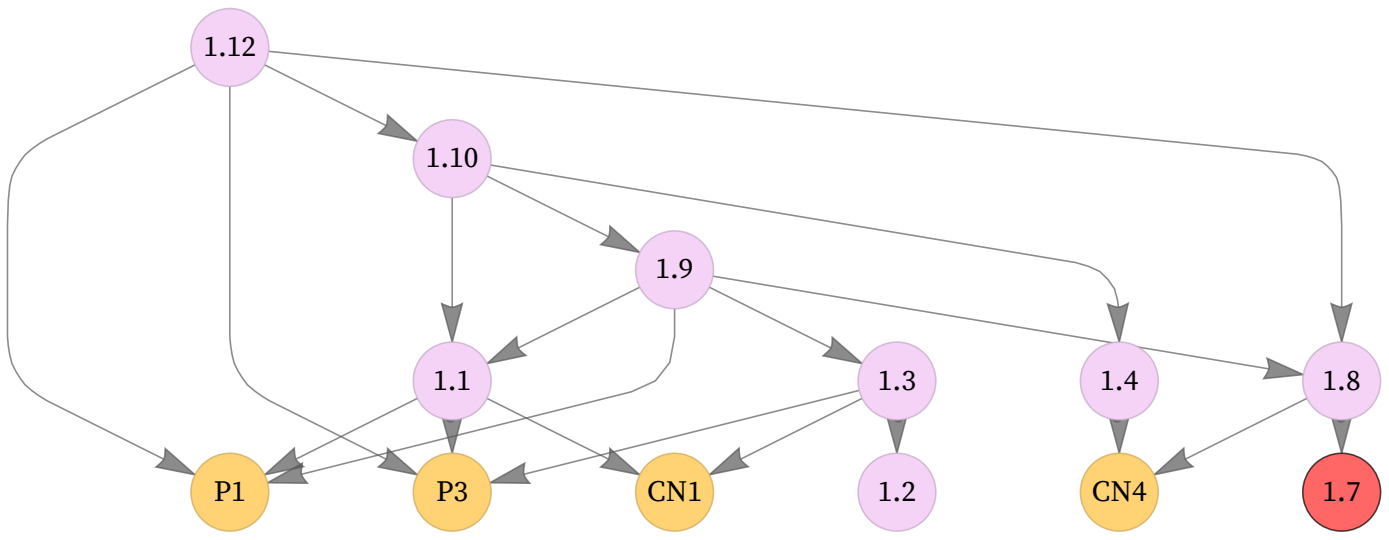

What does adding 1.7 as a superaxiom do for the proofs of other theorems? Here's how much it shortens each of them:

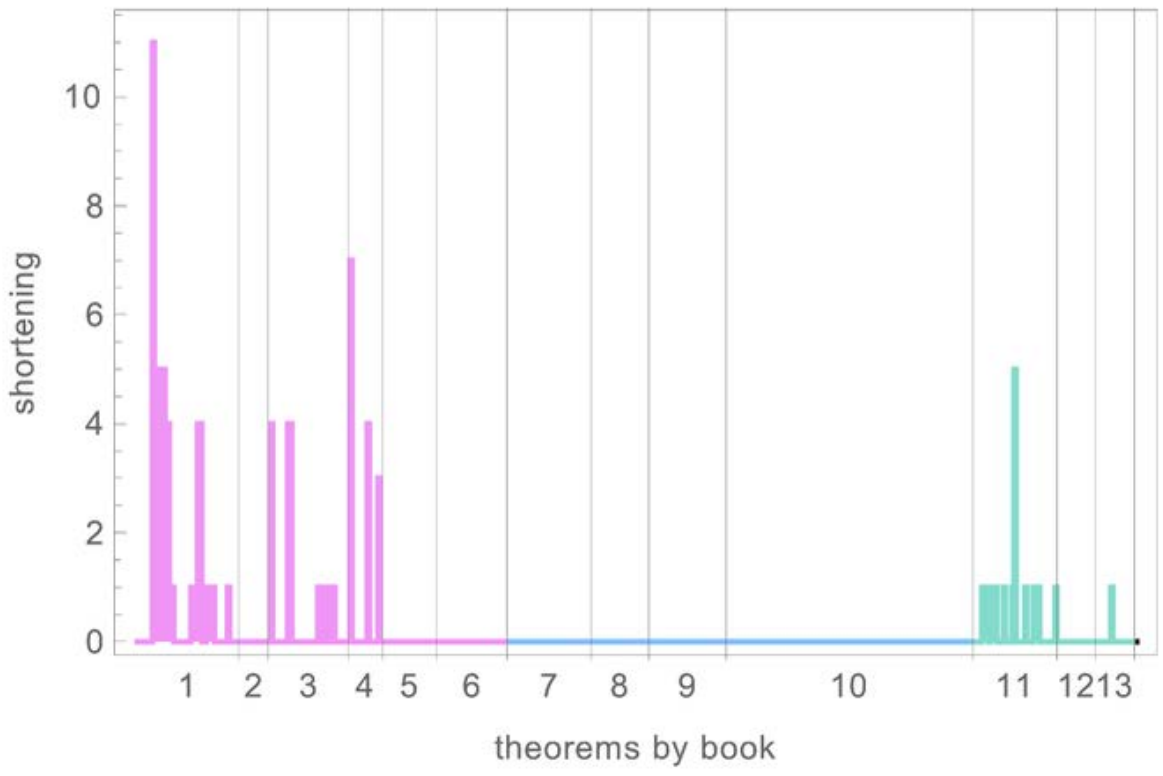

(The largest shortening is for 1.8, followed by 4.1.)

So what are the "best" superaxioms to add? Here's a plot of the average amount of shortening achieved by adding each possible individual theorem as a superaxiom:

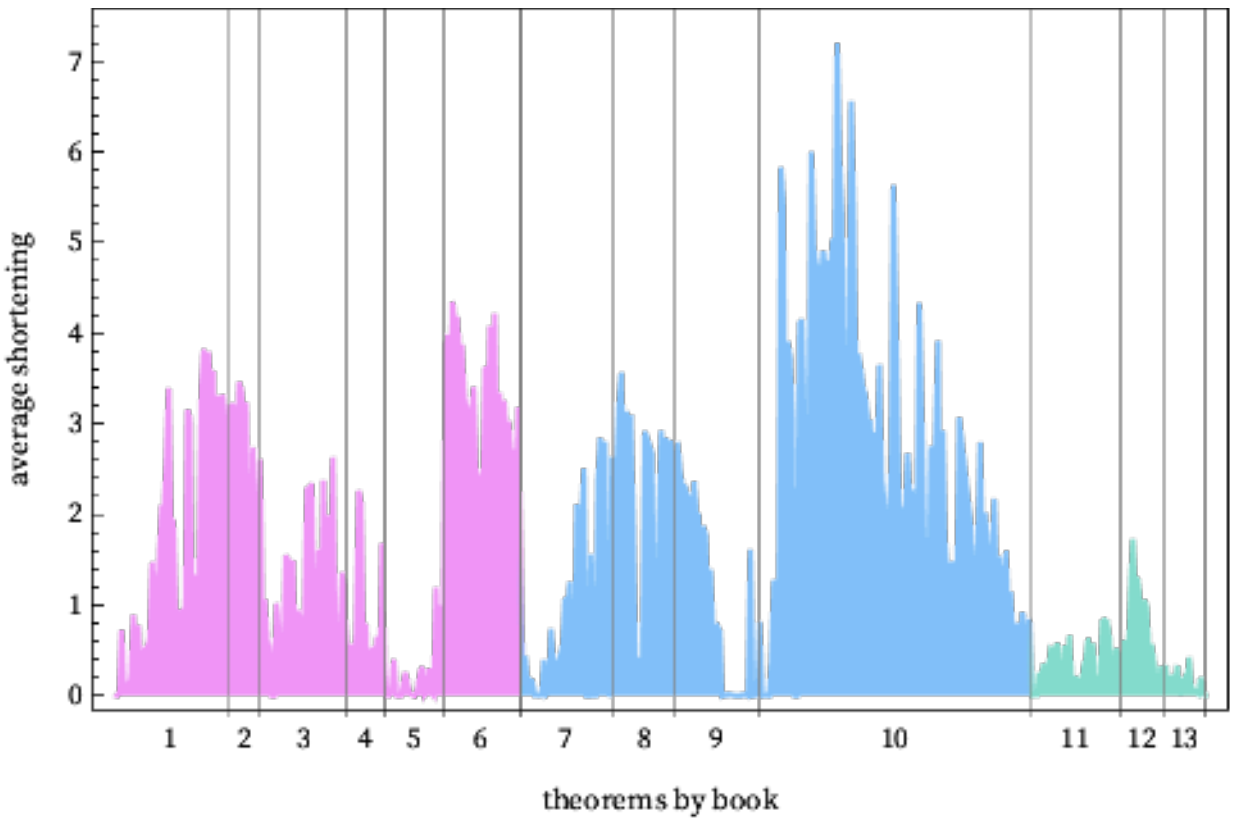

The rather unimpressive best result—an average shortening of 7.2—is achieved with 10.33 (which says that it's possible to come up with numbers $x$ and $y$ such that $\frac{x}{y}$ and $\sqrt{x\,y}$ are irrational, while $x\,y$ and $x + y$ are rational).

The maximum shortenings are more impressive—with 10.41 and 10.78 achieving the maximum shortening of 165

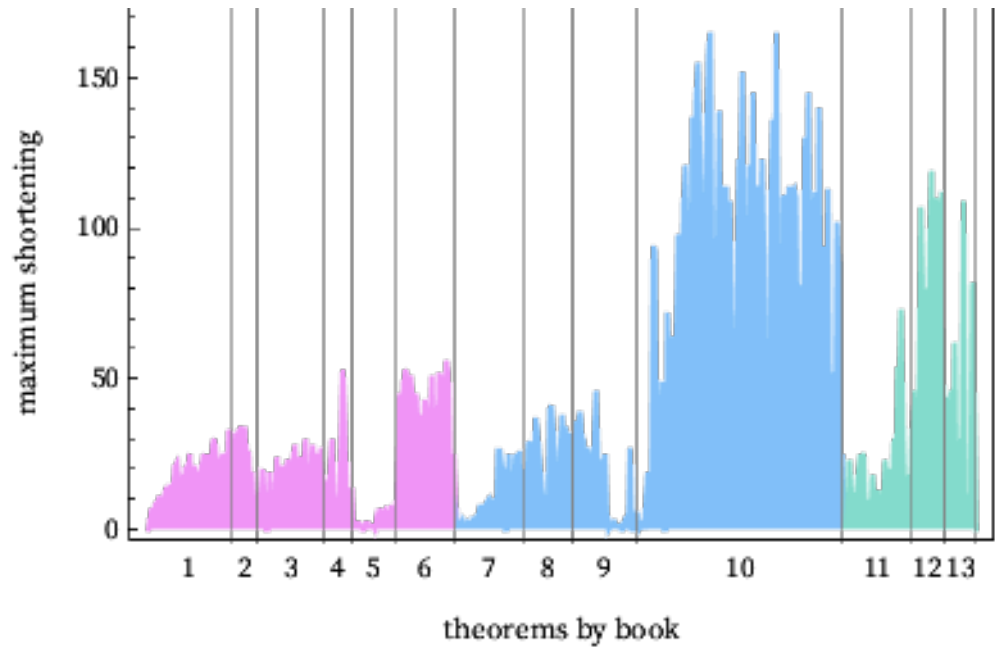

although this shortening is very concentrated around "nearby theorems":

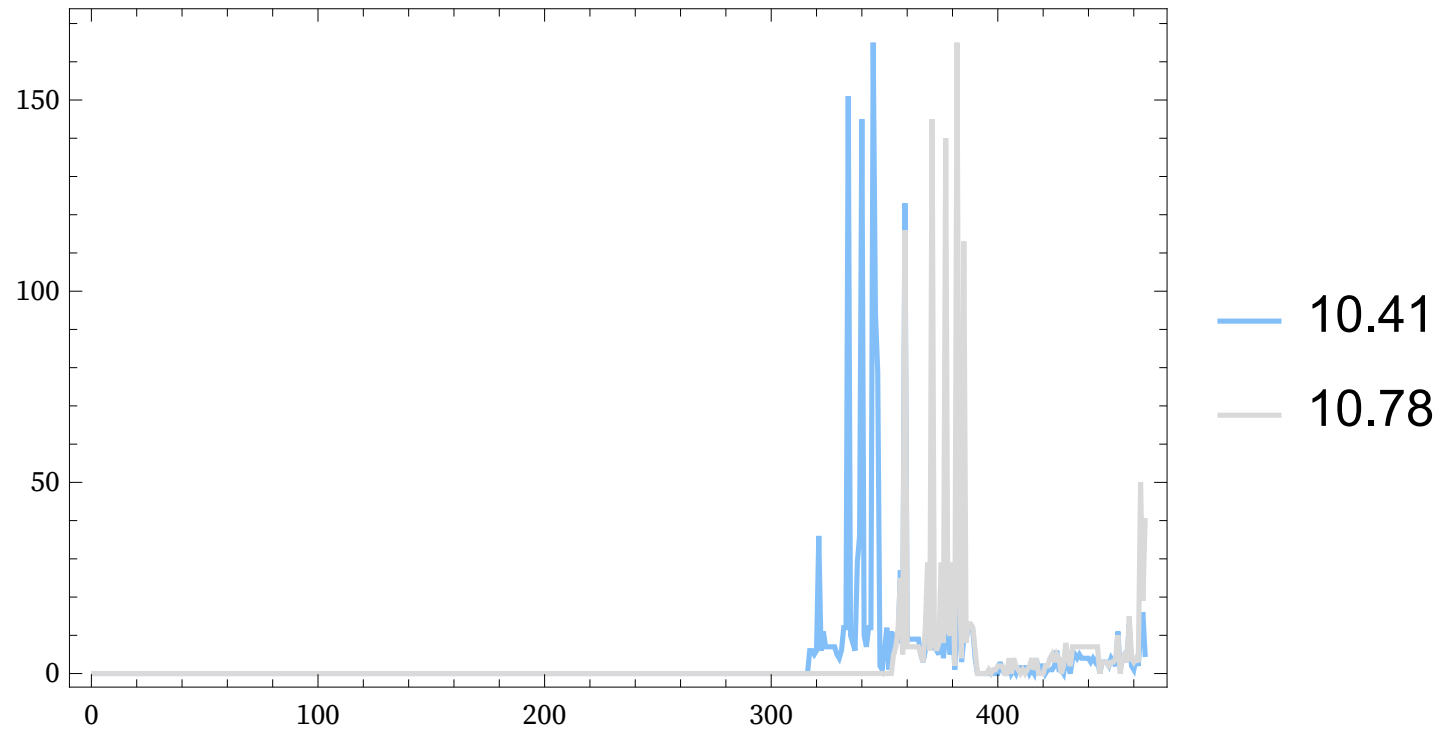

By the way, adding a superaxiom can not only decrease the number of intermediate theorems used in a proof, it can also decrease the "depth" of the proof, i.e. the longest needed to reach an axiom (or superaxiom). Here is the average depth reduction achieved by adding each possible theorem as a superaxiom:

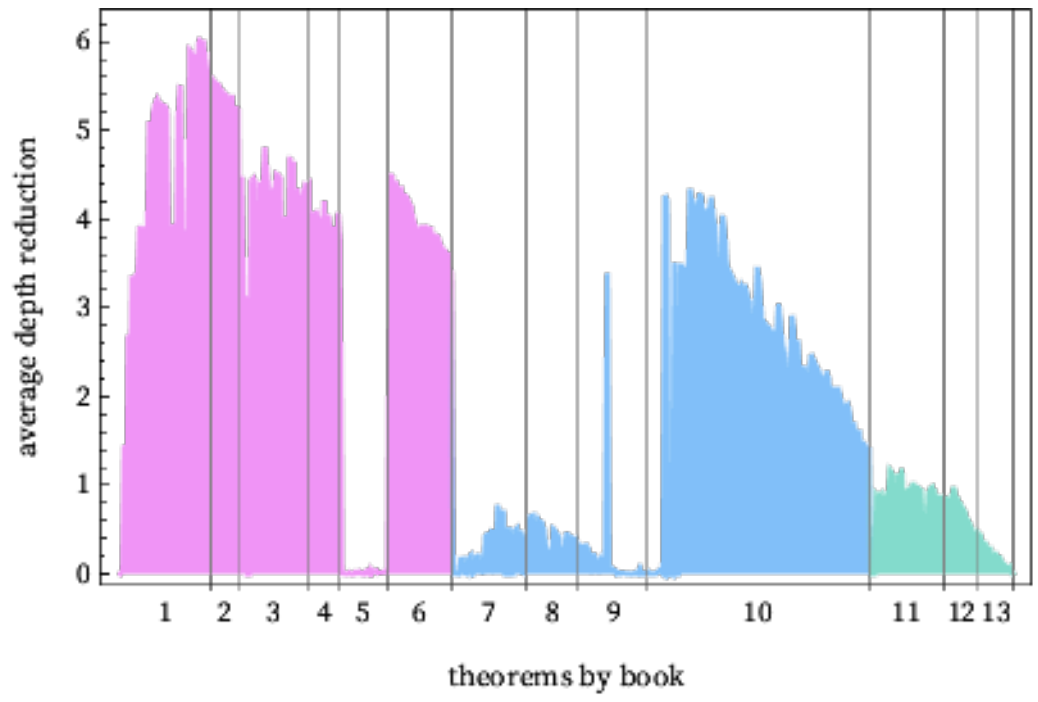

(The peak in Book 9 is 9.15, which reduces the depth of many subsequent theorems by 10 steps, though—in a possible goof—is not actually used by Euclid in the proofs of any of them.)

Here is the maximum depth reduction achieved by adding each possible theorem:

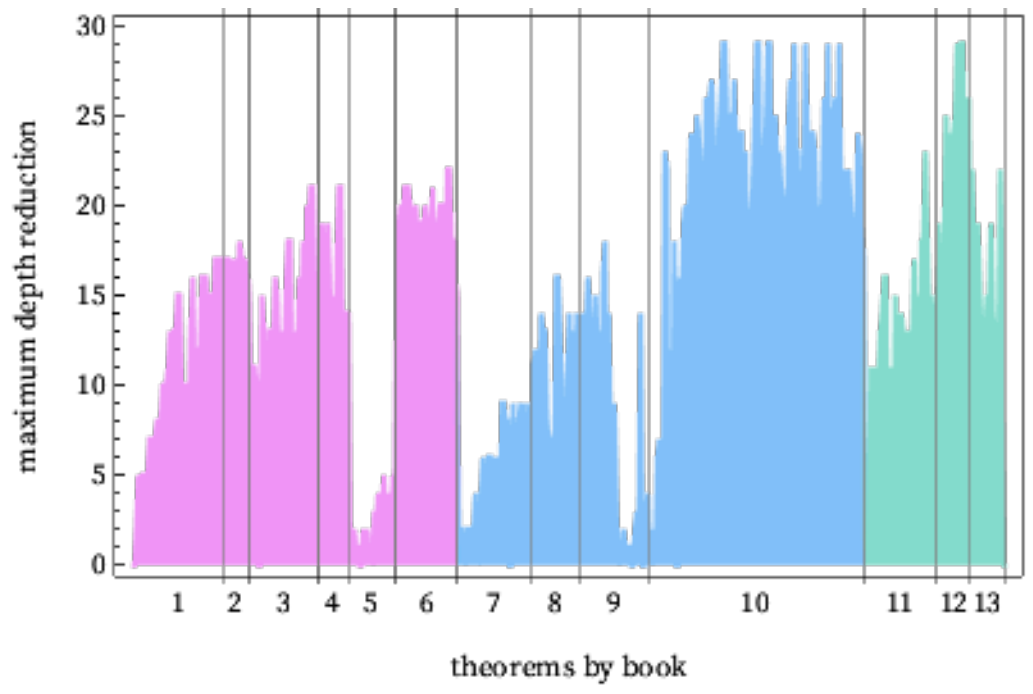

# Formalizing Euclid

Everything we've discussed so far is basically derived from the original text of Euclid's *Elements*. But what if we look instead at the pure "mathematical content" of Euclid? We've now got a way to represent this in the Wolfram Language. Consider Euclid's 3.16. It asserts that:

The perpendicular ($\overline{BI}$) to the diameter ($\overline{AB}$) of a circle at one of its endpoints ($B$) is tangent to the circle at that point.

Well, we can now give a "computational translation" of this:

In[∘]:= **Euclid book 3 proposition 16** GEOMETRIC SCENE **["Scene"]**

```
Out[∘]= GeometricScene[{{A, B, C, I}, {}}, {CircleThrough[{A, B}, C], Line[{A, C, B}],
    GeometricAssertion[{InfiniteLine[{B, I}], Line[{A, B}]}, Perpendicular]},
  {GeometricAssertion[{CircleThrough[{A, B}, C], InfiniteLine[{B, I}]}, Tangent]}]
```

And this is all we need to say to define that theorem in Euclid. Given the definition of the Wolfram Language, this is completely self-contained, and ready to be understood by both computers and humans. And from this form, we can now for example compute a random instance of the theorem:

In[∘]:= **RandomInstance[%]**

Out[∘]=

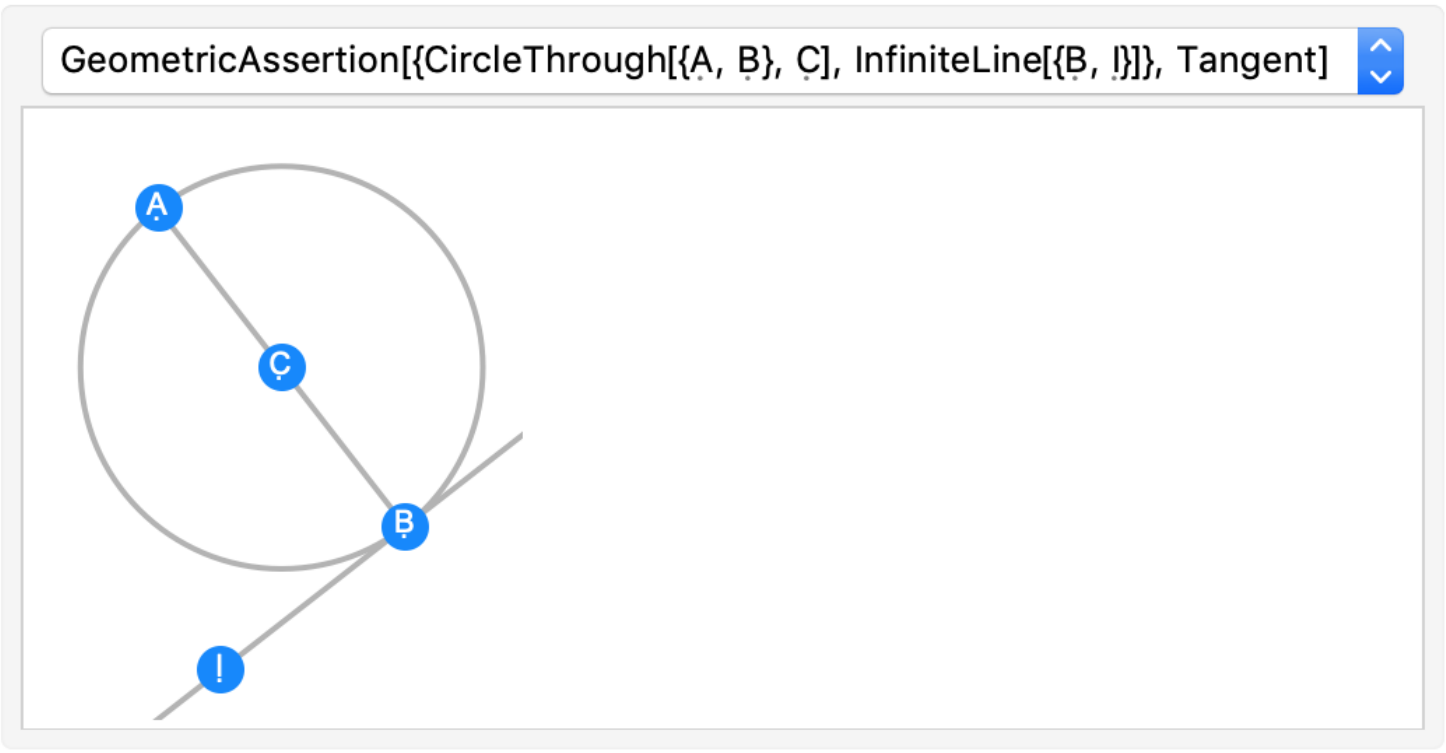

As another example, here's Euclid's 4.2:

To inscribe inside a given circle a triangle similar to a given triangle.

This is now asking for a construction—or, effectively, stating the theorem that it's possible to do such a construction with ruler and compass. And again we can give a computable version of this in the Wolfram Language, including the construction:

```
In[•]:= Euclid book 4 proposition 2 GEOMETRIC SCENE ["Scene"]

Out[•]= GeometricScene[{{A, O, B, C, D, E, F, G, H}, {}}, {{CircleThrough[{A}, O], Triangle[{D, E, F}]},
          {GeometricAssertion[{Line[{G, A, H}], CircleThrough[{A}, O]}, Tangent]},
          {C ∈ CircleThrough[{A}, O], Line[{A, C}], PlanarAngle[{C, A, H}] == PlanarAngle[{D, E, F}]},
          {B ∈ CircleThrough[{A}, O], Line[{A, B}], PlanarAngle[{B, A, G}] == PlanarAngle[{D, F, E}]}, {Line[{B, C}]}},
        {GeometricAssertion[{Triangle[{A, B, C}], Triangle[{D, E, F}]}, Similar]}]

In[•]:= RandomInstance[%]
```

Out[•]=

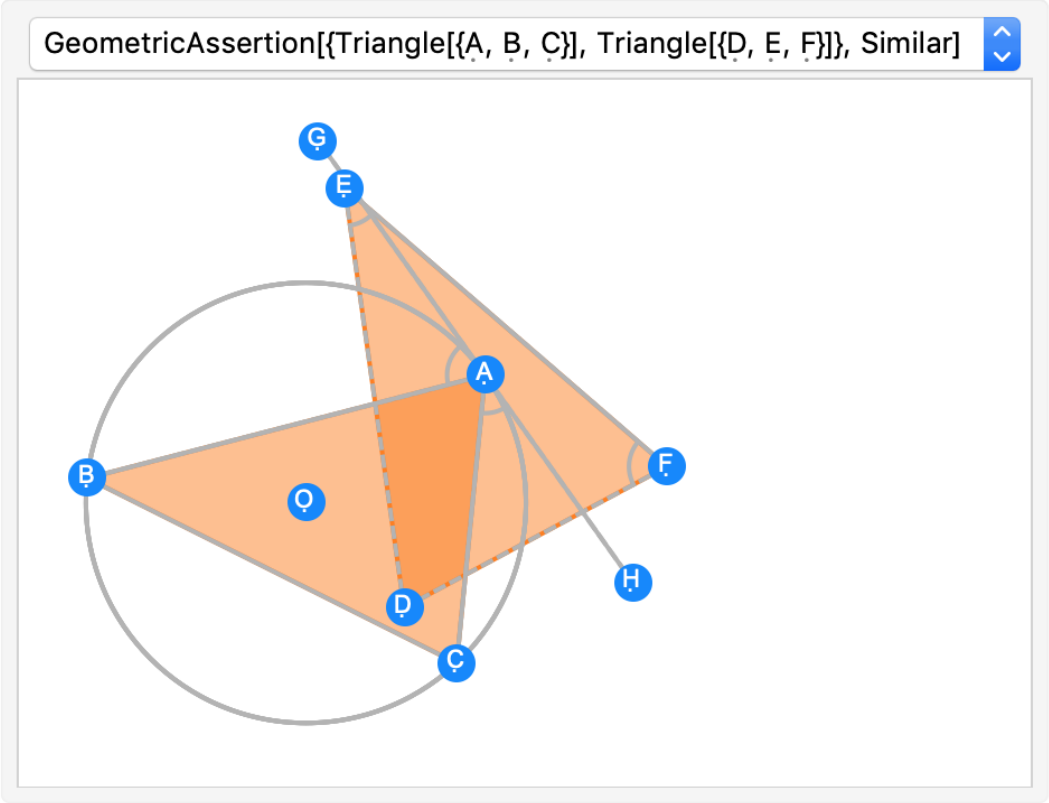

It's interesting to see, though, how the computable versions of theorems compare to their textual ones. Here are length comparisons for 2D geometry theorems:

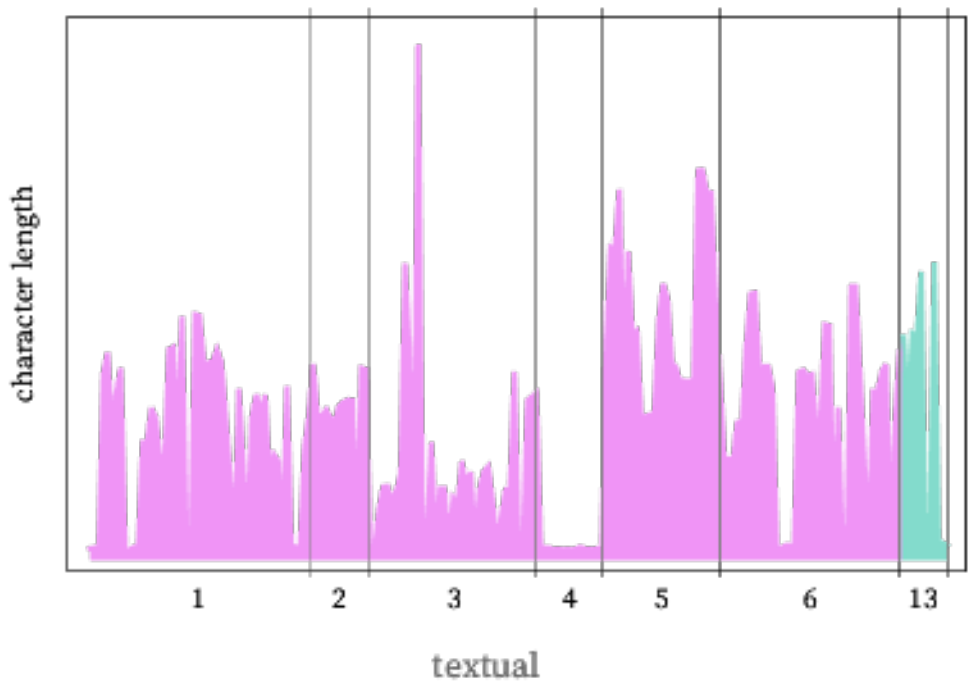

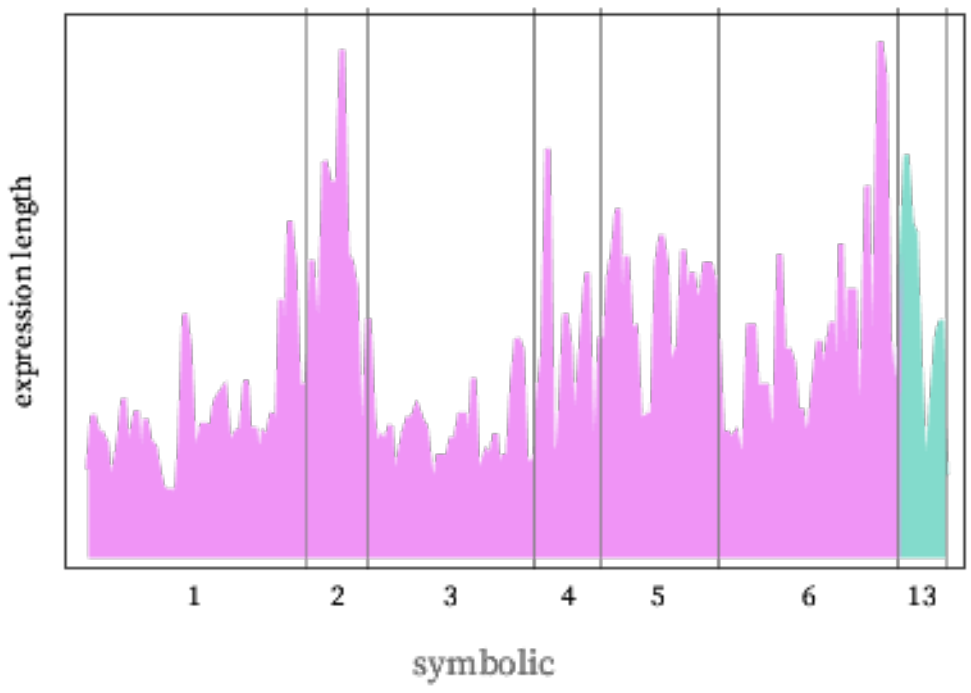

And we see that there is indeed at least some correlation between the lengths of textual and symbolic representations of theorems (the accumulation of points on the left is associated with constructions, where the text just says what's wanted, and the symbolic form also says how to do it):

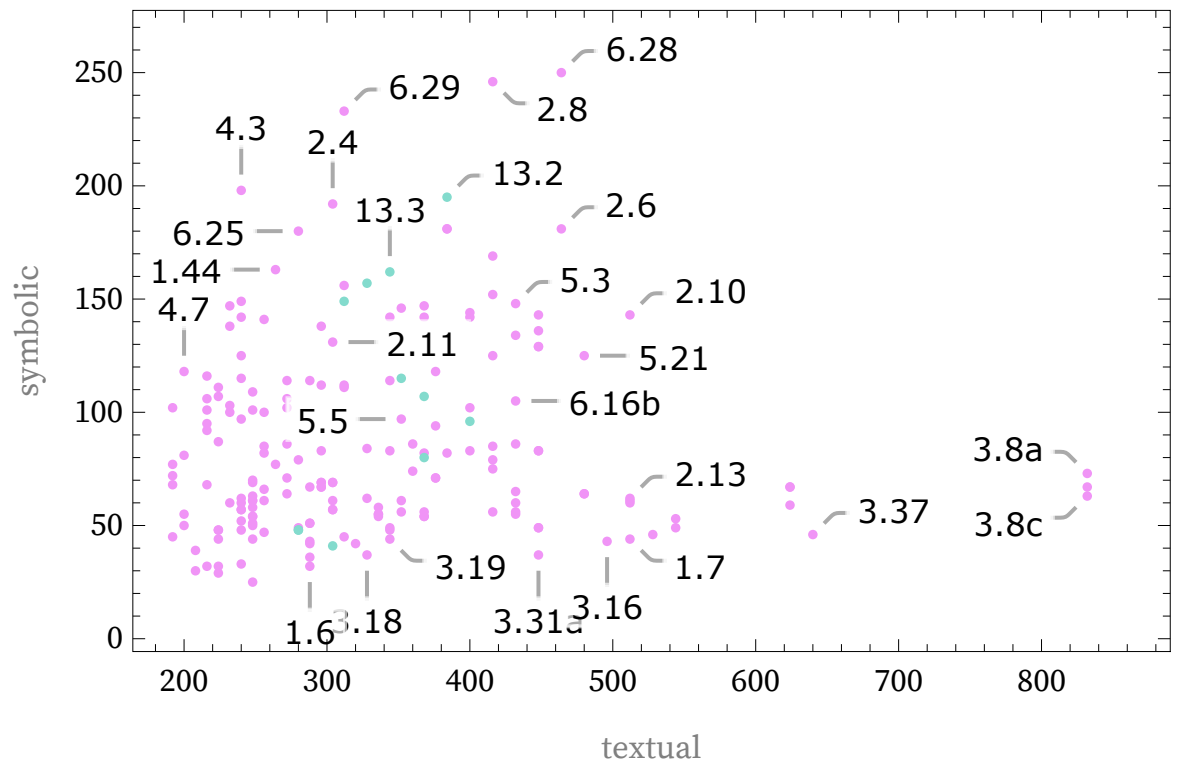

In the Wolfram Language representation we've just been discussing, there's a built-in Wolfram Language meaning to things like **CircleThrough** and **PlanarAngle**—and we can in a sense do general computations with these.

But at some level we can view what Euclid did as something purely formal. Yes, he talks about lines and planes. But we can think of these things just as formal constructs, without any externally known properties. Many centuries after Euclid, this became a much more familiar way to think about mathematics. And in the Wolfram Language we capture it with **AxiomaticTheory** and related functions.

For example, we can ask for an axiom system for Boolean algebra, or group theory:

```
In[]:= AxiomaticTheory["BooleanAxioms"]
```

Out[]= $\{\forall_{\{a,b\}}\ a\otimes b == b\otimes a, \forall_{\{a,b\}}\ a\oplus b == b\oplus a, \forall_{\{a,b\}}\ a\otimes(b\oplus\overline{b}) == a, \forall_{\{a,b\}}\ a\oplus b\otimes\overline{b} == a, \forall_{\{a,b,c\}}\ a\otimes(b\oplus c) == a\otimes b\oplus a\otimes c, \forall_{\{a,b,c\}}\ a\oplus b\otimes c == (a\oplus b)\otimes(a\oplus c)\}$

```
In[]:= AxiomaticTheory["GroupAxioms"]
```

Out[]= $\{\forall_{\{a,b,c\}}\ a\otimes(b\otimes c) == (a\otimes b)\otimes c, \forall_{a}\ a\otimes\bar{1} == a, \forall_{a}\ a\otimes\bar{a} == \bar{1}\}$

What does the ⊗ mean? We're not saying. We're just formally defining certain properties it's supposed to have. In the case of Boolean algebra, we can interpret it as **And**. In the case of group theory, it's group multiplication—though we're not saying what particular group it's for. And, yes, we could as well write the group theory axioms for example as:

$\{\forall_{\{a,b,c\}}\ f[a, f[b, c]] == f[f[a, b], c], \forall_{a}\ f[a, e] == a, \forall_{a}\ f[a, c[a]] == e\}$

OK, so can we do something similar for Euclid's geometry? It's more complicated, but thanks particularly to work by David Hilbert and Alfred Tarski in the first half of the 1900s, we can—and here's a version of the result:

| |
|---|
| $\forall_{\{x,y,z\}}$ implies[congruent[line[x, y], line[z, z]], congruent[x, y]] |
| $\forall_{\{x,y,z,u,v,w\}}$ implies[<br>and[congruent[line[x, y], line[z, u]], congruent[line[x, y], line[v, w]]], congruent[line[z, u], line[v, w]]] |
| $\forall_{\{x,y,z\}}$ implies[between[x, y, z], equal[x, y]] |
| $\forall_{\{x,y,z,u,v\}}$ implies[and[between[x, u, z], between[y, v, z]], $\exists_a$ and[between[u, a, y], between[v, a, x]]] |
| $\forall_{\{x,y,z,u,v\}}$ implies[and[and[and[congruent[line[x, u], line[x, v]], congruent[line[y, u], line[y, v]]],<br>congruent[line[z, u], line[z, v]]], not[equal[u, v]]],<br>or[or[between[x, y, z], between[y, z, x]], between[z, x, y]]] |
| $\forall_{\{x,y,z,u,v,w\}}$ implies[and[and[and[between[x, y, w], congruent[line[x, y], line[y, w]]],<br>and[between[x, u, v], congruent[line[x, u], line[u, v]]]],<br>and[between[y, u, z], congruent[line[y, u], line[z, u]]]], congruent[line[y, z], line[v, w]]] |
| $\forall_{\{x,y,z,a,b,c,u,v\}}$ implies[and[and[and[and[and[and[not[equal[x, y]], between[x, y, z]], between[a, b, c]],<br>congruent[line[x, y], line[a, b]]], congruent[line[y, z], line[b, c]]], congruent[line[x, u], line[a, v]]],<br>congruent[line[y, u], line[b, v]]], congruent[line[z, u], line[c, v]]] |
| $\forall_{\{x,y\}}$ implies[equal[x, y], equal[y, x]] |
| $\forall_{\{x,y,z\}}$ implies[and[equal[x, y], equal[y, z]], equal[x, z]] |
| $\forall_{x}$ equal[x, x] |
| $\forall_{\{a,b\}}$ and[a, b] == and[b, a] |
| $\forall_{\{a,b\}}$ or[a, b] == or[b, a] |
| $\forall_{\{a,b\}}$ and[a, or[b, not[b]]] == a |
| $\forall_{\{a,b\}}$ or[a, and[b, not[b]]] == a |
| $\forall_{\{a,b,c\}}$ and[a, or[b, c]] == or[and[a, b], and[a, c]] |
| $\forall_{\{a,b,c\}}$ or[a, and[b, c]] == and[or[a, b], or[a, c]] |
| $\forall_{\{a,b\}}$ implies[a, b] == or[not[a], b] |
| $\forall_{\{\alpha,\beta,y,z\}}$ implies[$\exists_x$ implies[and[$\alpha$[y], $\beta$[z]], between[x, y, z]], $\exists_u$ implies[and[$\alpha$[y], $\beta$[z]], between[y, u, z]]] |

Once again, this is all just a collection of formal statements. The fact that we're calling an operator between is just for our convenience and understanding. All we can really say for sure is that this is some ternary operator; any properties it has have to be defined by the axioms.

To get to this formalization of Euclid, quite a bit of tightening up had to be done. Euclid's theorems often had implicit assumptions, and it sometimes wasn't even clear exactly what their logical structure was supposed to be. But the mathematical content is presumably the same, and indeed some of Euclid's axioms (like CN1) say basically exactly the same as these. (An important addition to what Euclid explicitly said is the last axiom above, which states Euclid's implicit assumption—that I now believe to be incorrect for the physical universe—that space is continuous. Unlike other axioms, which just make statements "true for all values of …", this axiom makes a statement "true for all functions …".)

So what can we do with these axioms? Well, in principle we can prove any theorem in Euclidean geometry. Appending to the axioms (that we refer to—ignoring the last axiom—as geometry) an assertion that we can interpret as saying that if a point $y$ is between $x$ and $z$ and between $x$ and $w$, then either $z$ is between $y$ and $w$ or $w$ is between $y$ and $z$:

```
In[∘]:= FindEquationalProof[or[between[y, z, w], between[y, w, z]],
          Append[geometry, and[between[x, y, z], between[x, y, w]]]]
```

```
Out[∘]= ProofObject[ $Failed   Logic: Predicate/EquationalLogic   Steps: 225
                               Theorem: or[between[y, z, w], between[y, w, z]] ]
```

Here's a graph representing this proof:

```
In[∘]:= %["ProofGraph"]
```

Out[∘]=

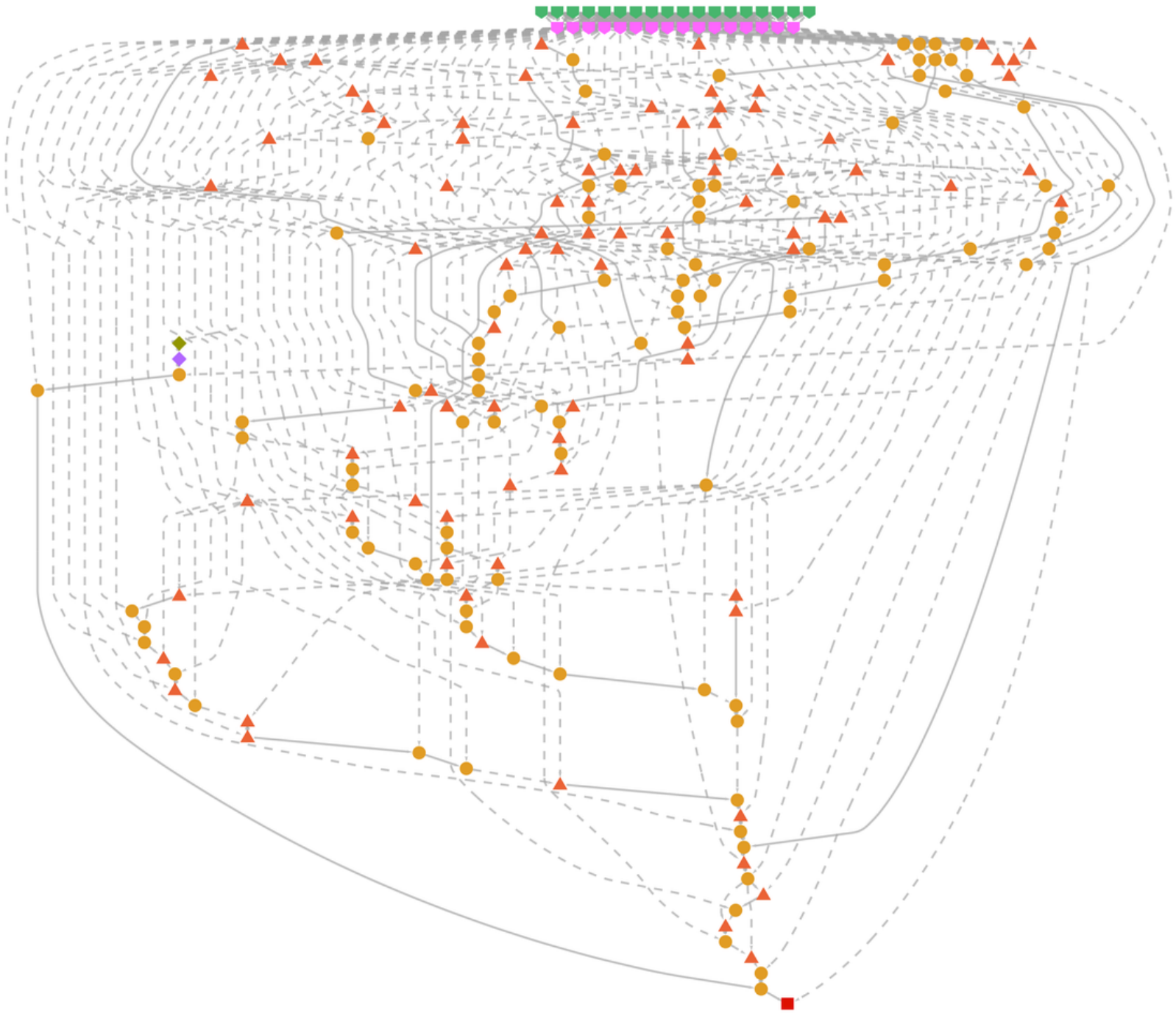

The axioms (including the "setup assertion") are at the top—and the proof, with all its various intermediate lemmas, establishes that our "hypothesis" (represented by a little purple diamond on the left) eventually leads to "true" at the bottom.

As a more complicated example, we can look at Euclid's very first theorem, 1.1, which asserts that there's a ruler-and-compass way to construct an equilateral triangle on any line segment. In the Wolfram Language, the construction is:

In[∘]:= Euclid book 1 proposition 1 GEOMETRIC SCENE ["Scene"]

Out[∘]= GeometricScene[{{A, B, C}, {}},
{{Line[{A, B}]}, {GeometricAssertion[{CircleThrough[{B}, A], CircleThrough[{A}, B]}, {Concurrent, C}]},
{Line[{C, A}], Line[{C, B}]}}, {GeometricAssertion[Triangle[{A, B, C}], Equilateral]}]

In[∘]:= **RandomInstance[%]**

Out[∘]=

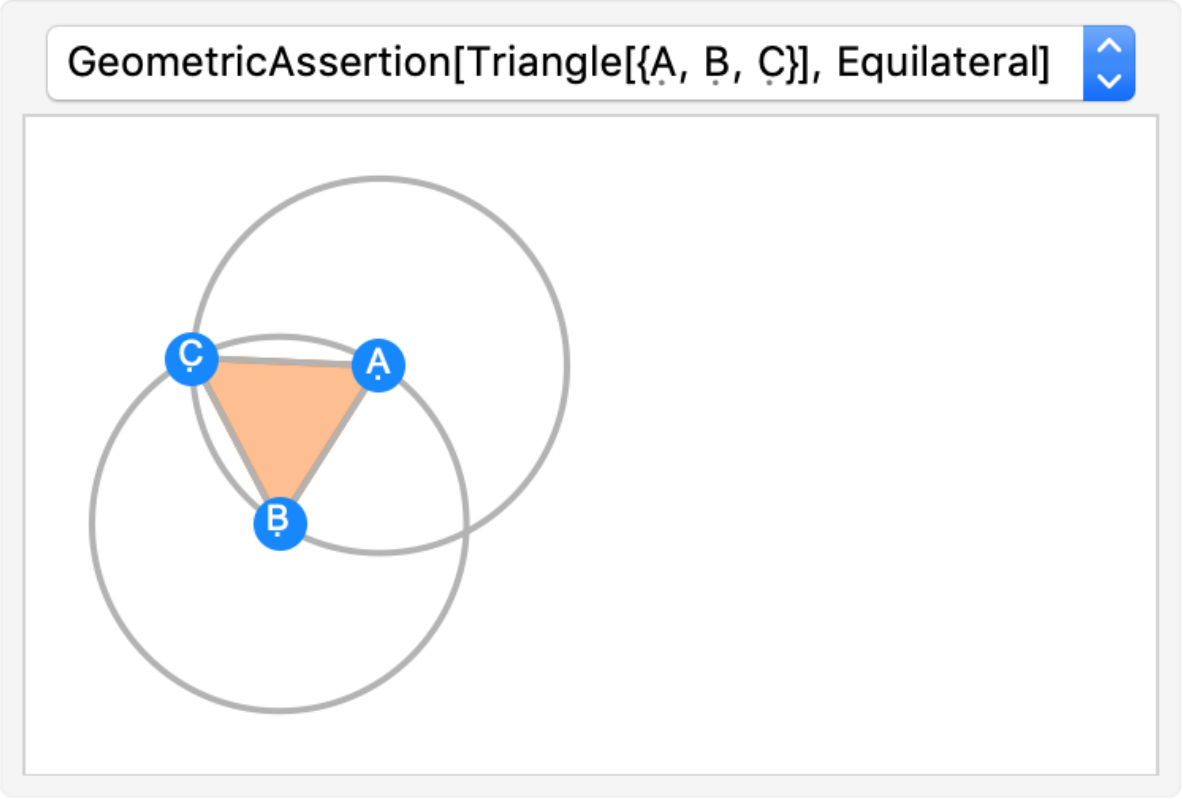

And now we can write this directly in terms of our low-level constructs. First we need a definition of what circles are (Euclid has this as Definition 1.15)—basically saying that two circles centered at *a* that go through *b* and *c* are equal if the lines from *a* to *b* and *a* to *c* are congruent:

$\forall_{\{a,b,c\}}$ implies[equal[circle[a, b], circle[a, c]], congruent[line[a, b], line[a, c]]]

We'll call this definition circles. We're going to do a construction that involves having circles that overlap, as specified by the assertions:

{equal[circle[a, b], circle[a, c]], equal[circle[b, a], circle[b, c]]}

And then our goal is to show that we get an equilateral triangle, for which the following is true:

and[congruent[line[a, b], line[a, c]], congruent[line[b, a], line[b, c]]]

Putting this all together we can prove Euclid's 1.1:

In[∘]:= **FindEquationalProof[and[congruent[line[a, b], line[a, c]], congruent[line[b, a], line[b, c]]],**
**Join[geometry, {$\forall_{\{a,b,c\}}$ implies[equal[circle[a, b], circle[a, c]], congruent[line[a, b], line[a, c]]]},**
**{equal[circle[a, b], circle[a, c]], equal[circle[b, a], circle[b, c]]}]]**

Out[∘]= ProofObject[ Logic: PredicateEquationalLogic Steps: 272 Theorem: and[congruent[line[a, b], line[a, c]], congruent[line[b, a], line[b, c]]] ]

And, yes, it took 272 steps—and here's a graphical representation of the proof that got generated, with all its intermediate lemmas:

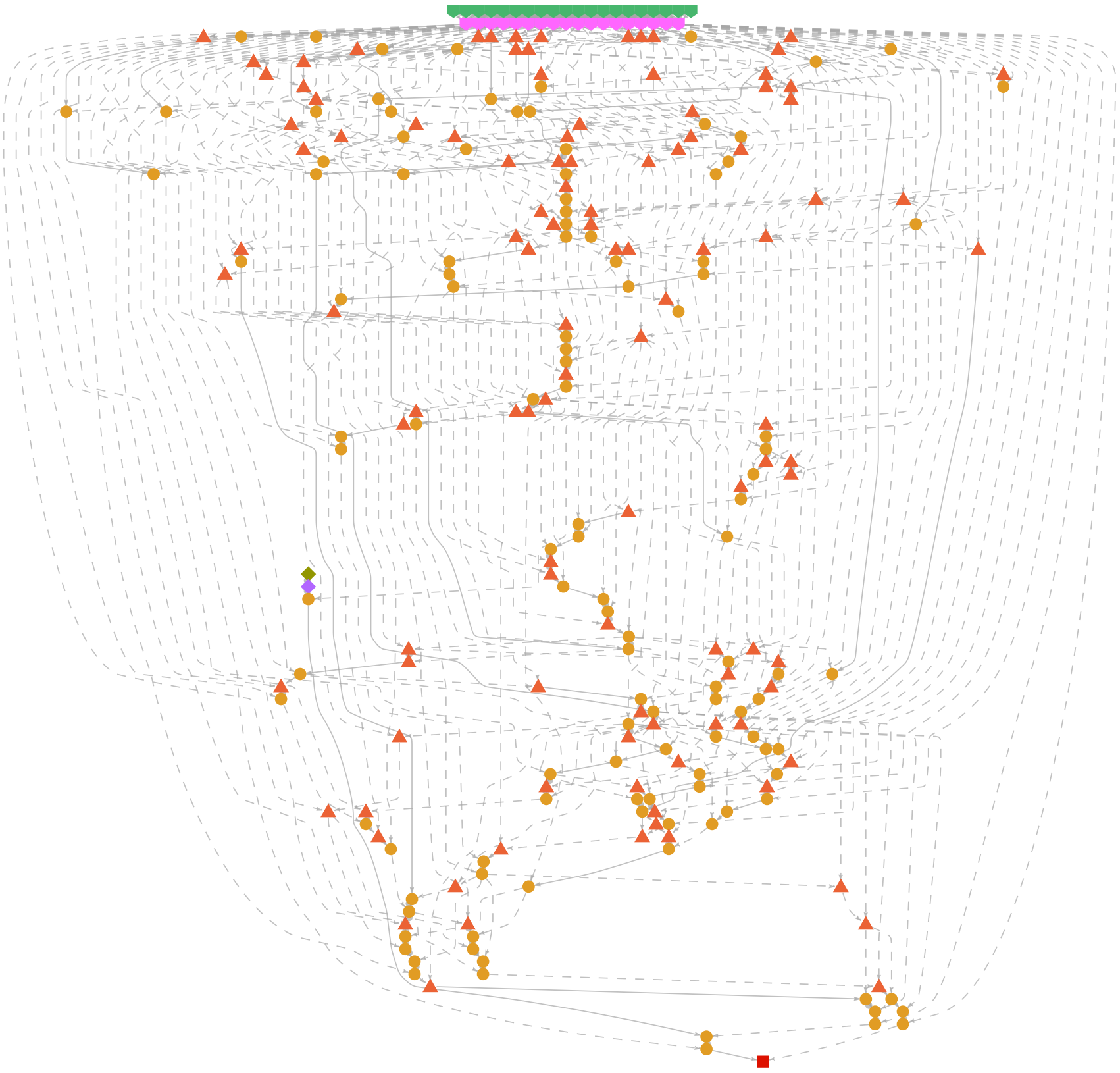

We can go on and prove Euclid's 1.2 as well, all the way from the lowest-level axioms. This time it takes us 330 steps, with proof graph:

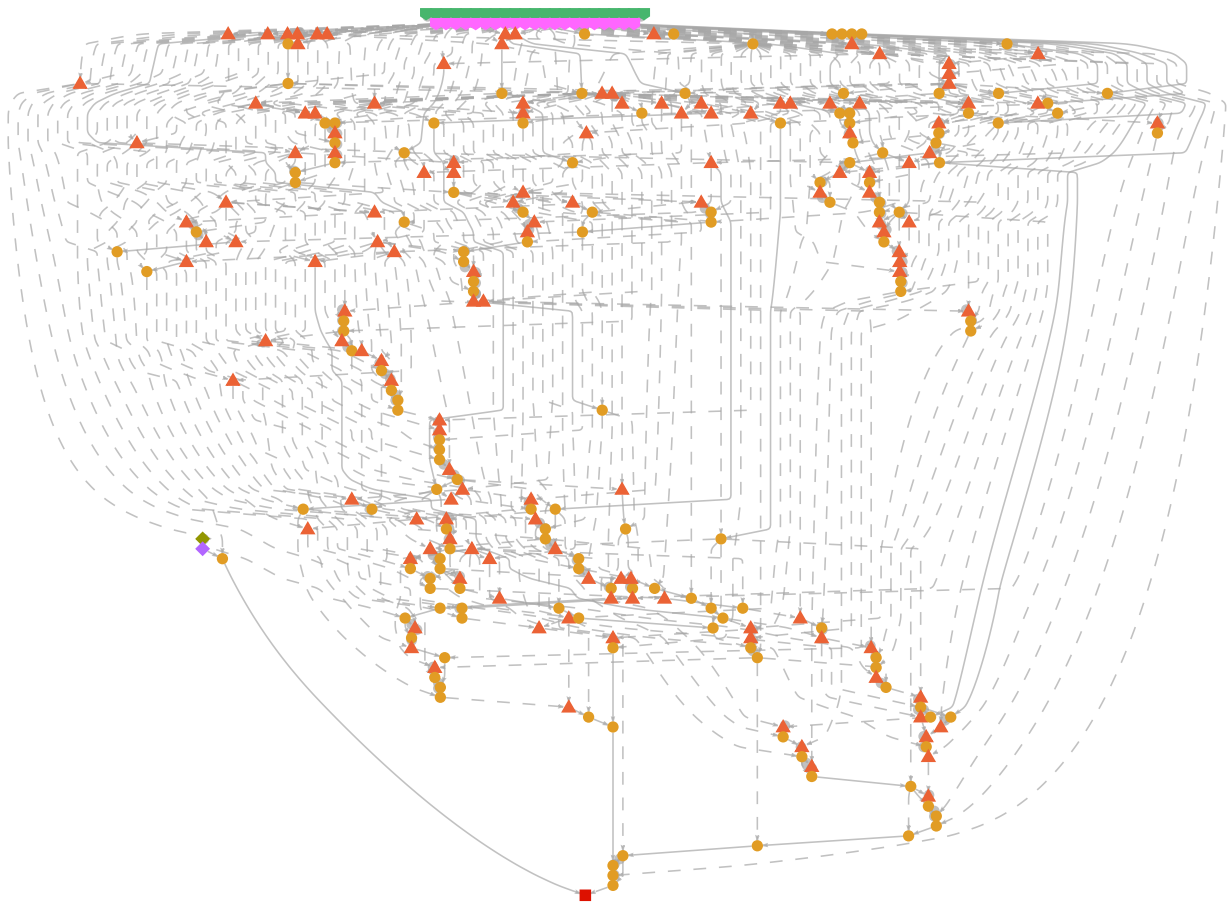

These graphs are conceptually similar to, but concretely rather different from, our "empirical metamathematics" graphs above. There are differences at the level of how interdependence of theorems is defined. But, more important, this graph is generated by automated theorem proving methods; the intermediate theorems (or lemmas) it involves are produced "on the fly" for the convenience of the computer, not because they help in any way to explain the proof to a human. In our empirical metamathematics on Euclid's *Elements*, however, we're dealing with the theorems that Euclid chose to define, and that have served as a basis for explaining his proofs to humans for more than two thousand years.

By the way, if our goal is simply to find out what's true in geometry—rather than to write out step-by-step proofs—then we now know how to do that. Essentially it involves turning geometric assertions into algebraic ones—and then systematically solving the polynomial equations and inequalities that result. It can be computationally expensive, but in the Wolfram Language we now have one master function, **CylindricalDecomposition**, that ultimately does the job. And, yes, given Gödel's theorem, one might wonder whether this kind of finite procedure for solving any Euclid-style geometry problem was even possible. But it turns out that—unlike arithmetic, for which Gödel's theorem was originally proved—Euclid-style geometry, like basic logic, is decidable, in the sense that there is ultimately a finite procedure for deciding whether any given statement is true or not. In principle, this procedure could be based on theorem proving from the axioms, but **CylindricalDecomposition** effectively leverages a tower of more sophisticated mathematics to provide a much more efficient approach.

## All Possible Theorems

From the axioms of geometry one can in principle derive an infinite number of true theorems—of which Euclid picked just 465 to include in his *Elements*. But why these theorems, and not others? Given a precise symbolic representation of geometry—as in the axioms above—one can just start enumerating true theorems.

One way to do this is to use a multiway system, with the axioms defining transformation rules that one can apply in all possible ways. In effect this is like constructing every possible proof, and seeing what gets proved. Needless to say, the network that gets produced quickly becomes extremely large—even if its structure is interesting for our attempt to find a "bulk theory of metamathematics".

Here's an example of doing it, not for the full geometry axioms above, but for basic logic (which is actually part of the axiom system we've used for geometry). We can either start with expressions, or with statements. Here we start with the expression $x \wedge y$, and then progressively find all expressions equal to it. Here's the first, rather pedantic step:

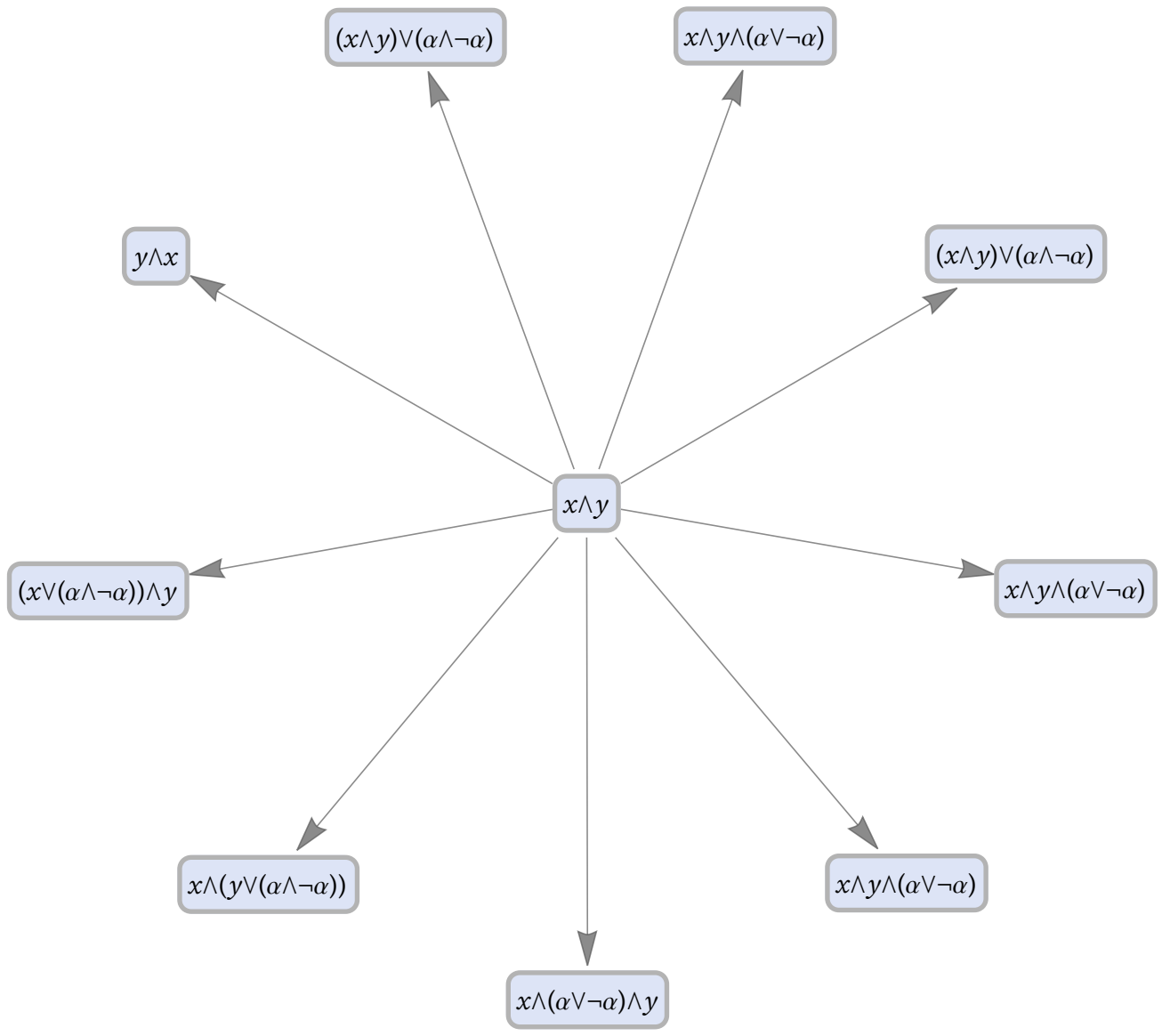

And here's the second step:

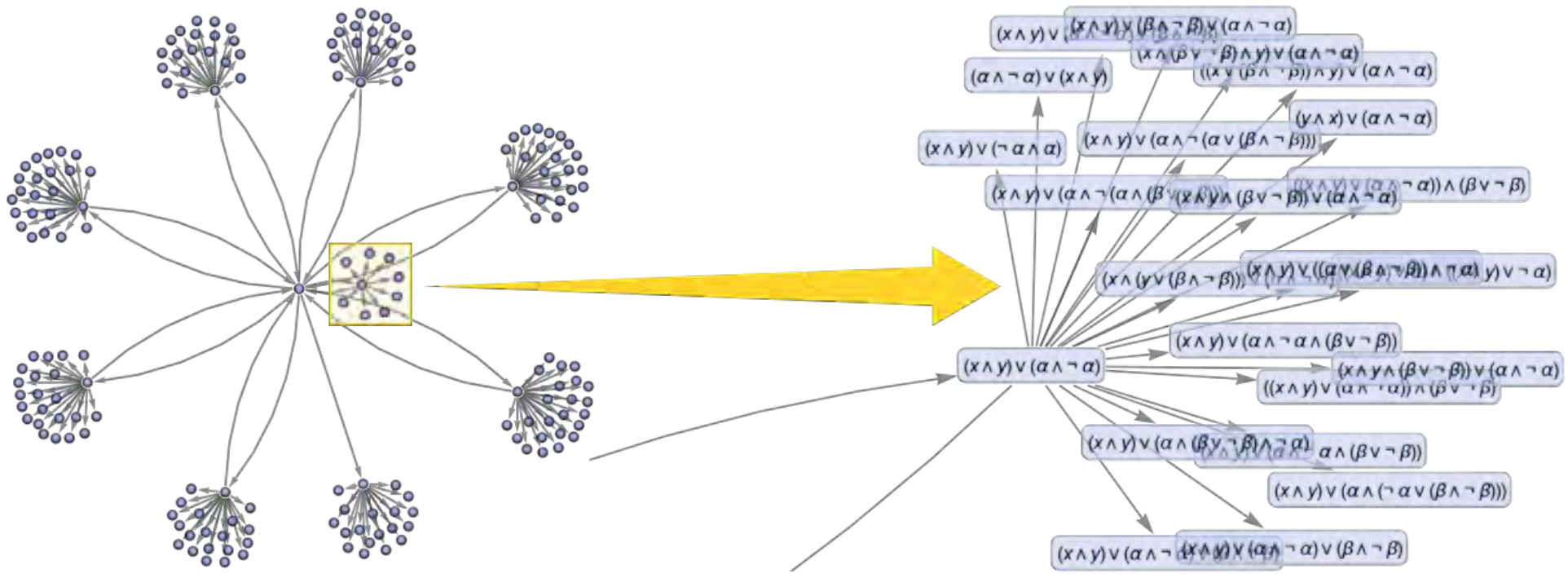

Every path in this graph is a proof that its endpoint expressions are equal. And while eventually this approach will give us every possible theorem (in this case about equalities involving $x \wedge y$), it'll obviously take a while, generating huge numbers of long and uninteresting results on its way to anything interesting.

As a different approach, we can consider just enumerating short possible statements, then picking out ones that we determine are true. In principle we could determine truth by explicitly proving theorems using the axioms (and, yes, if there was undecidability we wouldn't always be able to do this). But in practice for the case of basic logic that we're using as an example here, we can basically just explicitly construct truth tables to find out what's true and what's not.

Here are some statements in logic, sorted in increasing order of complexity (as measured by depth and number of symbols):

$$
\begin{aligned}
&\{a = b,\ a = (\neg a),\ a = (\neg b),\ (\neg a) = (\neg b),\ a = (a \land a),\ (\neg a) = (a \land a),\ a = (a \lor a),\ (\neg a) = (a \lor a),\ (a \land a) = (a \lor a),\\
&a = (a \land b),\ (\neg a) = (a \land b),\ (a \land a) = (a \land b),\ (a \lor a) = (a \land b),\ a = (a \lor b),\ (\neg a) = (a \lor b),\ (a \land a) = (a \lor b),\\
&(a \lor a) = (a \lor b),\ (a \land b) = (a \lor b),\ (a \land b) = (a \land c),\ (a \lor b) = (a \land c),\ (a \land b) = (a \lor c),\ (a \lor b) = (a \lor c),\\
&a = (b \land a),\ (\neg a) = (b \land a),\ (a \land a) = (b \land a),\ (a \lor a) = (b \land a),\ (a \land b) = (b \land a),\ (a \lor b) = (b \land a),\\
&a = (b \lor a),\ (\neg a) = (b \lor a),\ (a \land a) = (b \lor a),\ (a \lor a) = (b \lor a),\ (a \land b) = (b \lor a),\ (a \lor b) = (b \lor a),\\
&a = (b \land b),\ (\neg a) = (b \land b),\ (a \land a) = (b \land b),\ (a \lor a) = (b \land b),\ (a \land b) = (b \land b),\ (a \lor b) = (b \land b),\\
&a = (b \lor b),\ (\neg a) = (b \lor b),\ (a \land a) = (b \lor b),\ (a \lor a) = (b \lor b),\ (a \land b) = (b \lor b),\ (a \lor b) = (b \lor b),\ a = (b \land c),\\
&(\neg a) = (b \land c),\ (a \land a) = (b \land c),\ (a \lor a) = (b \land c),\ (a \land b) = (b \land c),\ (a \lor b) = (b \land c),\ a = (b \lor c),\ (\neg a) = (b \lor c),\\
&(a \land a) = (b \lor c),\ (a \lor a) = (b \lor c),\ (a \land b) = (b \lor c),\ (a \lor b) = (b \lor c),\ (a \land b) = (c \land a),\ (a \lor b) = (c \land a),\\
&(a \land b) = (c \lor a),\ (a \lor b) = (c \lor a),\ (a \land b) = (c \land b),\ (a \lor b) = (c \land b),\ (a \land b) = (c \lor b),\ (a \lor b) = (c \lor b),\\
&(a \land b) = (c \land c),\ (a \lor b) = (c \land c),\ (a \land b) = (c \lor c),\ (a \lor b) = (c \lor c),\ a = (\neg (\neg a)),\ (\neg a) = (\neg (\neg a)),\\
&(a \land a) = (\neg (\neg a)),\ (a \lor a) = (\neg (\neg a)),\ (a \land b) = (\neg (\neg a)),\ (a \lor b) = (\neg (\neg a)),\ a = (\neg (\neg b)),\ (\neg a) = (\neg (\neg b)),\\
&(a \land a) = (\neg (\neg b)),\ (a \lor a) = (\neg (\neg b)),\ (a \land b) = (\neg (\neg b)),\ (a \lor b) = (\neg (\neg b)),\ (\neg (\neg a)) = (\neg (\neg b)),\\
&(a \land b) = (\neg (\neg c)),\ (a \lor b) = (\neg (\neg c)),\ a = (\neg a \land a),\ (\neg a) = (\neg a \land a),\ (a \land a) = (\neg a \land a),\ (a \lor a) = (\neg a \land a),\\
&(a \land b) = (\neg a \land a),\ (a \lor b) = (\neg a \land a),\ (\neg (\neg a)) = (\neg a \land a),\ a = (a \land \neg a),\ (\neg a) = (a \land \neg a),\ (a \land a) = (a \land \neg a),\\
&(a \lor a) = (a \land \neg a),\ (a \land b) = (a \land \neg a),\ (a \lor b) = (a \land \neg a),\ (\neg (\neg a)) = (a \land \neg a),\ (\neg a \land a) = (a \land \neg a)\}
\end{aligned}
$$

Many (like $a=b$) are very obviously not true, at least not for all possible values of each variable. But—essentially by using truth tables—we can readily pick out ones that are always true:

$$
\begin{aligned}
&\{a = b,\ a = (\neg a),\ a = (\neg b),\ (\neg a) = (\neg b),\ \boxed{a = (a \land a)},\ (\neg a) = (a \land a),\ \boxed{a = (a \lor a)},\ (\neg a) = (a \lor a),\\
&\boxed{(a \land a) = (a \lor a)},\ a = (a \land b),\ (\neg a) = (a \land b),\ (a \land a) = (a \land b),\ (a \lor a) = (a \land b),\ a = (a \lor b),\\
&(\neg a) = (a \lor b),\ (a \land a) = (a \lor b),\ (a \lor a) = (a \lor b),\ (a \land b) = (a \lor b),\ (a \land b) = (a \land c),\\
&(a \lor b) = (a \land c),\ (a \land b) = (a \lor c),\ (a \lor b) = (a \lor c),\ a = (b \land a),\ (\neg a) = (b \land a),\ (a \land a) = (b \land a),\\
&(a \lor a) = (b \land a),\ \boxed{(a \land b) = (b \land a)},\ (a \lor b) = (b \land a),\ a = (b \lor a),\ (\neg a) = (b \lor a),\ (a \land a) = (b \lor a),\\
&(a \lor a) = (b \lor a),\ (a \land b) = (b \lor a),\ \boxed{(a \lor b) = (b \lor a)},\ a = (b \land b),\ (\neg a) = (b \land b),\ (a \land a) = (b \land b),\\
&(a \lor a) = (b \land b),\ (a \land b) = (b \land b),\ (a \lor b) = (b \land b),\ a = (b \lor b),\ (\neg a) = (b \lor b),\ (a \land a) = (b \lor b),\\
&(a \lor a) = (b \lor b),\ (a \land b) = (b \lor b),\ (a \lor b) = (b \lor b),\ a = (b \land c),\ (\neg a) = (b \land c),\ (a \land a) = (b \land c),\\
&(a \lor a) = (b \land c),\ (a \land b) = (b \land c),\ (a \lor b) = (b \land c),\ a = (b \lor c),\ (\neg a) = (b \lor c),\ (a \land a) = (b \lor c),\\
&(a \lor a) = (b \lor c),\ (a \land b) = (b \lor c),\ (a \lor b) = (b \lor c),\ (a \land b) = (c \land a),\ (a \lor b) = (c \land a),\\
&(a \land b) = (c \lor a),\ (a \lor b) = (c \lor a),\ (a \land b) = (c \land b),\ (a \lor b) = (c \land b),\ (a \land b) = (c \lor b),\\
&(a \lor b) = (c \lor b),\ (a \land b) = (c \land c),\ (a \lor b) = (c \land c),\ (a \land b) = (c \lor c),\ (a \lor b) = (c \lor c),\\
&\boxed{a = (\neg (\neg a))},\ (\neg a) = (\neg (\neg a)),\ \boxed{(a \land a) = (\neg (\neg a))},\ \boxed{(a \lor a) = (\neg (\neg a))},\ (a \land b) = (\neg (\neg a)),\\
&(a \lor b) = (\neg (\neg a)),\ a = (\neg (\neg b)),\ (\neg a) = (\neg (\neg b)),\ (a \land a) = (\neg (\neg b)),\ (a \lor a) = (\neg (\neg b)),\\
&(a \land b) = (\neg (\neg b)),\ (a \lor b) = (\neg (\neg b)),\ (\neg (\neg a)) = (\neg (\neg b)),\ (a \land b) = (\neg (\neg c)),\ (a \lor b) = (\neg (\neg c)),\\
&a = (\neg a \land a),\ (\neg a) = (\neg a \land a),\ (a \land a) = (\neg a \land a),\ (a \lor a) = (\neg a \land a),\ (a \land b) = (\neg a \land a),\\
&(a \lor b) = (\neg a \land a),\ (\neg (\neg a)) = (\neg a \land a),\ a = (a \land \neg a),\ (\neg a) = (a \land \neg a),\ (a \land a) = (a \land \neg a),\\
&(a \lor a) = (a \land \neg a),\ (a \land b) = (a \land \neg a),\ (a \lor b) = (a \land \neg a),\ (\neg (\neg a)) = (a \land \neg a),\ \boxed{(\neg a \land a) = (a \land \neg a)}\}
\end{aligned}
$$

OK, so now we can get a list of true theorems:

$$\begin{gathered}
\{a=(a\wedge a),\ a=(a\vee a),\ (a\wedge a)=(a\vee a),\ (a\wedge b)=(b\wedge a),\ (a\vee b)=(b\vee a),\ a=(\neg(\neg a)),\ (a\wedge a)=(\neg(\neg a)),\\
(a\vee a)=(\neg(\neg a)),\ (\neg a\wedge a)=(a\wedge\neg a),\ (\neg a\vee a)=(a\vee\neg a),\ (\neg a)=(\neg(a\wedge a)),\ (\neg a)=(\neg(a\vee a)),\\
(\neg(a\wedge a))=(\neg(a\vee a)),\ (a\wedge\neg b)=(\neg b\wedge a),\ (\neg a\wedge b)=(b\wedge\neg a),\ (a\vee\neg b)=(\neg b\vee a),\ (\neg a\vee b)=(b\vee\neg a),\\
(\neg(a\wedge b))=(\neg(b\wedge a)),\ (\neg(a\vee b))=(\neg(b\vee a)),\ (\neg a\wedge a)=(\neg b\wedge b),\ (a\wedge\neg a)=(\neg b\wedge b),\\
(\neg a\wedge a)=(b\wedge\neg b),\ (a\wedge\neg a)=(b\wedge\neg b),\ (\neg a\vee a)=(\neg b\vee b),\ (a\vee\neg a)=(\neg b\vee b),\ (\neg a\vee a)=(b\vee\neg b),\\
(a\vee\neg a)=(b\vee\neg b),\ (\neg a)=(\neg a\wedge\neg a),\ (\neg(a\wedge a))=(\neg a\wedge\neg a),\ (\neg(a\vee a))=(\neg a\wedge\neg a),\ (\neg a)=(\neg a\vee\neg a),\\
(\neg(a\wedge a))=(\neg a\vee\neg a),\ (\neg(a\vee a))=(\neg a\vee\neg a),\ (\neg a\wedge\neg a)=(\neg a\vee\neg a),\ (\neg(a\vee b))=(\neg a\wedge\neg b),\\
(\neg(a\wedge b))=(\neg a\vee\neg b),\ (\neg(a\vee b))=(\neg b\wedge\neg a),\ (\neg a\wedge\neg b)=(\neg b\wedge\neg a),\ (\neg(a\wedge b))=(\neg b\vee\neg a),\\
(\neg a\vee\neg b)=(\neg b\vee\neg a),\ a=((a\wedge a)\wedge a),\ (a\wedge a)=((a\wedge a)\wedge a),\ (a\vee a)=((a\wedge a)\wedge a),\\
(\neg(\neg a))=((a\wedge a)\wedge a),\ a=((a\vee a)\wedge a),\ (a\wedge a)=((a\vee a)\wedge a),\ (a\vee a)=((a\vee a)\wedge a),\ (\neg(\neg a))=((a\vee a)\wedge a),\\
((a\wedge a)\wedge a)=((a\vee a)\wedge a),\ a=(a\wedge(a\wedge a)),\ (a\wedge a)=(a\wedge(a\wedge a)),\ (a\vee a)=(a\wedge(a\wedge a)),\\
(\neg(\neg a))=(a\wedge(a\wedge a)),\ ((a\wedge a)\wedge a)=(a\wedge(a\wedge a)),\ ((a\vee a)\wedge a)=(a\wedge(a\wedge a)),\ a=(a\wedge(a\vee a)),\\
(a\wedge a)=(a\wedge(a\vee a)),\ (a\vee a)=(a\wedge(a\vee a)),\ (\neg(\neg a))=(a\wedge(a\vee a)),\ ((a\wedge a)\wedge a)=(a\wedge(a\vee a))\}
\end{gathered}$$

Some are "interesting". Others seem repetitive, overly complicated, or otherwise not terribly interesting. But if we want to "channel Euclid" we somehow have to decide which are the interesting theorems that we're going to write down. And although Euclid himself didn't explicitly discuss logic, we can look at textbooks of logic from the last couple of centuries—and we find that there's a very consistent set of theorems that they end up picking out from the list, and giving names to:

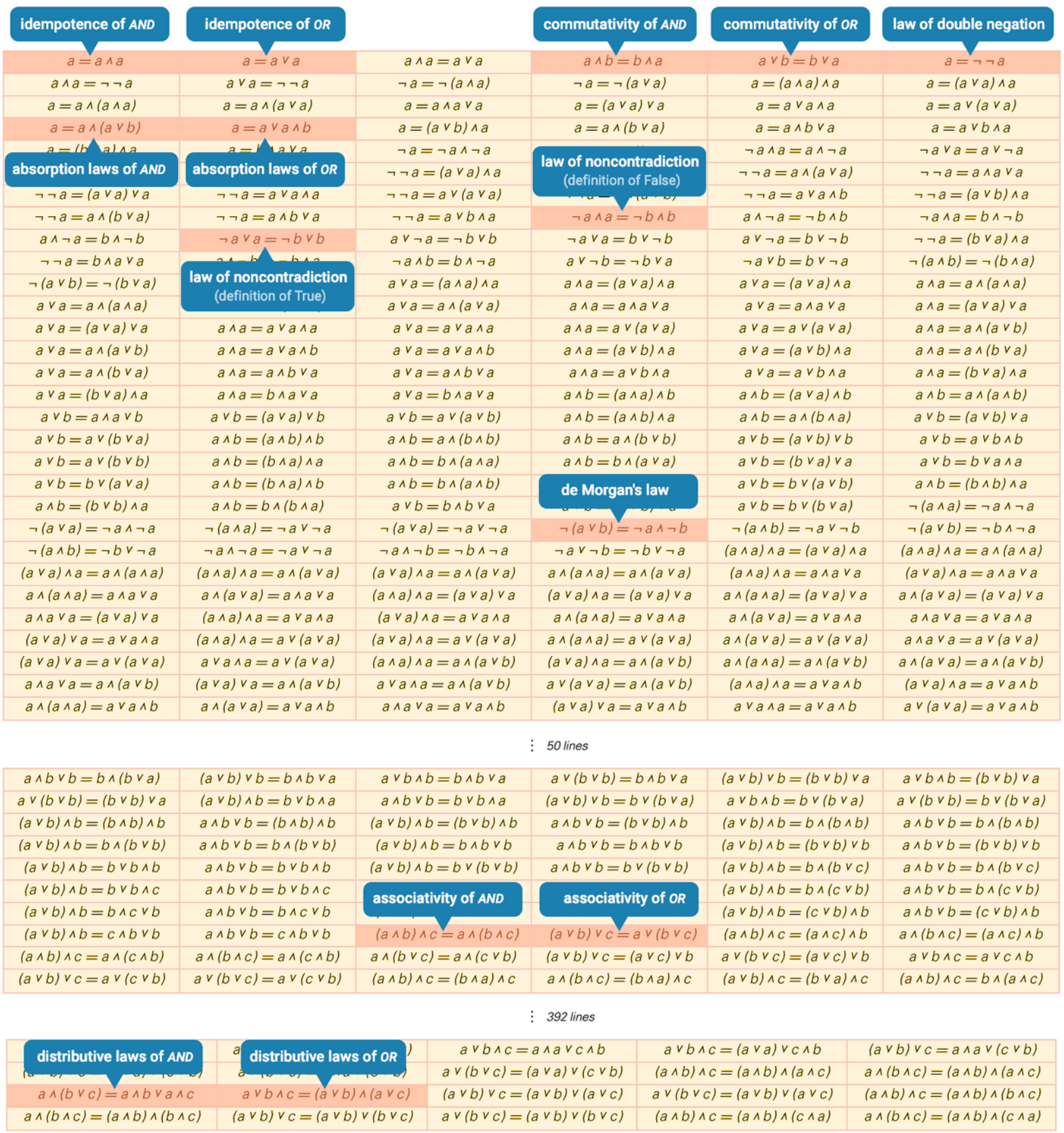

One might assume that these named theorems were just the result of historical convention. But when I was writing *A New Kind of Science* I discovered something quite surprising. With all the theorems written out in "order of complexity", I tried seeing which theorems I could prove just from theorems earlier in the list. Many were easy to prove. But some simply couldn't be proved. And it turned out that these were essentially precisely the "named theorems":

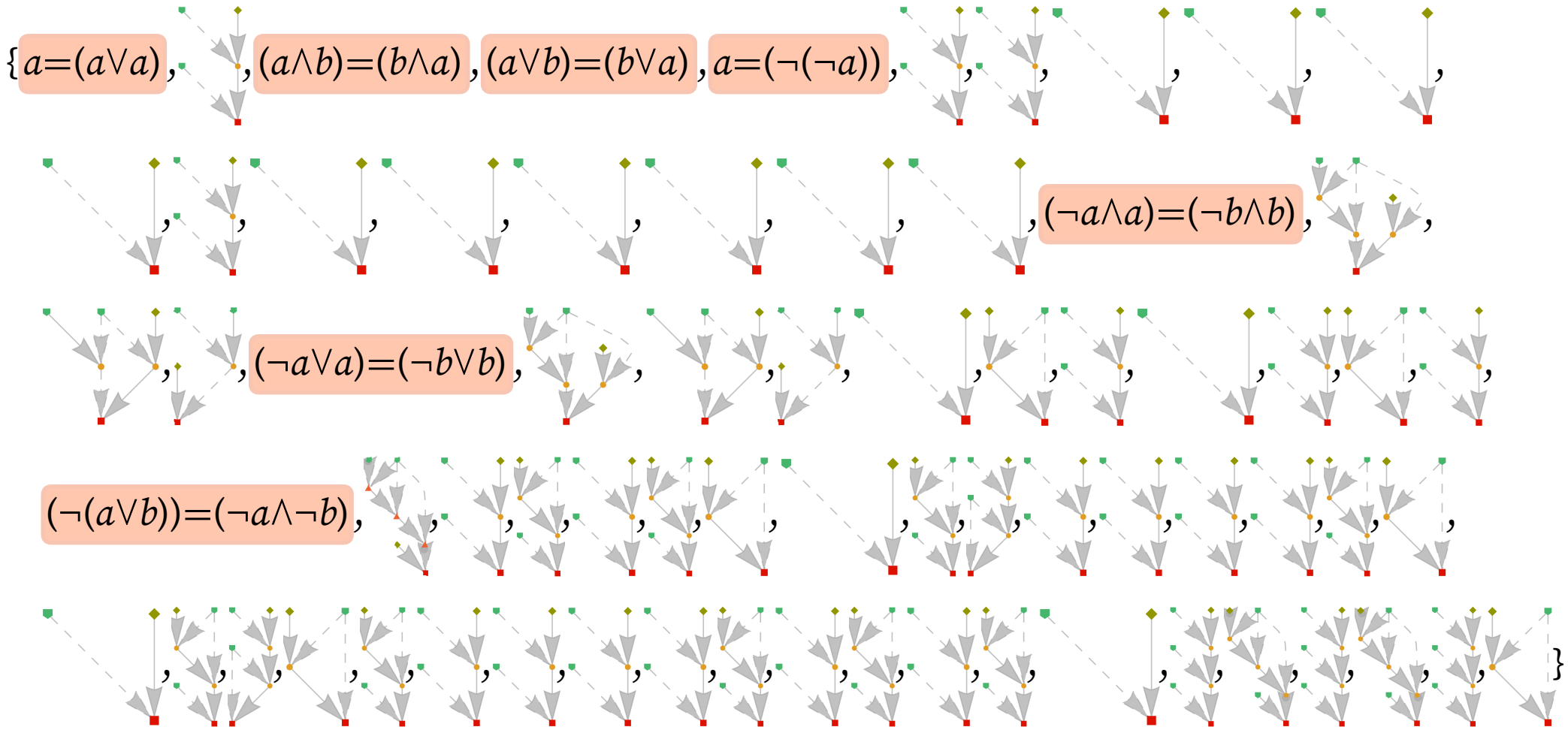

In other words, the "named theorems" are basically the simplest statements of new facts about logic, that can't be established from "simpler facts". Eventually as one's going through the list of theorems, one will have accumulated enough to fill out what can serve as full axioms for logic—so that then all subsequent theorems can be proved from "existing facts".

Now of course the setup we've just used relies on the idea that one's separately got a list of true theorems. To do something more like Euclid, we'd have to pick certain theorems to serve as axioms, then derive all others from these.

Back in 2000 I figured out the very simplest possible axiom system for logic, written in terms of Nand, just the single axiom:

$$\{\forall_{\{a,b,c\}}((b\cdot c)\cdot a)\cdot(b\cdot((b\cdot a)\cdot b))=a\}$$

So now writing And, Or and Not in terms of Nand according to

$$\{(\neg a)=a\cdot a,(a\wedge b)=(a\cdot b)\cdot(a\cdot b),(a\vee b)=(a\cdot a)\cdot(b\cdot b)\}$$

we can, for example, derive the notable theorems of logic from my axiom. FindEquationalProof gives automated proofs of these theorems, though most of them involve quite a few steps (the — indicates a theorem that is trivially true after substituting the forms for And, Or and Not):

| | |
|---|---|
| $a=(a\wedge a)$ | 54 |
| $a=(a\vee a)$ | 54 |
| $(a\wedge b)=(b\wedge a)$ | 103 |
| $(a\vee b)=(b\vee a)$ | 102 |
| $(\neg(\neg a))=a$ | 54 |
| $(\neg a\wedge a)=(\neg b\wedge b)$ | 95 |
| $(\neg a\vee a)=(\neg b\vee b)$ | 92 |
| $(\neg(a\vee b))=(\neg a\wedge\neg b)$ | 132 |
| $(\neg(a\wedge b))=(\neg a\vee\neg b)$ | 143 |
| $(a\wedge b)=(a\wedge(a\wedge b))$ | 91 |
| $a=(a\wedge(a\vee b))$ | 328 |
| $a=(a\vee(a\wedge b))$ | 274 |
| $(\neg a\vee b)=(\neg(\neg b)\vee\neg a)$ | 958 |
| $((a\wedge b)\wedge c)=(a\wedge(b\wedge c))$ | 1502 |
| $((a\vee b)\vee c)=(a\vee(b\vee c))$ | |
| $(a\vee b)=(a\vee(\neg a\wedge b))$ | 56 |
| $a=(a\wedge(b\vee\neg b))$ | 131 |
| $a=(a\vee(b\wedge\neg b))$ | 130 |
| $(a\vee(b\wedge c))=((a\vee b)\wedge(a\vee c))$ | 120 |
| $(a\wedge(b\vee c))=((a\wedge b)\vee(a\wedge c))$ | 103 |

The longer cases here involve first proving the lemma $a{\cdot}b{=}b{\cdot}a$ which takes 102 steps. Including this lemma as an axiom, the minimal axiom system (as I also found in 2000) is:

$$\{\forall_{\{a,b,c\}}(a{\cdot}b){\cdot}(a{\cdot}(b{\cdot}c))=a,\forall_{\{a,b\}}a{\cdot}b=b{\cdot}a\}$$

And with this axiom system **FindEquationalProof** succeeds in finding shorter proofs for the notable theorems of logic, even though now the definitions for **And**, **Or** and **Not** are just treated as theorems:

| | |
|---|---|
| $a=(a\wedge a)$ | 21 |
| $a=(a\vee a)$ | 15 |
| $(a\wedge b)=(b\wedge a)$ | 8 |
| $(a\vee b)=(b\vee a)$ | 9 |
| $(\neg(\neg a))=a$ | 17 |
| $(\neg a\wedge a)=(\neg b\wedge b)$ | 130 |
| $(\neg a\vee a)=(\neg b\vee b)$ | 119 |
| $(\neg(a\vee b))=(\neg a\wedge\neg b)$ | 9 |
| $(\neg(a\wedge b))=(\neg a\vee\neg b)$ | 28 |
| $(a\wedge b)=(a\wedge(a\wedge b))$ | 43 |
| $a=(a\wedge(a\vee b))$ | 32 |
| $a=(a\vee(a\wedge b))$ | 26 |
| $(\neg a\vee b)=(\neg(\neg b)\vee\neg a)$ | 20 |
| $((a\wedge b)\wedge c)=(a\wedge(b\wedge c))$ | 249 |
| $((a\vee b)\vee c)=(a\vee(b\vee c))$ | 239 |
| $(a\vee b)=(a\vee(\neg a\wedge b))$ | 89 |
| $a=(a\wedge(b\vee\neg b))$ | 129 |
| $a=(a\vee(b\wedge\neg b))$ | 129 |
| $(a\vee(b\wedge c))=((a\vee b)\wedge(a\vee c))$ | 328 |
| $(a\wedge(b\vee c))=((a\wedge b)\vee(a\wedge c))$ | 338 |

Actually looking at these proofs is not terribly illuminating; they certainly don't have the same kind of "explanatory feel" as Euclid. But combining the graphs for all these proofs is more interesting, because it shows us the common lemmas that were used in these proofs, and effectively defines a network of interdependencies between theorems:

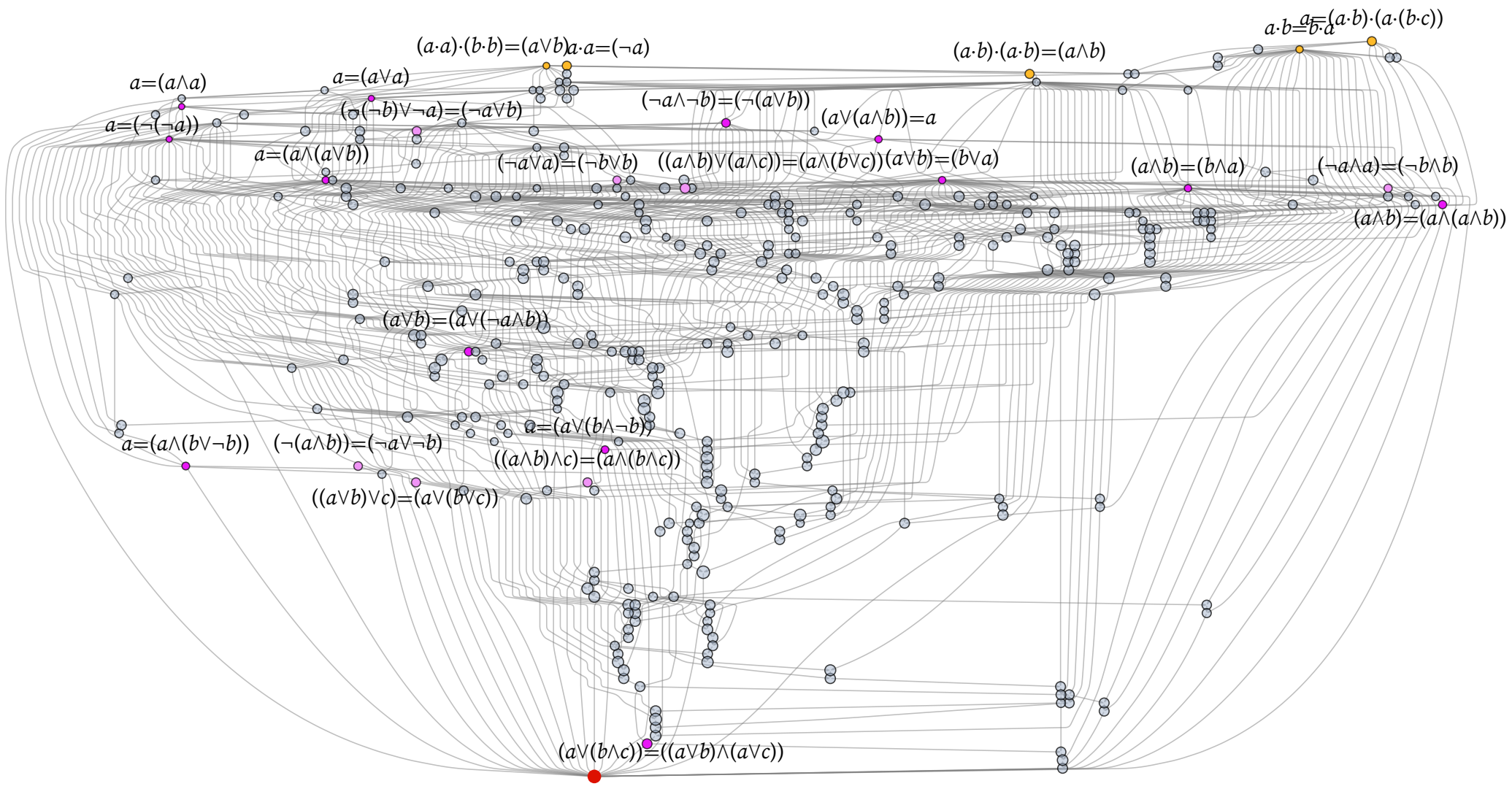

There are 361 lemmas (i.e. automatically generated intermediate theorems) here. It's a fair number, given that we're only proving 20 theorems—but it's definitely much less than the total of 1978 that would be involving in proving each of the theorems separately.

In our graph here—like in our Euclid theorem-dependency graphs above—the axioms are shown (in yellow) at the top. The "notable theorems" that we're proving are shown in pink. But the structure of the graph is a little different from our earlier Euclid theorem-dependency graphs, and this alternative layout makes it clearer:

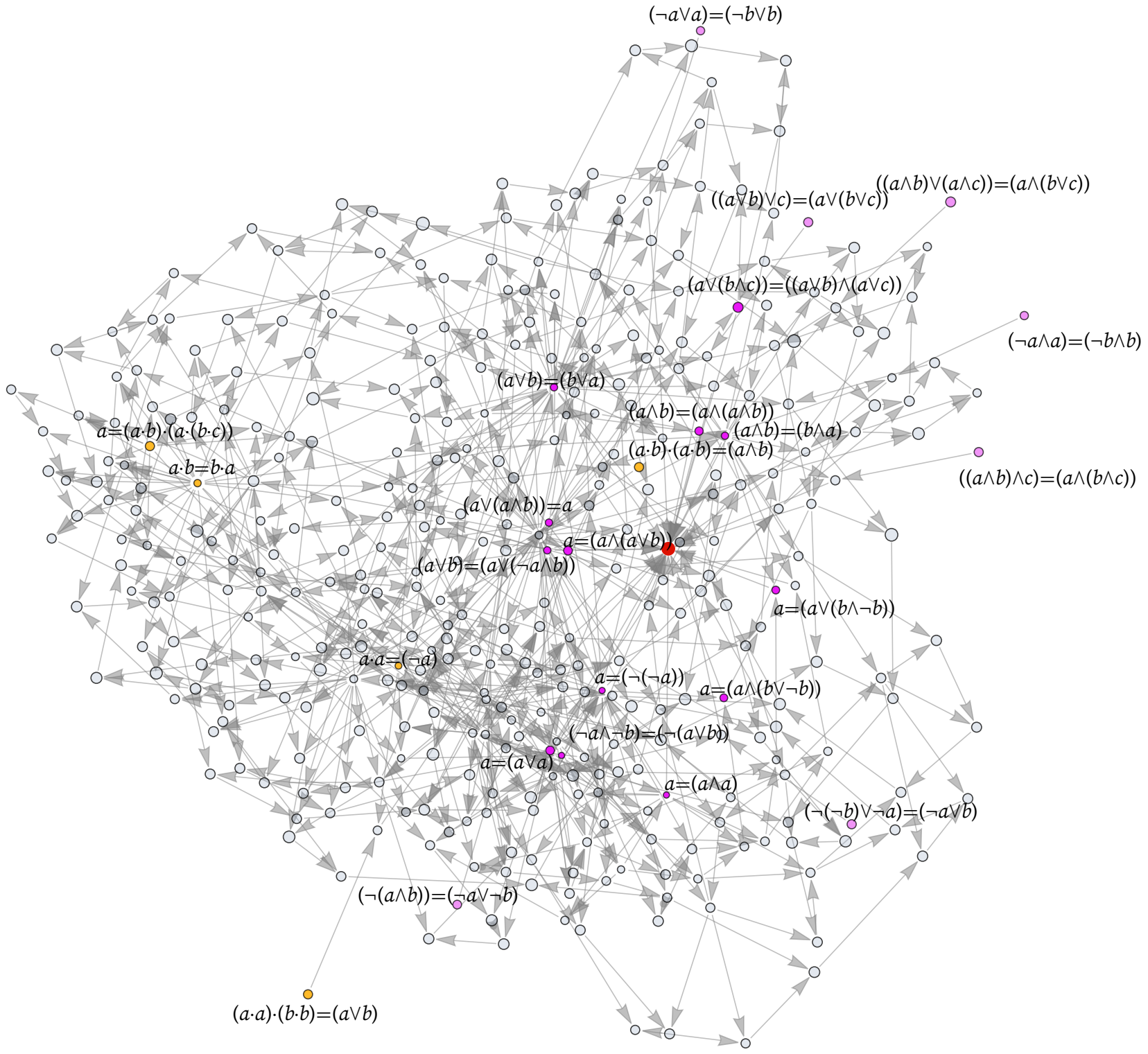

In Euclid, a given theorem is proved on the basis of other theorems, and ultimately on the basis of axioms. But here the automated theorem-proving process creates lemmas that ultimately allow one to show that the theorems one's trying to prove are equivalent to "true" (i.e. to a tautology)—shown as a red node.

We can ask other questions, such as how long the lemmas are. Here are the distributions of lengths of the final notable theorems, and of the intermediate lemmas used to prove them:

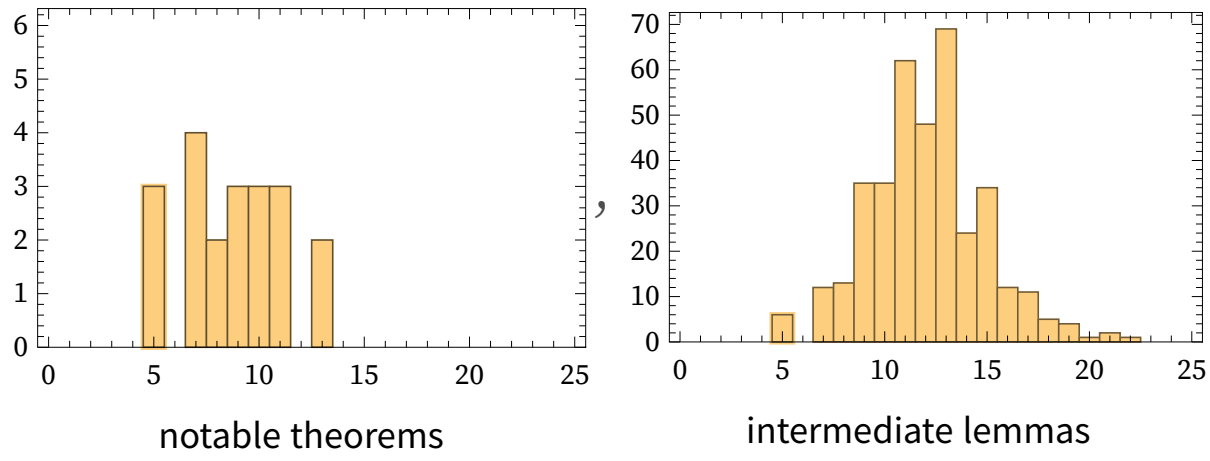

We get something slightly more in the spirit of Euclid if we elide the lemmas, and just find the implied effective dependency graph between notable theorems:

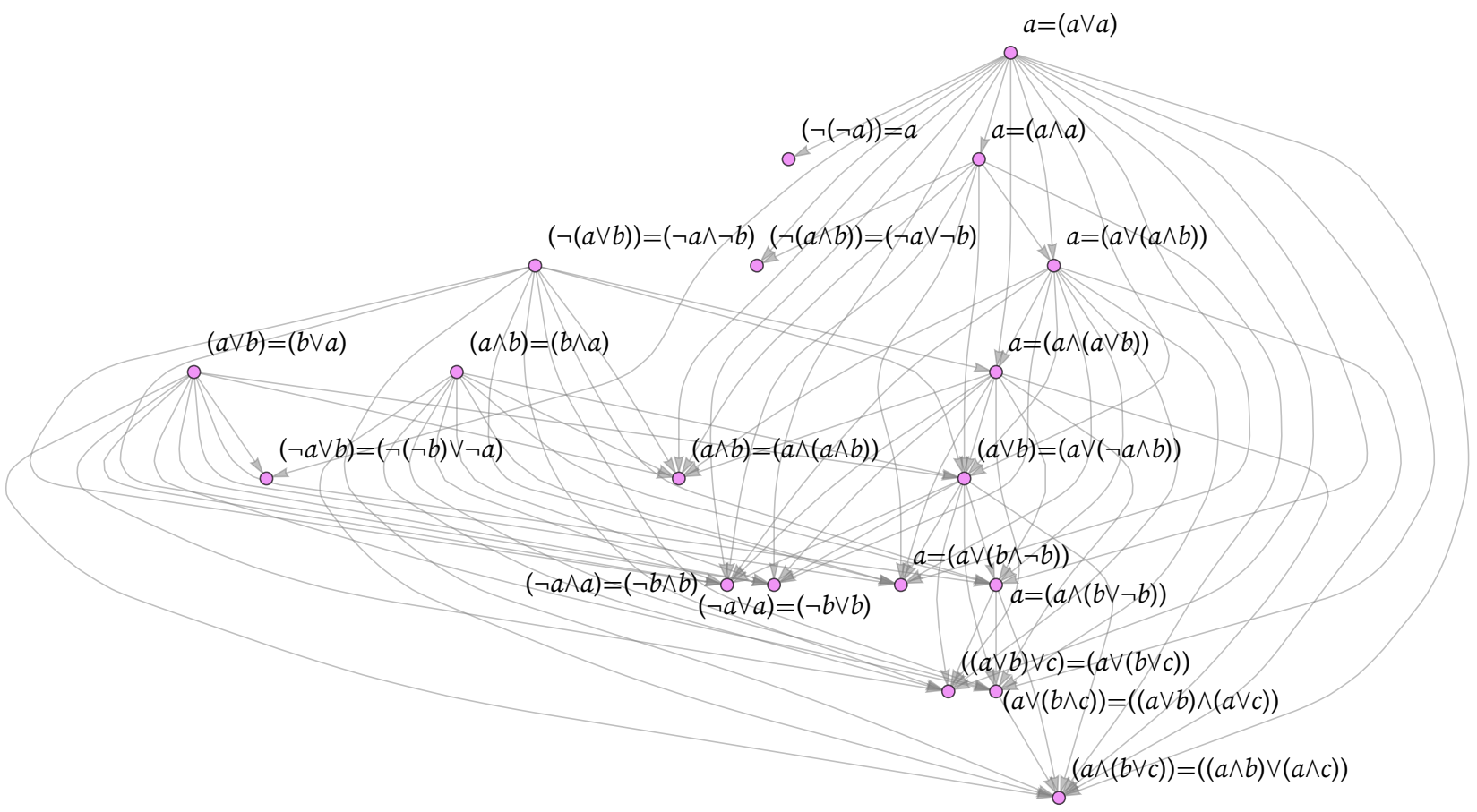

Transitive reduction then gives:

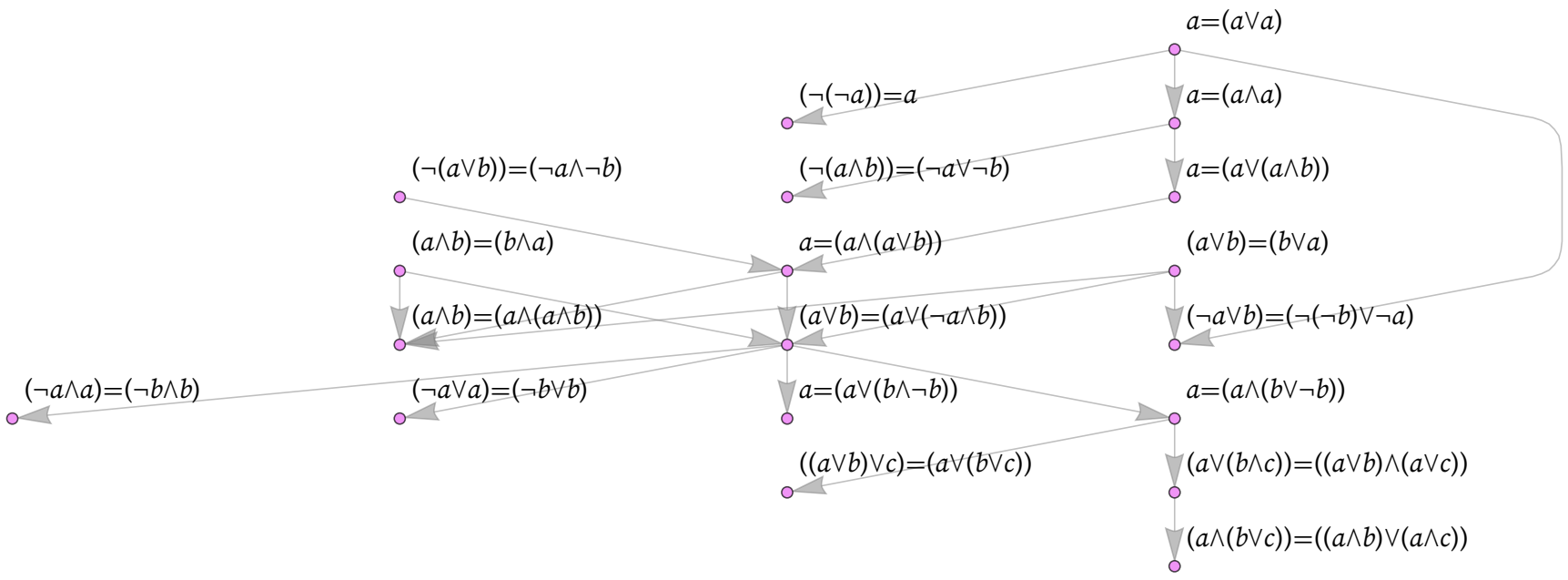

By omitting intermediate lemmas, we're in a sense just getting a shadow of the dependencies of the notable theorems, in the "environment" defined by our particular choice of axioms. But with this setup, it's interesting to see the distributive law be the "hardest theorem"—kind of the metamathematical analog of Euclid's 13.18 about the Platonic solids.

OK, but what we're doing so far with logic is still fundamentally a bit different from how most of Euclid works. Because what Euclid typically does is to say something like "imagine such-and-such a geometrical setup; then the following theorem will be true about it". And the analog of that for logic would be to take axioms of logic, then append some logical assertion, and ask if with the axioms and this assertion some particular statement is true. In other words, there are some statements—like the axioms—that will be true in "pure logic", but there are more statements that will be true with particular setups (or, in the case of logic, particular possible values for variables).

For example, in "pure logic" $a \lor b == b \lor b$ is not necessarily true (i.e. it is not a tautology). But if we assert that $a == (a \land b)$ is true, then this implies the following possible choices for $a$ and $b$

```
In[•]:= SatisfiabilityInstances[a == (a ∧ b), {a, b}, All]
Out[•]= {{True, True}, {False, True}, {False, False}}
```

and in all these cases $a\lor b$==$b\lor b$ is true. So, in a Euclid tradition, we could say "imagine a setup where $a$==($a\land b$); then we can prove from the axioms of logic the theorem that $a$==($a\land b$)".

Above we looked at which statements in logic are true for all values of variables:

$\{a=b,\ a=(\neg a),\ a=(\neg b),\ (\neg a)=(\neg b),\ a=(a\land a),\ (\neg a)=(a\land a),\ a=(a\lor a),\ (\neg a)=(a\lor a),\ (a\land a)=(a\lor a),\ a=(a\land b),$
$(\neg a)=(a\land b),\ (a\land a)=(a\land b),\ (a\lor a)=(a\land b),\ a=(a\lor b),\ (\neg a)=(a\lor b),\ (a\land a)=(a\lor b),\ (a\lor a)=(a\lor b),\ (a\land b)=(a\lor b),$
$(a\land b)=(a\land c),\ (a\lor b)=(a\land c),\ (a\land b)=(a\lor c),\ (a\lor b)=(a\lor c),\ a=(b\land a),\ (\neg a)=(b\land a),\ (a\land a)=(b\land a),\ (a\lor a)=(b\land a),$
$(a\land b)=(b\land a),\ (a\lor b)=(b\land a),\ a=(b\lor a),\ (\neg a)=(b\lor a),\ (a\land a)=(b\lor a),\ (a\lor a)=(b\lor a),\ (a\land b)=(b\lor a),\ (a\lor b)=(b\lor a),$
$a=(b\land b),\ (\neg a)=(b\land b),\ (a\land a)=(b\land b),\ (a\lor a)=(b\land b),\ (a\land b)=(b\land b),\ (a\lor b)=(b\land b),\ a=(b\lor b),\ (\neg a)=(b\lor b),$
$(a\land a)=(b\lor b),\ (a\lor a)=(b\lor b),\ (a\land b)=(b\lor b),\ (a\lor b)=(b\lor b),\ a=(b\land c),\ (\neg a)=(b\land c),\ (a\land a)=(b\land c),\ (a\lor a)=(b\land c)\}$

Now let's look at the ones that aren't always true. If we assume that some particular one of these statements is true, we can see which other statements it implies are true:

| | |
|---|---|
| $a = b$ | $\{a = b, a = (a\lor b), a = (b\lor a), a = (a\land b), a = (b\land a), (\neg a) = (\neg b), (a\lor a) = (a\lor b), (a\lor a) = (a\land b), (a\lor a) = (b\land a), (a\lor b) = (b\land a), (a\land a) = (a\lor b), (a\land a) = (a\land b), (a\land a) = (b\land a), (a\land b) = (a\lor b)\}$ |
| $a = (\neg a)$ | $\{\}$ |
| $a = (\neg b)$ | $\{a = (\neg b)\}$ |
| $a = (a\lor b)$ | $\{a = (a\lor b), a = (b\lor a), (a\lor a) = (a\lor b), (a\land a) = (a\lor b)\}$ |
| $a = (b\lor a)$ | $\{a = (a\lor b), a = (b\lor a), (a\lor a) = (a\lor b), (a\land a) = (a\lor b)\}$ |
| $a = (a\land b)$ | $\{a = (a\land b), a = (b\land a), (a\lor a) = (a\land b), (a\lor a) = (b\land a), (a\land a) = (a\land b), (a\land a) = (b\land a)\}$ |
| $a = (b\land a)$ | $\{a = (a\land b), a = (b\land a), (a\lor a) = (a\land b), (a\lor a) = (b\land a), (a\land a) = (a\land b), (a\land a) = (b\land a)\}$ |
| $(\neg a) = (\neg b)$ | $\{a = b, a = (a\lor b), a = (b\lor a), a = (a\land b), a = (b\land a), (\neg a) = (\neg b), (a\lor a) = (a\lor b), (a\lor a) = (a\land b), (a\lor a) = (b\land a), (a\lor b) = (b\land a), (a\land a) = (a\lor b), (a\land a) = (a\land b), (a\land a) = (b\land a), (a\land b) = (a\lor b)\}$ |
| $(\neg a) = (a\lor a)$ | $\{\}$ |
| $(\neg a) = (a\lor b)$ | $\{a = (\neg b), a = (a\land b), a = (b\land a), (\neg a) = (a\lor b), (\neg a) = (b\lor a), (a\lor a) = (a\land b), (a\lor a) = (b\land a), (a\land a) = (a\land b), (a\land a) = (b\land a), (a\land b) = (a\land c)\}$ |
| $(\neg a) = (b\lor a)$ | $\{a = (\neg b), a = (a\land b), a = (b\land a), (\neg a) = (a\lor b), (\neg a) = (b\lor a), (a\lor a) = (a\land b), (a\lor a) = (b\land a), (a\land a) = (a\land b), (a\land a) = (b\land a), (a\land b) = (a\land c)\}$ |
| $(\neg a) = (a\land a)$ | $\{\}$ |
| $(\neg a) = (a\land b)$ | $\{a = (\neg b), a = (a\lor b), a = (b\lor a), (\neg a) = (a\land b), (\neg a) = (b\land a), (a\lor a) = (a\lor b), (a\lor b) = (a\lor c), (a\land a) = (a\lor b)\}$ |
| $(\neg a) = (b\land a)$ | $\{a = (\neg b), a = (a\lor b), a = (b\lor a), (\neg a) = (a\land b), (\neg a) = (b\land a), (a\lor a) = (a\lor b), (a\lor b) = (a\lor c), (a\land a) = (a\lor b)\}$ |
| $(a\lor a) = (a\lor b)$ | $\{a = (a\lor b), a = (b\lor a), (a\lor a) = (a\lor b), (a\land a) = (a\lor b)\}$ |
| $(a\lor a) = (a\land b)$ | $\{a = (a\land b), a = (b\land a), (a\lor a) = (a\land b), (a\lor a) = (b\land a), (a\land a) = (a\land b), (a\land a) = (b\land a)\}$ |
| $(a\lor a) = (b\land a)$ | $\{a = (a\land b), a = (b\land a), (a\lor a) = (a\land b), (a\lor a) = (b\land a), (a\land a) = (a\land b), (a\land a) = (b\land a)\}$ |
| $(a\lor b) = (a\lor c)$ | $\{(a\lor b) = (a\lor c)\}$ |

Or on a larger scale, with a black dot when one statement implies another:

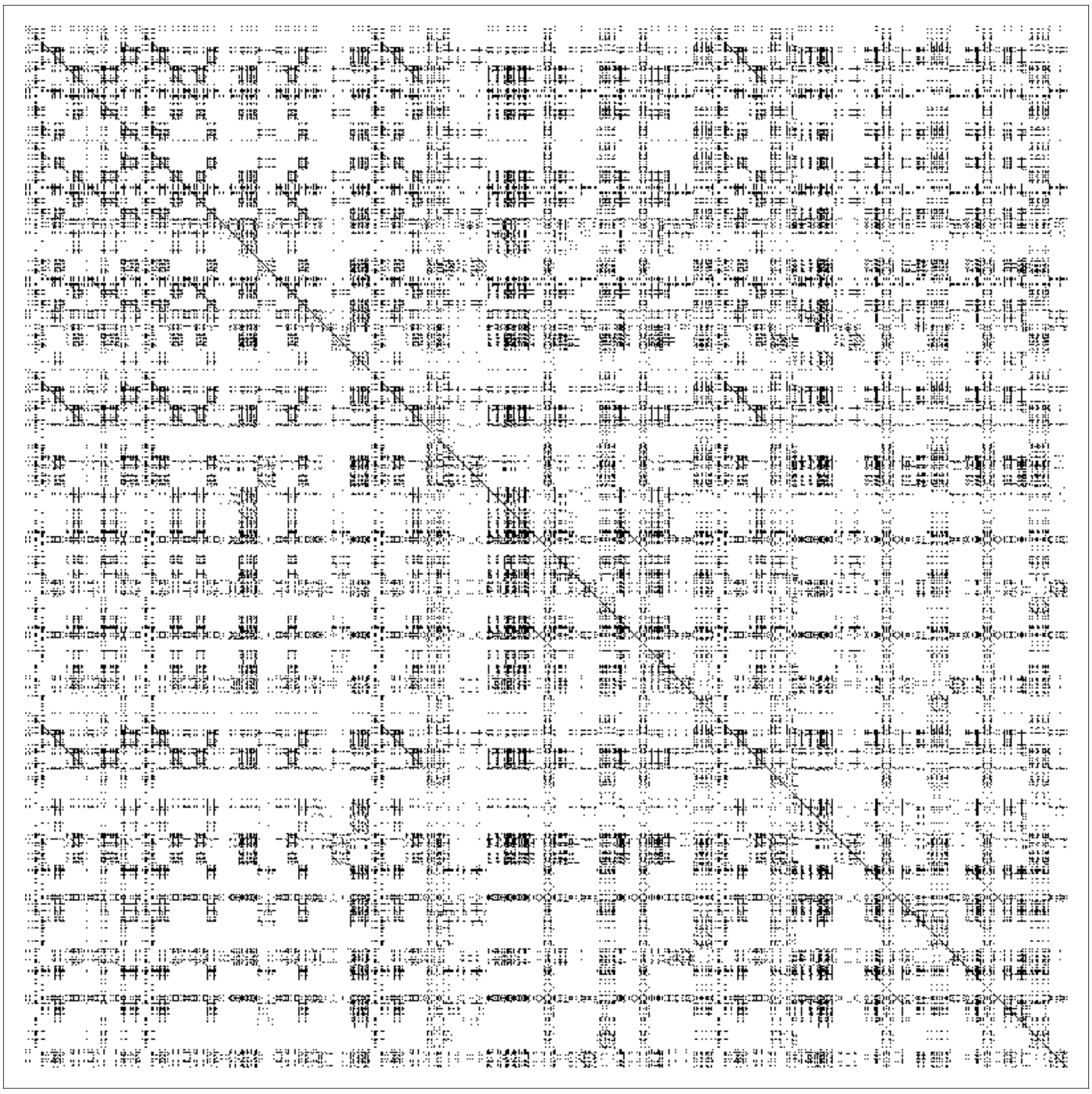

For each of these theorems we can in principle construct a proof, using the axioms:

```
In[•]:= FindEquationalProof[(a ⊕ b) == (b ⊕ b), Append[AxiomaticTheory["BooleanAxioms"], a == (a ⊗ b)],
          "ProofGraph"]
```

Out[•]=

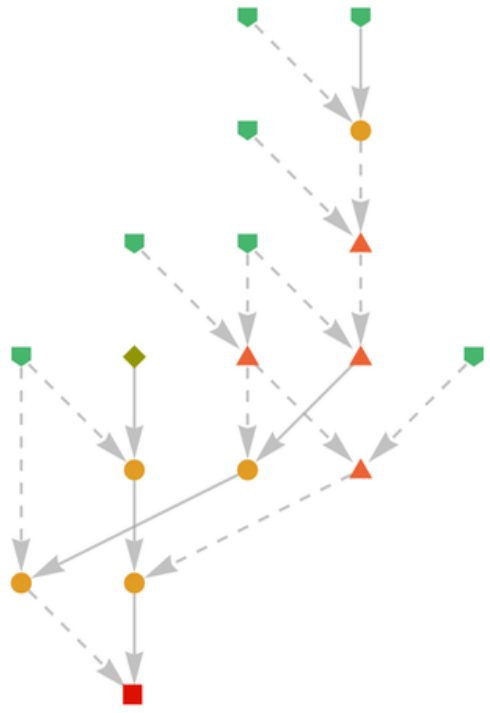

And now we could go through and find out which theorems are useful in proving other theorems—and in principle this would allow us to build up a theorem dependency network. But there are undoubtedly many ways to do this, and so we'd need additional criteria to find ones that have whatever attributes would make us say "that might have been how someone like Euclid would have done it".

OK, so could one look at geometry the same way? Basically, yes. Using the formalization we had above in terms of line, between, congruent, etc. we can again start by just enumerating possible statements. Unlike for logic, many of them won't even make "structural sense"; for example they might contain line[congruent[...],...], but it makes no sense to have a line whose endpoint is a truth value. But we can certainly get a list of "structurally meaningful" statements.

And then we can ask which are "tautologically true"—though it's in practice considerably harder to do this than for logic (the best known methods involve all sorts of elaborate algebraic computations, which Mathematica can certainly do, but which quickly become quite unwieldy). And after that, we can proceed like Euclid, and start saying "assert this, then you can prove this". And, yes, it's nice that after 2000+ years, we can finally imagine automating the process of producing generalizations of Euclid's *Elements*. Though this just makes it more obvious that part of what Euclid did was in a sense a matter of art—picking in some kind of aesthetic way which possible sequence of theorems would best "tell his story" of geometry.

## Math beyond Euclid

We've looked here at some of the empirical metamathematics of what Euclid did on geometry more than 2000 years ago. But what about more recent mathematics, and all those other areas of mathematics that have now been studied? In the history of mathematics, there have been perhaps 5 million research papers published, as well as probably hundreds of thousands of textbooks (though few quite as systematic as Euclid).

And, yes, in modern times almost all mathematics that's published is on the web in some form. A few years ago we scraped arXiv and identified about 2 million things described as theorems there (the most popular being the central limit theorem, the implicit function theorem and Fubini's theorem); we also scraped as much as we could of the visible web and found about 30 million theorems there. No doubt many were duplicates (though it's hard—and in principle undecidable!—which they are). But it's a reasonable estimate that there are a few million distinct theorems for which proofs have been published in the history of human mathematics.

It's a remarkable piece of encapsulated intellectual achievement—perhaps the largest coherent such one produced by our species. And I've long been interested in seeing just what it would take to make it computable, and to bring it into the whole computational knowledge framework we have in the Wolfram Language. A few years ago I hoped that we could mobilize the mathematics community to help make this happen. But formalization is hard work, and it's not at the center of what most mathematicians aspire to. Still, we've at least been slowly working—much as we have for Euclid-style geometry—to define the elements of computational language needed to represent theorems in various areas of mathematics.

For example, in the area of point-set topology, we have under development things like

```
In[∘]:= is Hausdorff TOPOLOGY CONCEPT ["Output"] // InputForm

Out[∘]= ForAll[{x, y}, Element[x, 𝒳["Elements"]] && Element[y, 𝒳["Elements"]] && x != y,
 Exists[{U, V}, Element[U, 𝒳["Topology"]] && Element[V, 𝒳["Topology"]] && SetIntersection[U, V] == EmptySet &&
  Element[x, U] && Element[y, V]]]
```

which in traditional mathematical notation becomes:

$$\forall_{\{x,y\},x\in\mathcal{X}\land y\in\mathcal{X}\land x\neq y}\exists_{\{U,V\}}(U\in\tau_{\mathcal{X}}\land V\in\tau_{\mathcal{X}}\land U\bigcap V=\emptyset\land x\in U\land y\in V)$$

So far we have encoded in computable form 742 "topology concepts", and 1687 theorems about them. Here are the connections recorded between concepts (dropping the concept of "topological spaces" that a third of all concepts are connected to, and labeling concepts with high betweenness centrality):

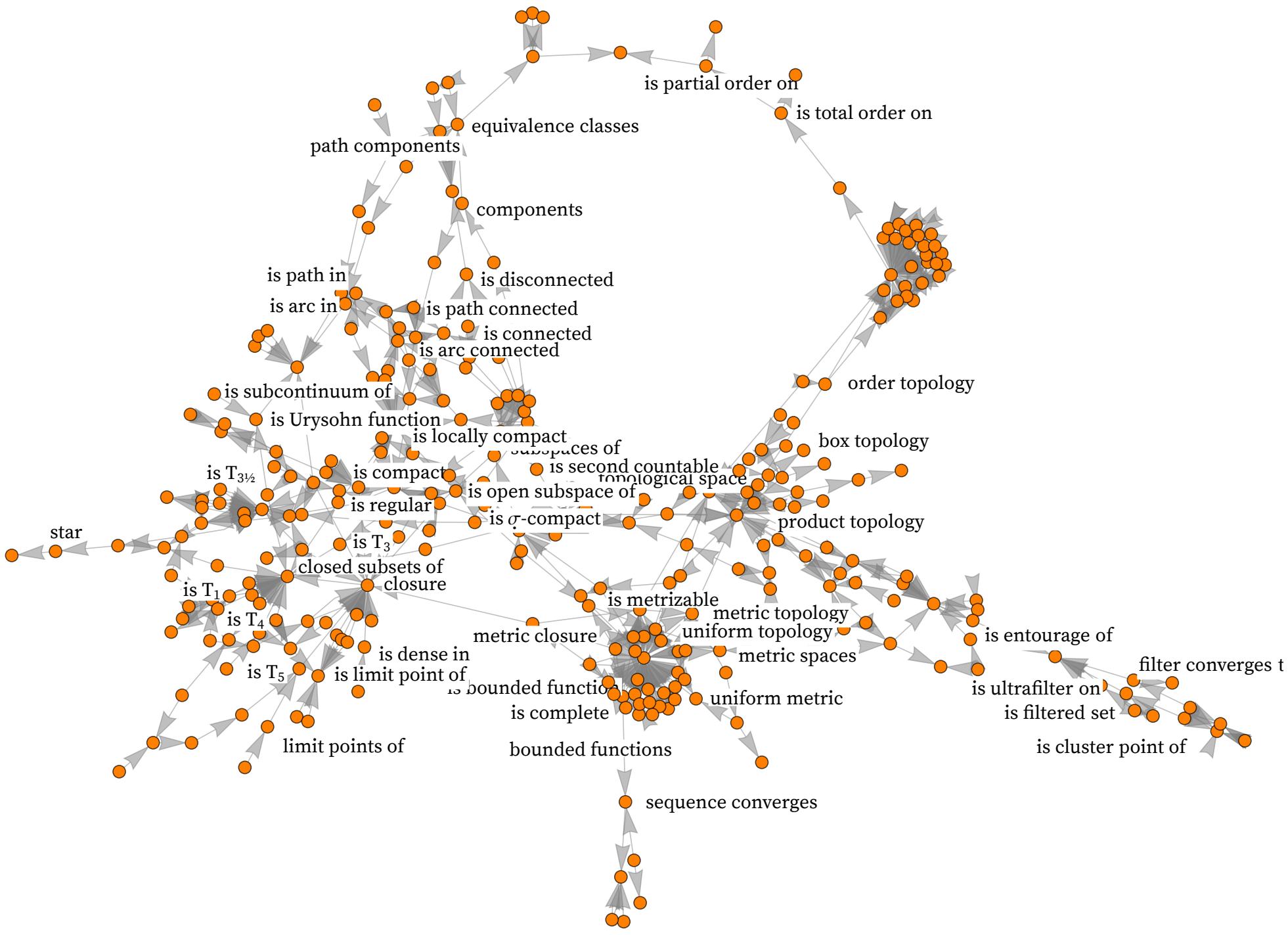

And here is the graph of what theorem references what in its description:

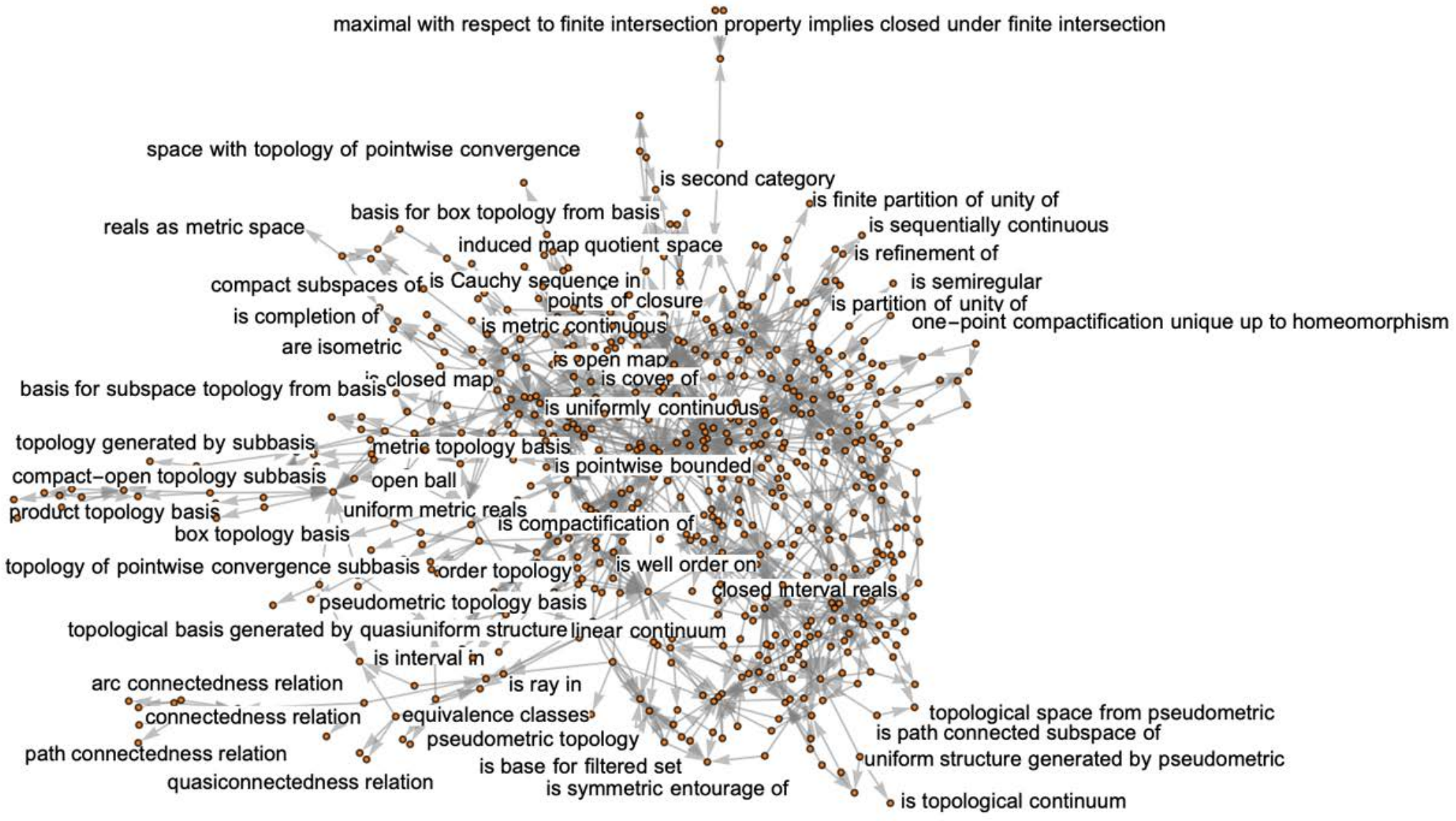

We haven't encoded proofs for these theorems, so we can't yet make the kind of theorem dependency graph that we did for Euclid. But we do have the dependency graph for 76 properties of topological spaces:

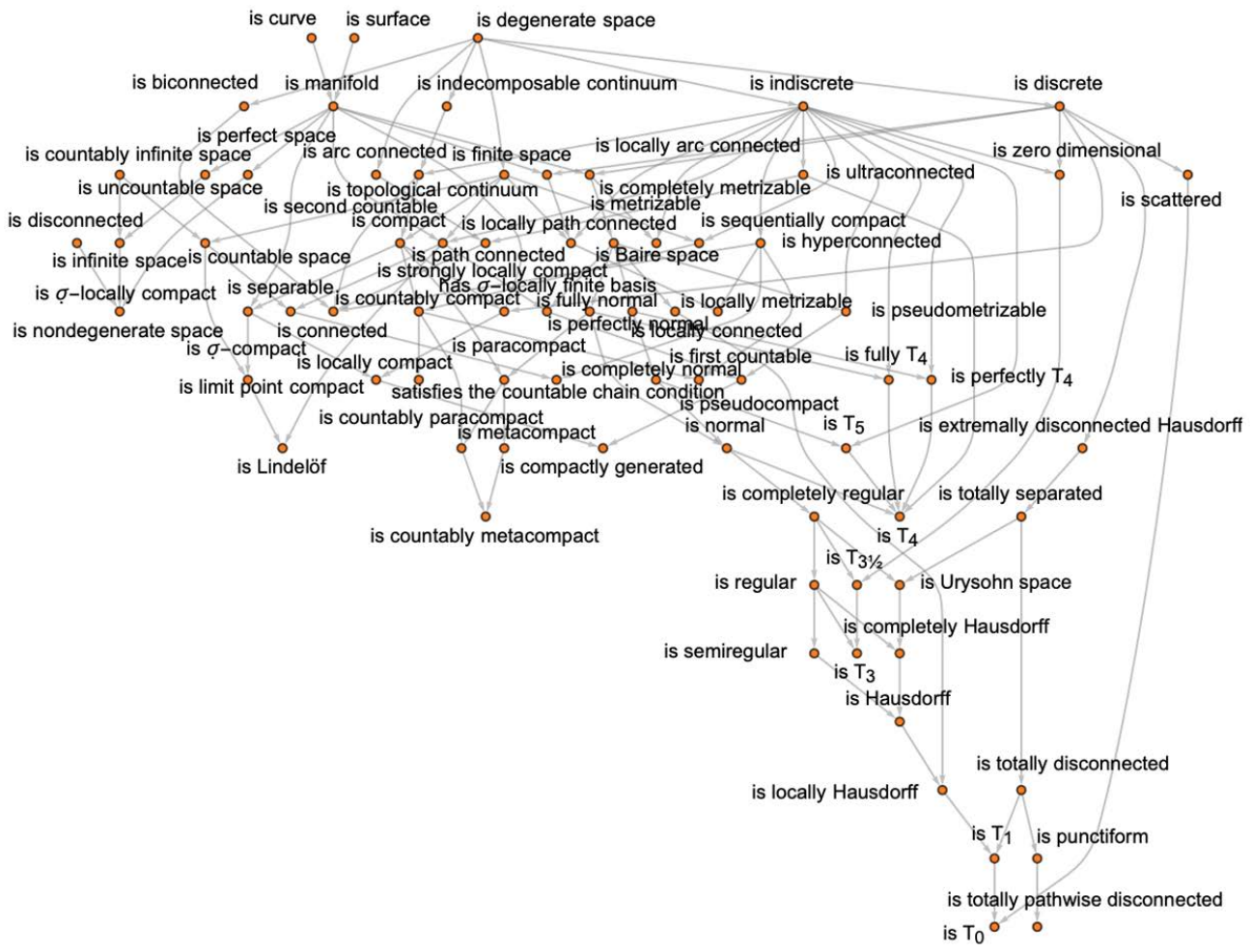

The longest path here (along with a similar one starting with *is curve*) is 14 steps:

*is surface* → *is manifold* → *is completely metrizable* → *is metrizable* →
*is perfectly normal* → *is completely normal* → *is normal* → *is completely regular* →
*is regular* → *is semiregular* → *is Hausdorff* → *is locally Hausdorff* → *is $T_1$* → *is $T_0$*

(And, yes, this isn't particularly profound; it's just an indication of what it looks like to make specific definitions in topology computable.)

So far, what we've discussed is being able to represent pure mathematical ideas and results in a high-level computable way, understandable to both humans and computers. But what if we want to just formalize everything, from the ground up, explicitly deriving and validating every theorem from the lowest-level foundations? Over the past few decades there have been a number of large-scale projects—like Mizar, Coq, Isabelle, HOL, Metamath, Lean—that have tried to do this (nowadays often in connection with creating "proof assistants").

Ultimately each project defines a certain "machine code" for mathematics. And yes, even though people might think that "mathematics is a universal language", if one's really going to give full, precise, formal specifications there are all sorts of choices to be made. Should things be based on set theory, type theory, higher-order logic, calculus of constructions, etc.? Should the law of excluded middle be assumed? The axiom of choice? What if one's axiomatic structure seems great, but implies a few silly results, 1/0 = 0? There's no perfect solution, but each of these projects has made a certain set of choices.

And the good news here is that for our purposes in doing large-scale empirical metamathematics—as in doing mathematics in the way mathematicians usually do it—it doesn't seem like the choices will matter much. But what's important for us is that these projects have accumulated tens of thousands of theorems (well, OK, some are "throwaway lemmas" or simple rearrangements), and that starting from axioms (or what amount to axioms), they've reached decently far into quite a few areas of mathematics.

Looking at them is a bit of a different experience from looking at Euclid. While the *Elements* has the feel of a "narrative textbook" (albeit from a different age), formalized mathematics projects tend to seem more like software codebases, with their theorem dependency graphs being more like function call graphs. But they still provide fascinating metamathematical corpuses, and there's undoubtedly lots about empirical metamathematics that one can learn from them.

Here I'm going to look at two examples: the Lean mathlib collection, which includes about 36,000 theorems (and 16,000 definitions) and the Metamath set.mm ("set theory") collection, which has about 44,000 theorems (and 1500 definitions).

To get a sense of what's in these collections, we can start by drawing interdependence graphs for the theorems they contain in different areas of mathematics. Just like for Euclid, we make the size of each node represent the number of theorems in a particular area, and the thickness of each edge represent the fraction of theorems from one area that directly reference another in their proof.

Leaving out theorems that effectively just do structural manipulation, rather than representing mathematical content (as well as "self-loop" connections within a single domain) here's the interdependence graph for Lean:

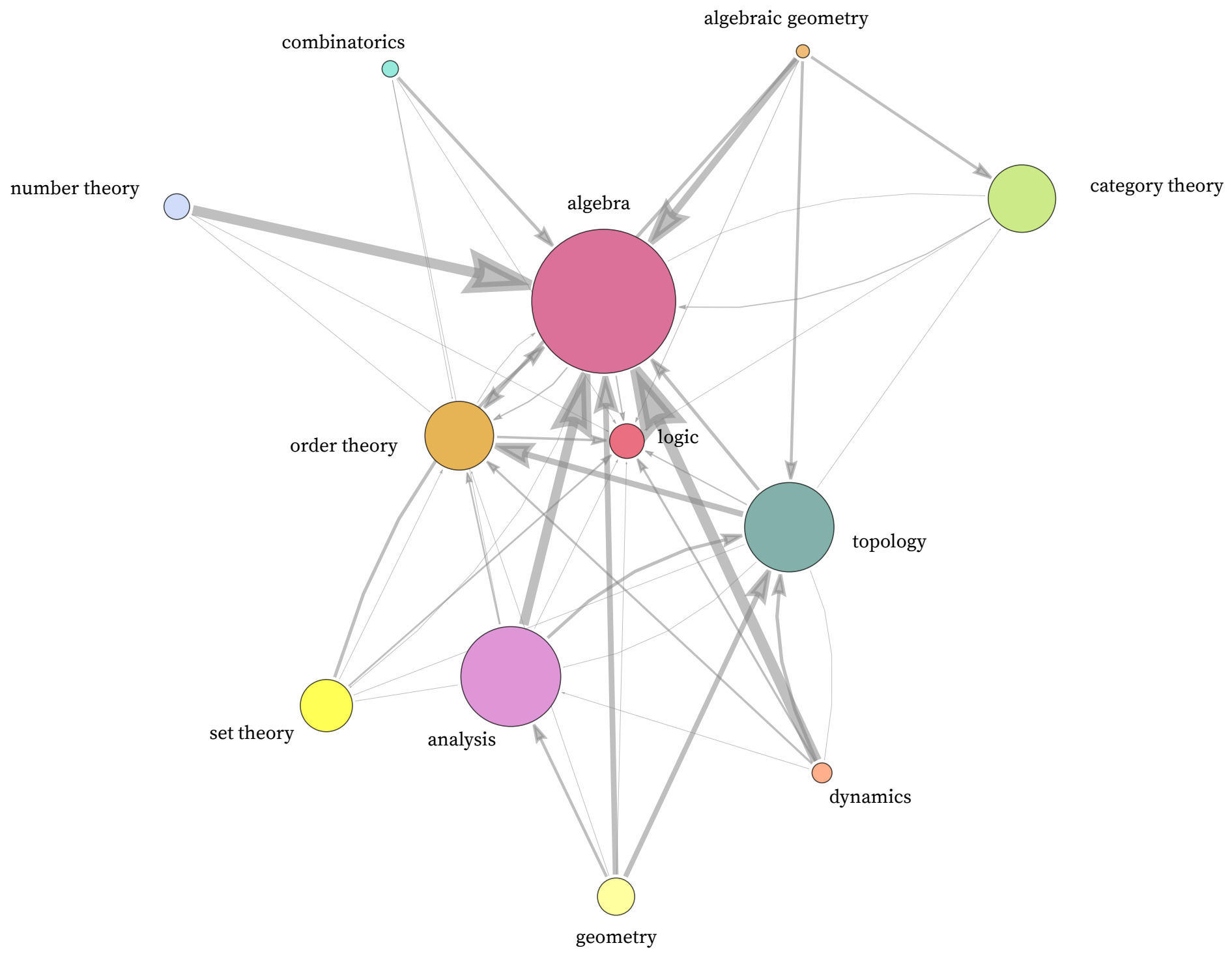

And here's the corresponding one for Metamath:

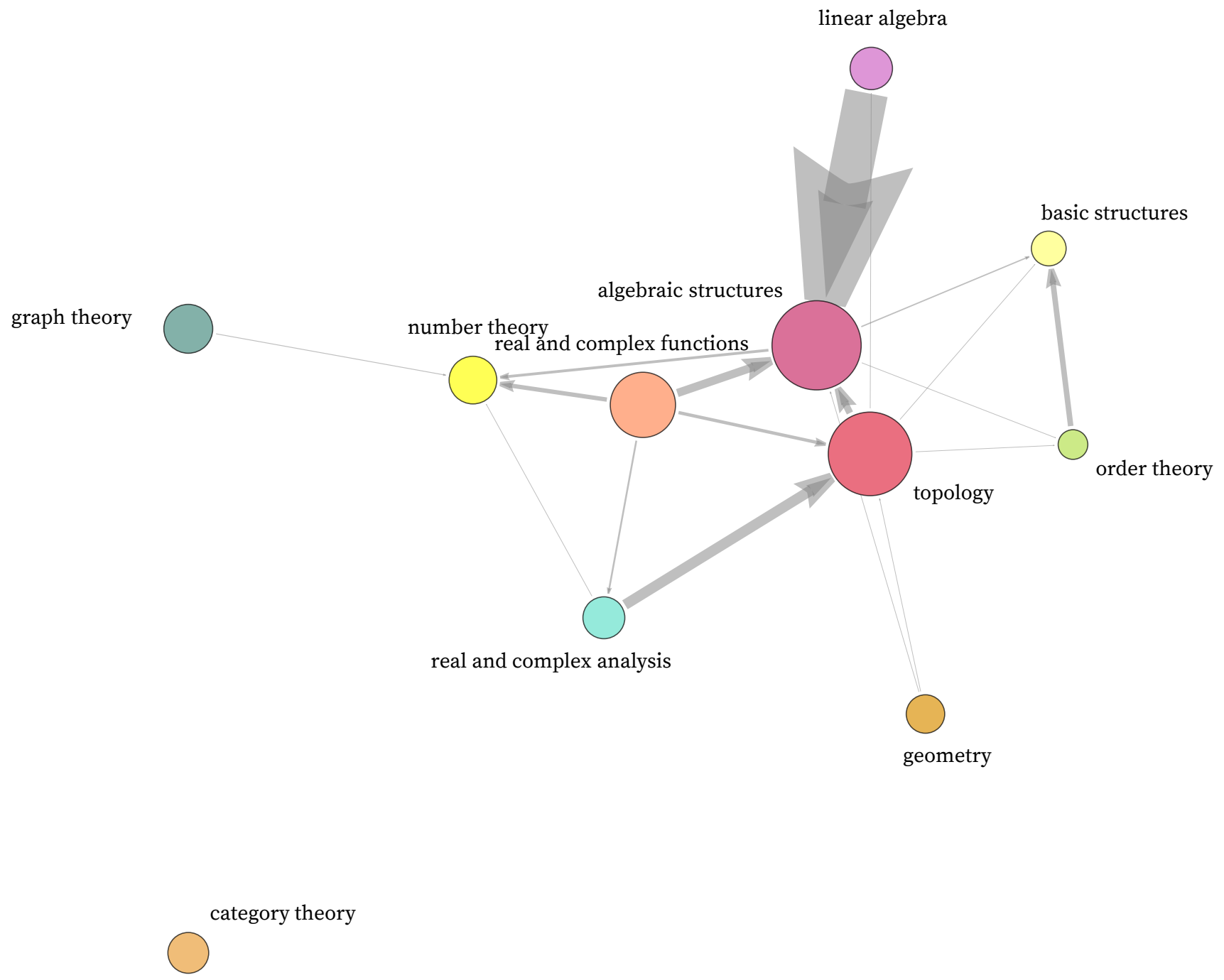

It's somewhat interesting to see how central algebra ends up being in both cases, and how comparatively "off on the side" category theory is. But it's clear that much of what one's seeing in these graphs is a reflection of the particular user communities of these systems, with some important pieces of modern mathematics (like the applications of algebraic geometry to number theory) notably missing.

But, OK, how do individual theorems work in these systems? As an example, let's consider the Pythagorean theorem. In Euclid, this is 1.47, and here's the first level of its dependency graph:

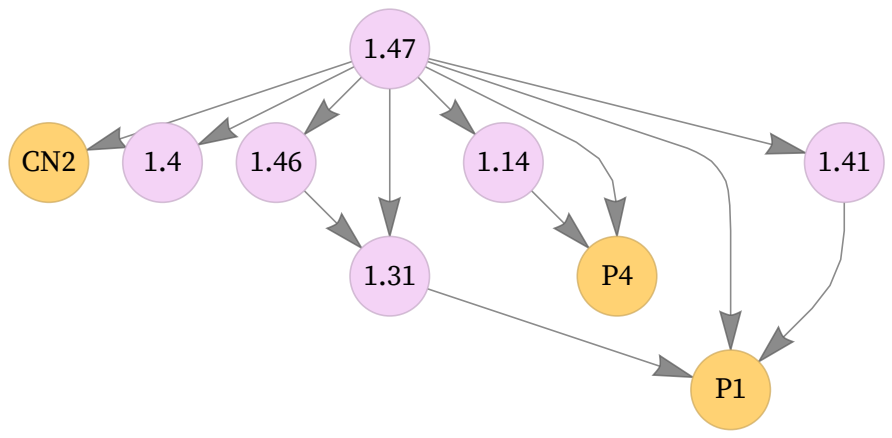

Here's the full graph involving a total of 39 elements (including, by the way, all 10 of the axioms), and having "depth" 20:

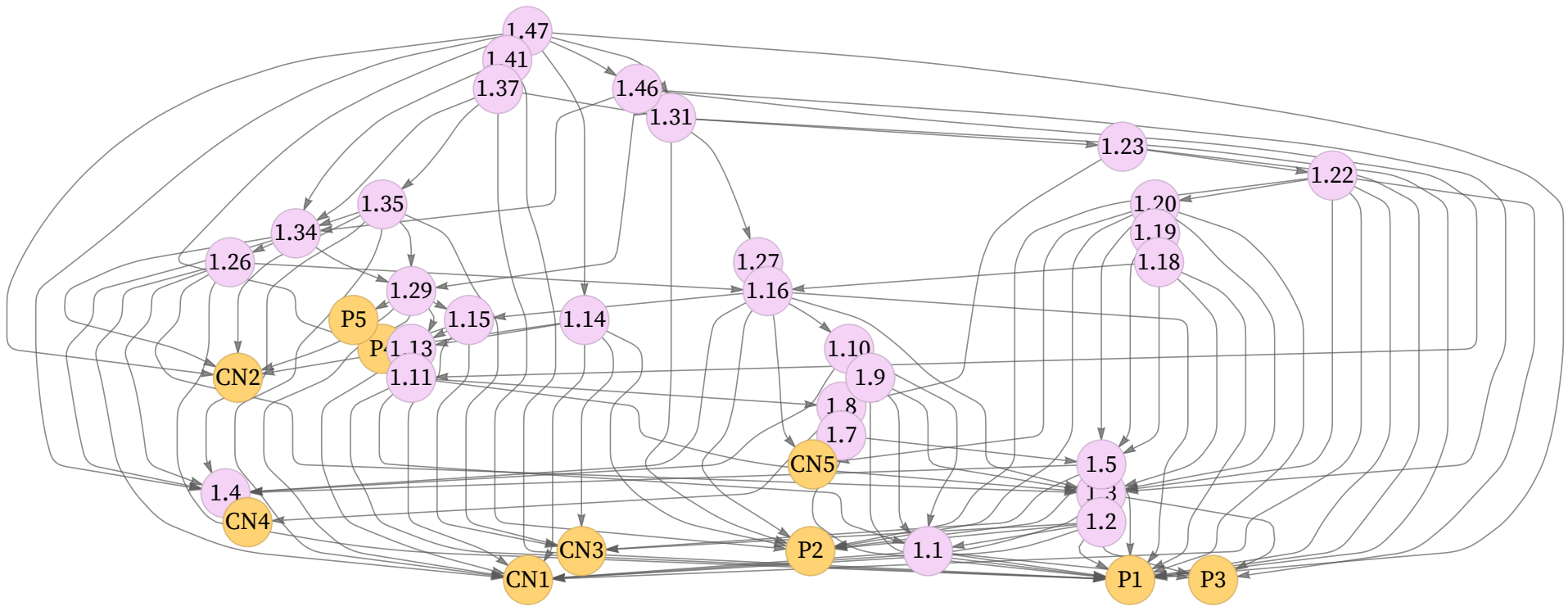

In Lean's mathlib, the theorem is called euclidean_geometry.dist_square_eq _dist _square _add _dist _square _iff _angle _eq _pi _div _two—and its stated proof directly involves 7 other theorems:

| |
|---|
| inner_product_geometry.norm_sub_square_eq_norm_square_add_norm_square_iff_angle_eq_pi_div_two |
| norm_neg |
| eq.symm |
| iff.refl |
| vsub_sub_vsub_cancel_right |
| dist_eq_norm_vsub |
| neg_vsub_eq_vsub_rev |

Going 3 steps, the theorem dependency graph looks like (where "init" and "tactic" basically refer to structure rather than mathematical content):

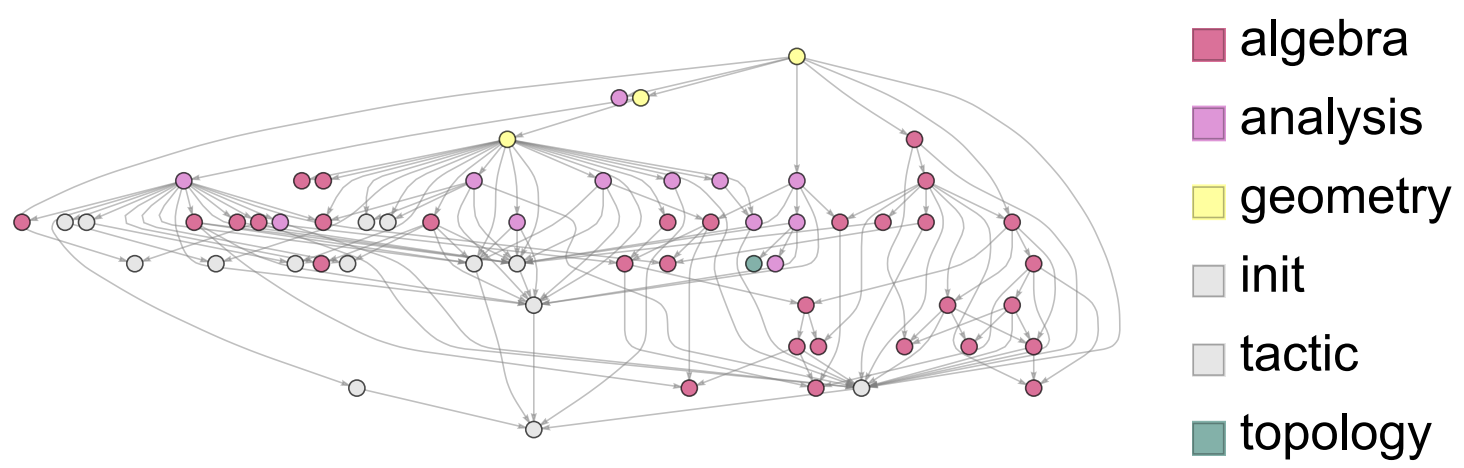

The full graph involves a total of 2850 elements (and has “depth” 84), and after transitive reduction has the form:

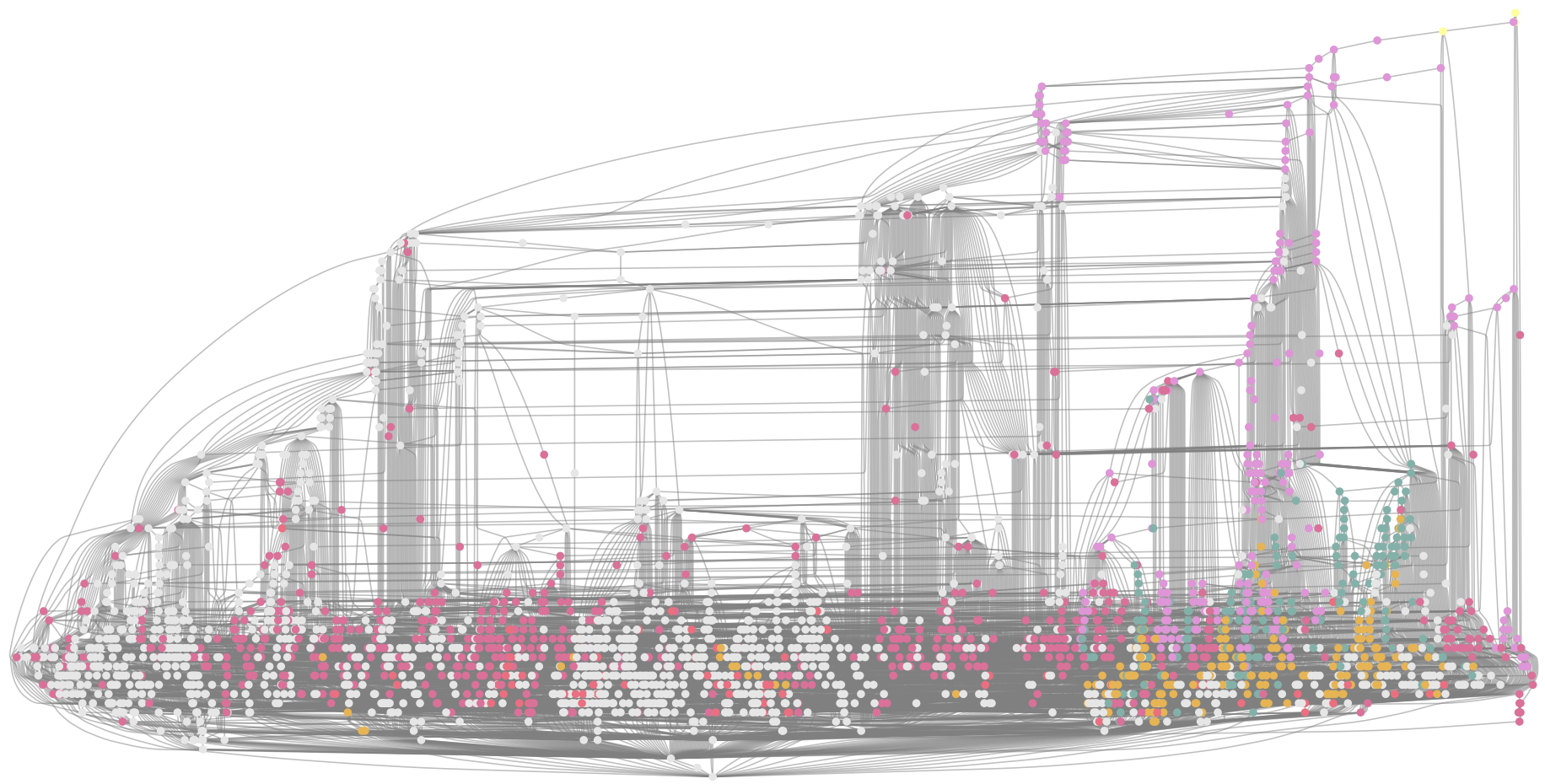

And, yes, this is considerably more complicated than Euclid’s version—but presumably that’s what happens if you insist on full formalization. Of the 2850 theorems used, 1503 are basically structural. The remainder bring mathematical content from different areas, and it’s notable in the picture above that different parts of the proof seem to “concentrate” on different areas. Curiously, theorems from geometry (which is basically all Euclid used) occupy only a tiny sliver of the pie chart of all theorems used:

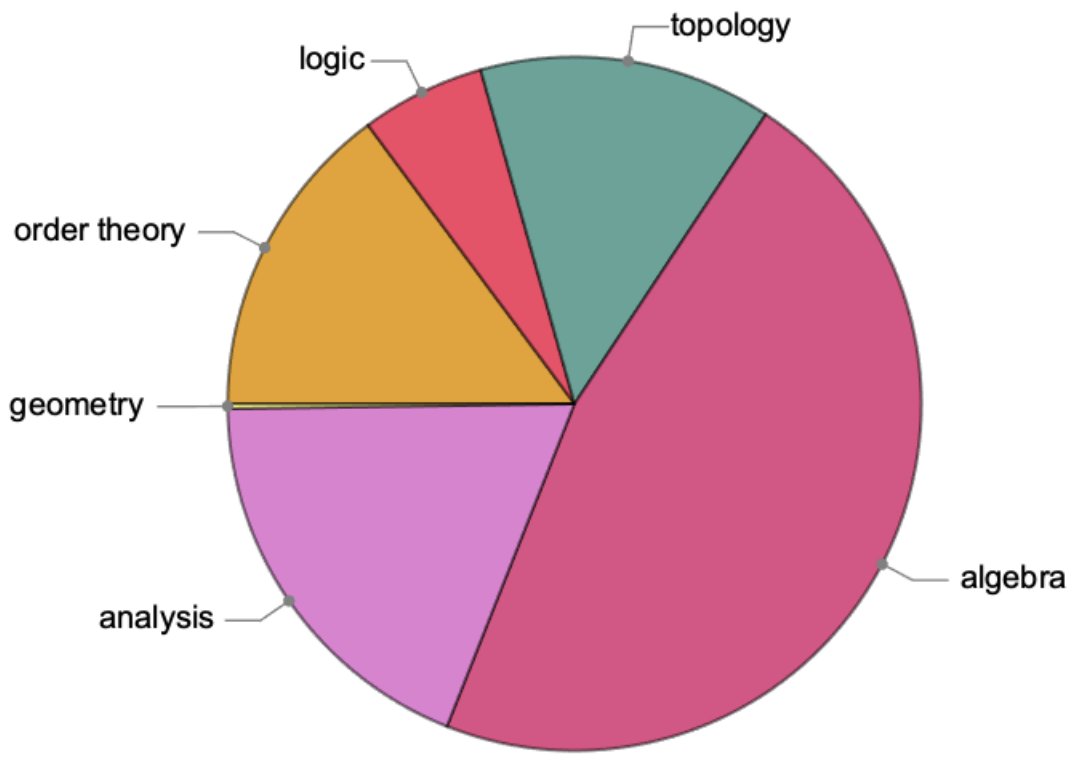

The Metamath set.mm version of the Pythagorean theorem is called `pythag`, and its proof directly depends on 26 other theorems:

| | | | | | |
|---|---|---|---|---|---|
| lawcos | 3adant3 | elpri | fveq2 | coshalfpi | syl6eq |
| cosneghalfpi | jaoi | syl | 3ad2ant3 | oveq2d | subcl |
| 3adant1 | 3ad2ant1 | abscld | recnd | syl5eqel | 3adant2 |
| mulcld | mul01d | eqtrd | 2t0e0 | sqcld | addcld |
| subid1d | 3eqtrd | | | | |

After 1 step, the theorem dependency graph is:

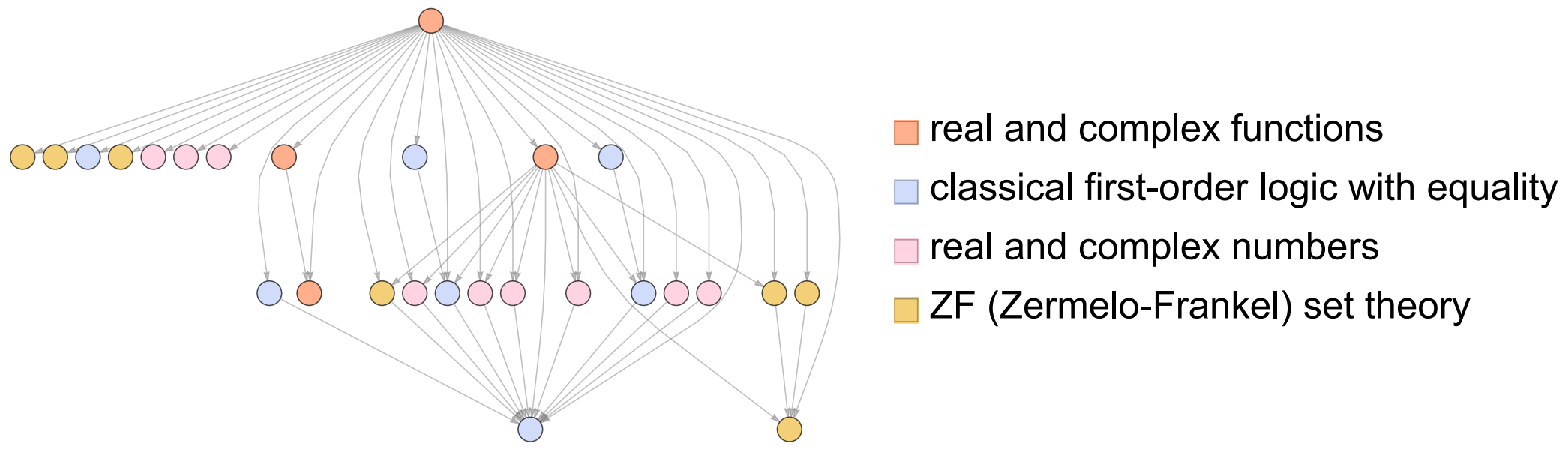

The full graph involves 7099 elements—and has depth 270. In other words, to get from the Pythagorean theorem all the way to the axioms can take as many as 270 steps.

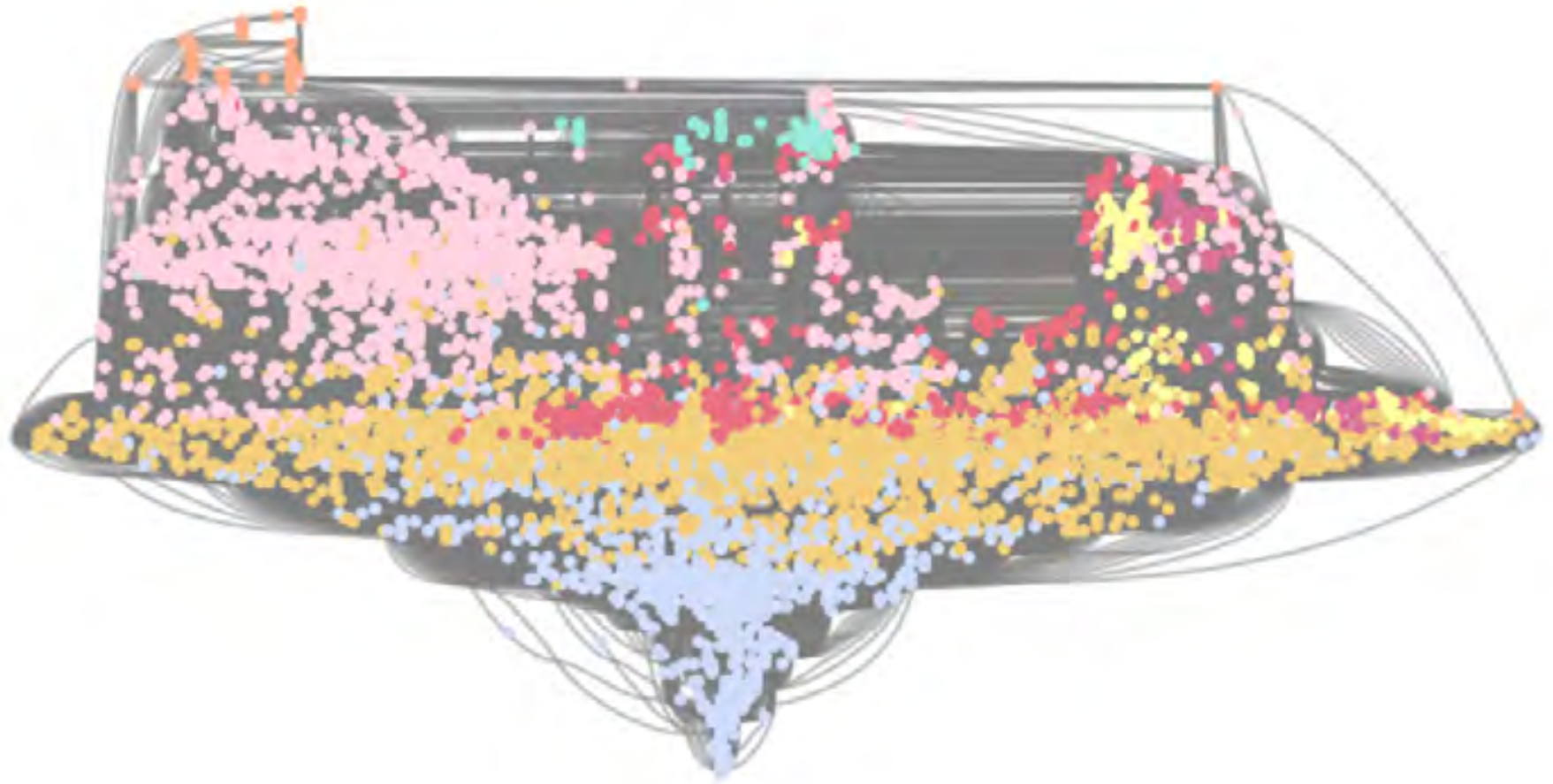

Given the complete Lean or Metamath corpuses, we can start doing the same kind of empirical metamathematics we did for Euclid's *Elements*—except now the higher level of formalization that's being used potentially allows us to go much further.

As a very simple example, here's the distribution of numbers of theorems directly referenced in the proof of each theorem in Lean, Metamath and Euclid:

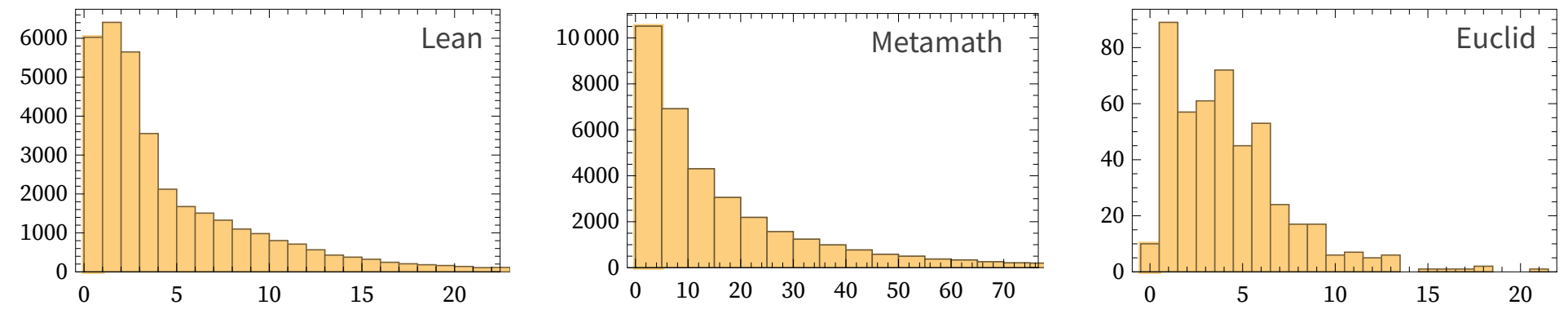

The differences presumably reflect different “hierarchical modularity conventions” in Lean and Metamath (and Euclid). But it’s interesting to note, for example, that in all three cases, the Pythagorean theorem is “above average” in terms of number of theorems referenced in its proof:

| | Lean | Metamath | Euclid |
|---|---|---|---|
| | 7 | 26 | 8 |
| mean | 4.9 | 18.7 | 4.3 |

What are the most popular theorems used in proofs? In terms of direct references, here are the top-5 lists:

| | |
|---|---|
| congr_arg | 9896 |
| congr | 9241 |
| eq.trans | 8632 |
| eq.symm | 6032 |
| eq_self_iff_true | 3711 |

Lean

| | |
|---|---|
| syl | 11524 |
| eqid | 8360 |
| adantr | 7434 |
| syl2anc | 6469 |
| a1i | 5727 |

Metamath

| | |
|---|---|
| 10.11 | 60 |
| 6.1 | 53 |
| 5.11 | 47 |
| 1.3 | 47 |
| 10.6 | 43 |

Euclid

Not surprisingly, for Lean and Metamath these are quite “structural”. For Lean, `congr_arg` is the “congruency” statement that if $a=b$ then $f(a)=f(b)$; `congr` is a variant that says if $a=b$ and $f=g$ then $f(a)=g(b)$; `eq.trans` is the transitivity statement if $a=b$ and $b=c$ then $a=c$ (Euclid’s CN1); `eq.symm` is the statement if $a=b$ then $b=a$; etc. For Metamath, `syl` is “transitive syllogism”: if $x \Rightarrow y$ and $y \Rightarrow z$ then $x \Rightarrow z$; `eqid` is about reflexity of equality; etc. In Euclid, these kinds of low-level results—if they are even stated at all—tend to be “many levels down” in the hierarchy of theorems, leaving the single most popular theorem, 10.11, to be one about proportion and rationality.

If one looks at all theorems directly and indirectly referenced by a given theorem, the distribution of total numbers of theorems is as follows (with Lean showing the most obviously exponential decay):

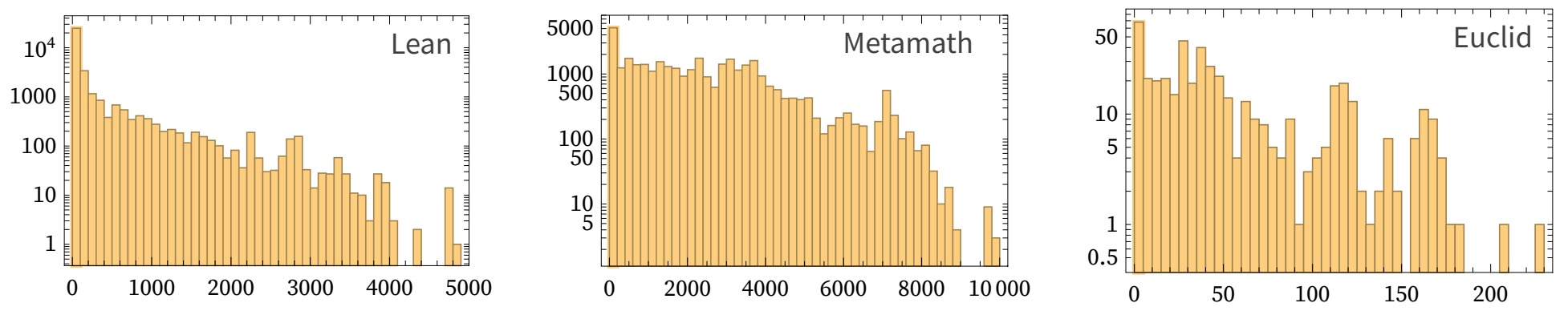

What about the overall structure of the Lean and Metamath dependency graphs? We can ask about effective dimension, about causal invariance, about "event horizons", and much more. But right now I'll leave that for another time...

## The Future of Empirical Metamathematics

I don't think empirical metamathematics has been much of a thing in the past. In fact, looking on the web as I write this, I'm surprised to see that essentially all references to the actual term "empirical metamathematics" seem to point directly or indirectly to that one note of mine on the subject in *A New Kind of Science*.

But as I hope this piece has made clear, there's a lot that can be done in empirical metamathematics. In everything I've written here, I haven't started analyzing questions like how one can recognize a powerful or a surprising theorem. And I've barely scratched the surface even of the empirical metamathematics that can be done on Euclid's *Elements* from 2000 years ago.

But what kind of a thing is empirical metamathematics? Assuming one's looking at theorems and proofs constructed by humans rather than by automated systems, it's about analyzing large-scale human output—a bit like doing data science on literary texts, or on things like websites or legal corpuses. But it's different. Because ultimately the theorems and proofs that are the subject of empirical metamathematics are derived not from features of the world, but from a formal system that defines some area of mathematics.

With computational language the goal is to be able to describe anything in formalized, computational terms. But in empirical metamathematics, things are in a sense "born formalized". Whatever the actual presentation of theorems and proofs may be there, their "true form" is ultimately something grounded in the formal structure of the mathematics being used.

Of course there is also a strong human element to the raw material of empirical metamathematics. It is (at least for now) humans who have chosen which of the infinite number of possible theorems should be considered interesting, and worthy of presentation. And at least traditionally, when humans write proofs, they usually do it less as a way to certify correctness, and more as a form of exposition: to explain to other humans why a particular theorem is true, and what structure it fits into.

In a sense, empirical metamathematics is a quite desiccated way to look at mathematics, in which all the elegant conceptual structure of its content has been removed. But if we're to make a "science of metamathematics", it's almost inevitable that we have to think this way. Part of what we need to do is to understand some of the human aesthetics of mathematics, and in effect to see to deduce laws by which it may operate.

In this piece I've mostly concentrated on doing fairly straightforward graph-oriented data science, primarily on Euclid's *Elements*. But in moving forward with empirical metamathematics a key question is what kind of model one should be trying to fit one's observations into.

And this comes back to my current motivation for studying empirical metamathematics: as a window onto a general “bulk” theory of metamathematics—and as the foundation for a science not just of how we humans have explored metamathematical space, but of what fundamentally is out there in metamathematical space, and what its overall structure may be.

No doubt there are already clues in what I’ve done here, but probably only after we have the general theory will we have the paradigm that’s needed to identify them. But even without this, there’s much to do in studying empirical metamathematics for its own sake—and of better characterizing the remarkable human achievement that is mathematics.

And for now, it’s interesting to be able to look at something as old as Euclid’s *Elements* and to realize what new perspectives modern computational thinking can give us about it. Euclid was a pioneer in the notion of building everything up from formal rules—and the seeds he sowed played an important role in leading us to the modern computational paradigm. So it’s something of a thrill to be able to come back two thousand years later and see that paradigm—now all grown up—applied not only to something like the fundamental theory of physics, but also to what Euclid did all those years ago.

## Thanks

*For help with various aspects of the content of this piece I’d like to thank Peter Barendse, Ian Ford, Jonathan Gorard, Rob Lewis, Jose Martin-Garcia, Norm Megill, James Mulnix, Nik Murzin, Mano Namuduri, Ed Pegg, Michael Trott, and Xiaofan Zhang, as well as Sushma Kini and Jessica Wong, and for past discussions about related topics, also Bruno Buchberger, Dana Scott and the various participants of our 2016 workshop on the Semantic Representation of Mathematical Knowledge.*

## Note Added

As I was working on this piece, I couldn’t help wondering whether—in 2300 years—anyone else had worked on the empirical metamathematics of Euclid before. Turns out (as Don Knuth pointed out to me) at least one other person did—more than 400 years ago.

The person in question was Thomas Harriot (1560–1621).

The only thing Thomas Harriot published in his lifetime was the book *A Briefe and True Report of the New Found Land of Virginia*, based on a trip that he made to America in 1585. But his papers show that he did all sorts of math and science (including inventing the · notation for multiplication, < and >, as well as drawing pictures of the Moon through a telescope before Galileo, etc.). He seems to have had a well-ahead-of-his-time interest in discrete mathematics, apparently making Venn diagrams a couple of centuries before Venn

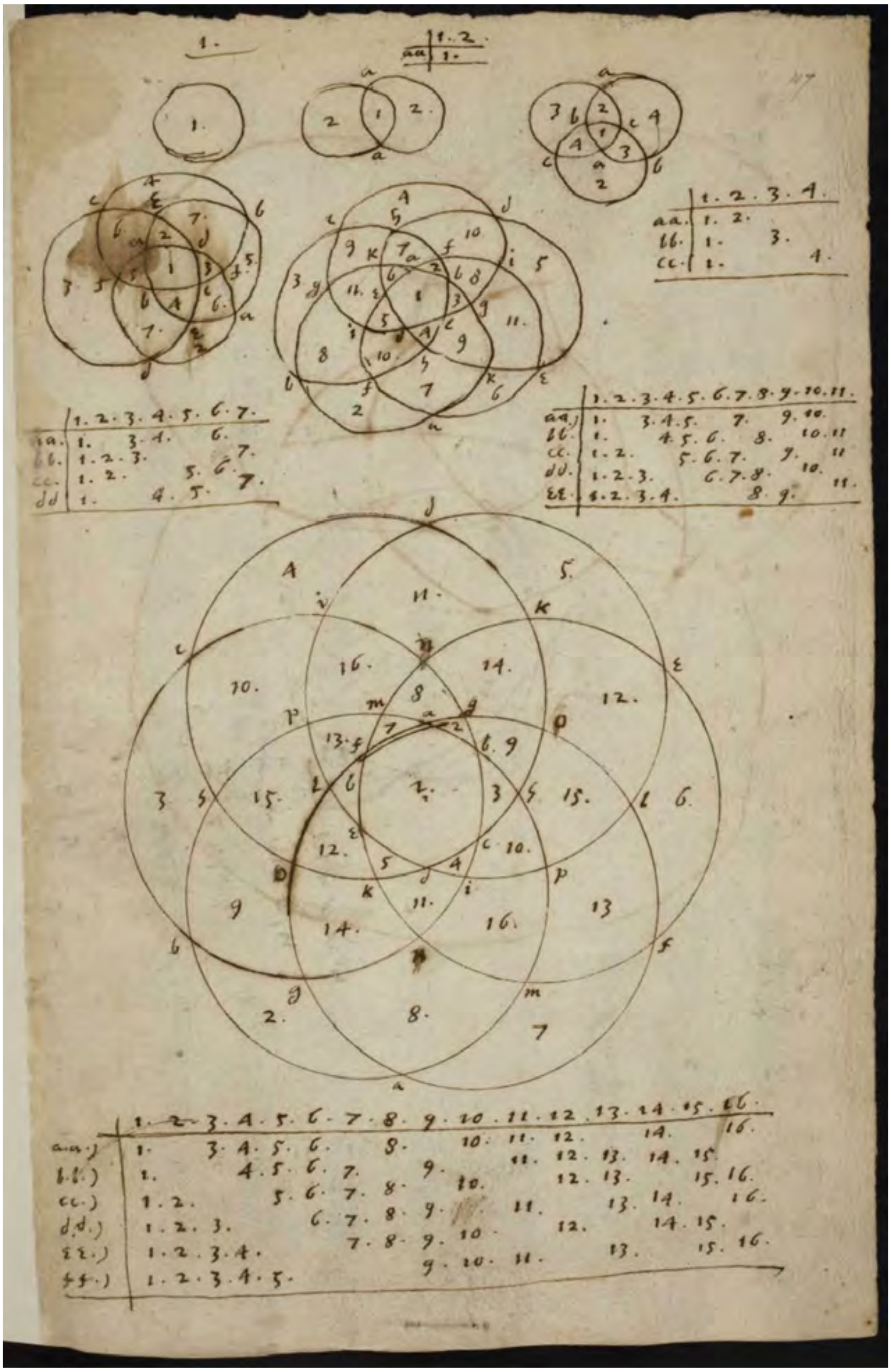

doing various enumerations of structures

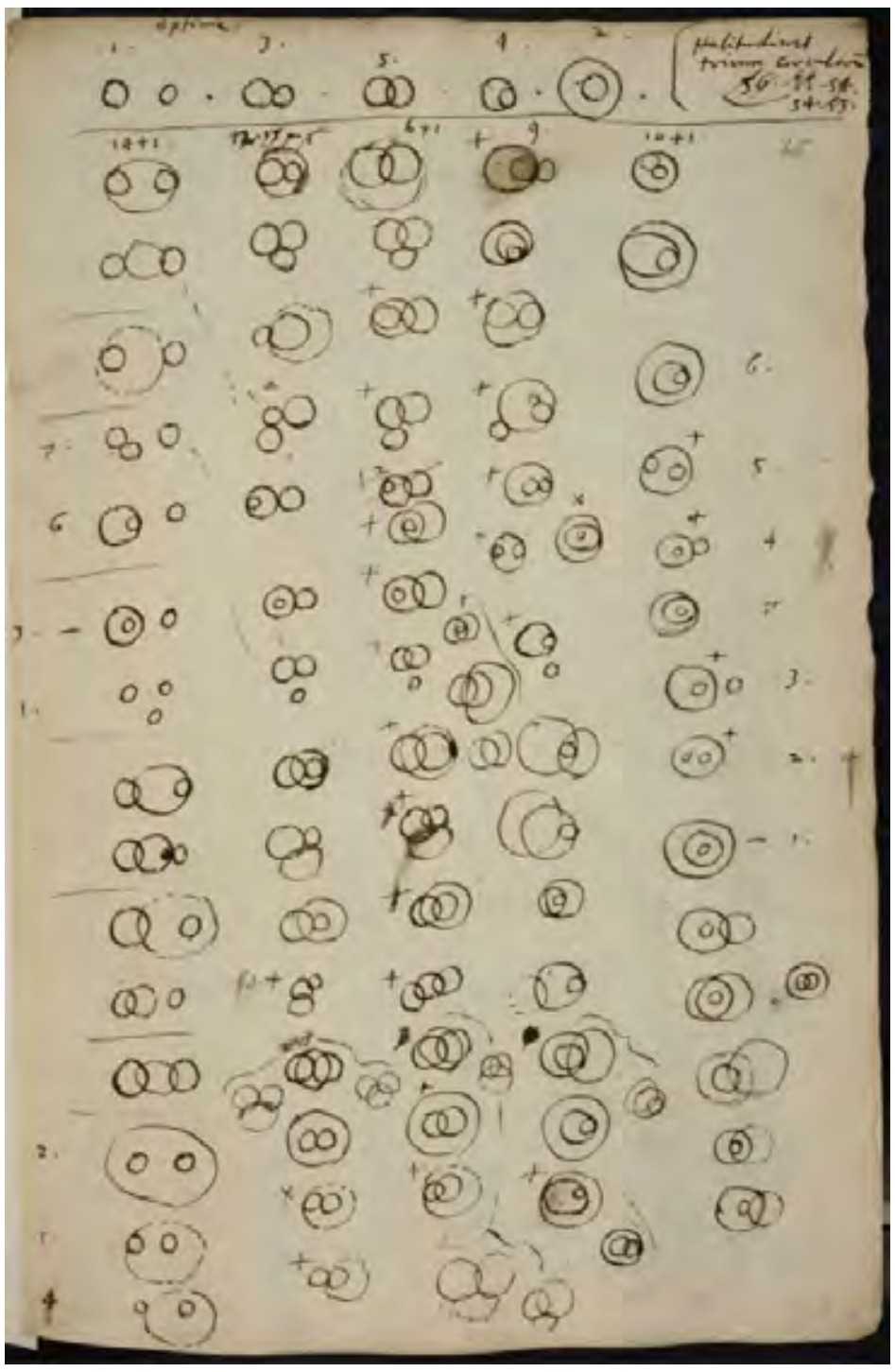

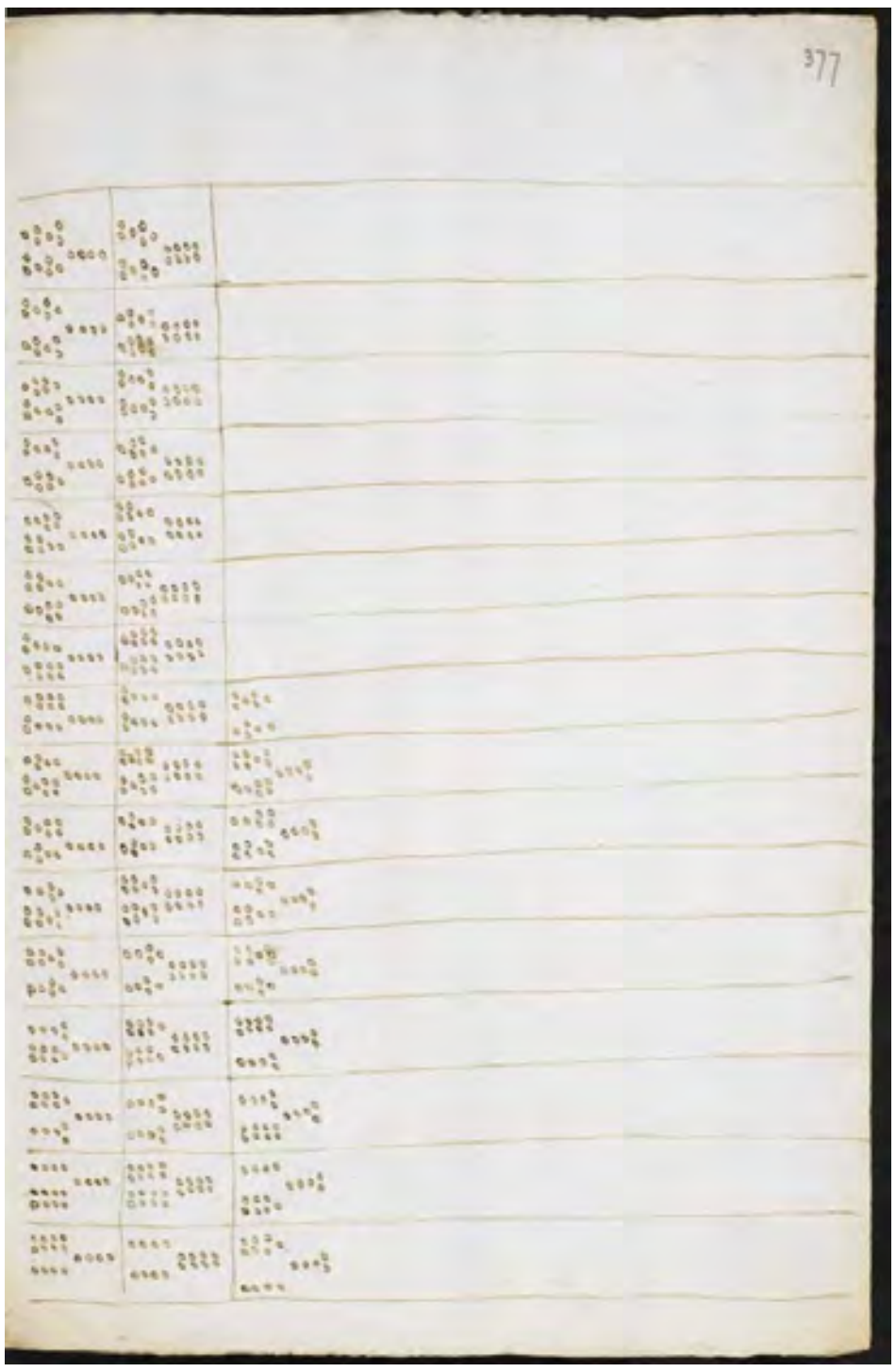

as well as various repeated computations (but no cellular automata, so far as I can tell!):

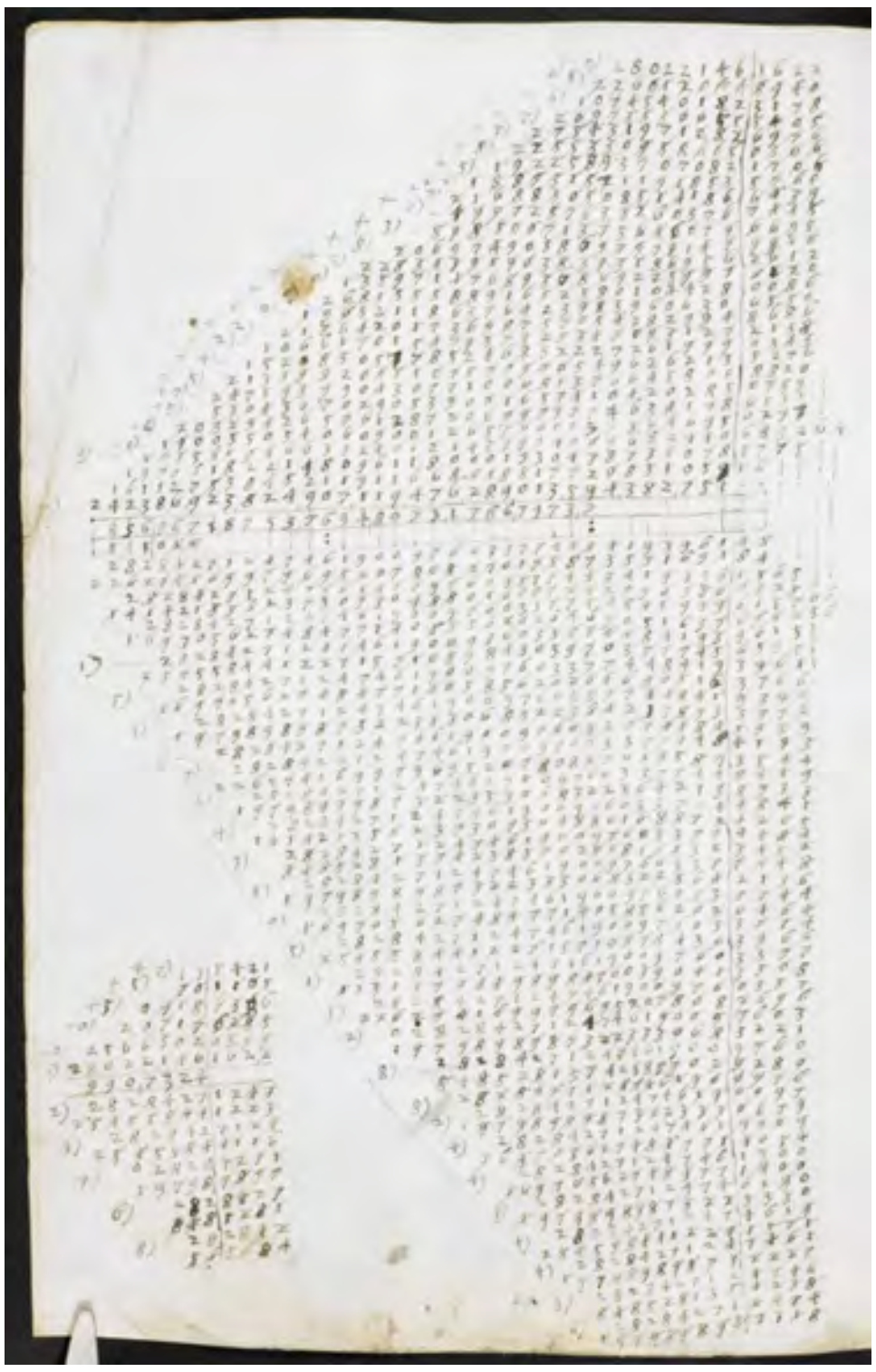

And he seems to have made a detailed study of Euclid's *Elements*, listing in detail (as I did) what theorems are used in each proof (this is for Book 1):

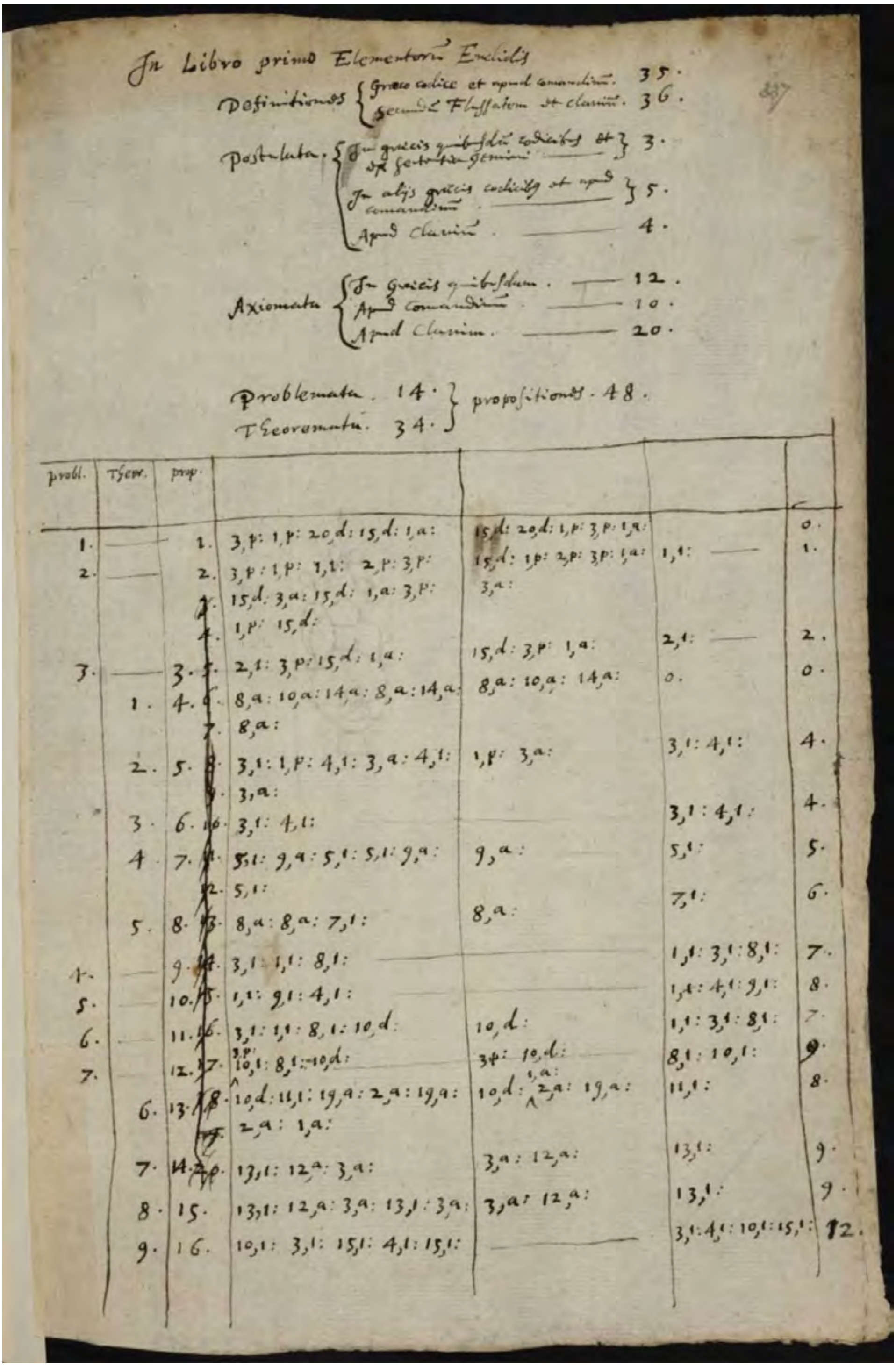

But then, in his "moment of empirical metamathematics" he lists out the full dependency table for theorems in Book 1, having computed what we'd now call the transitive closure:

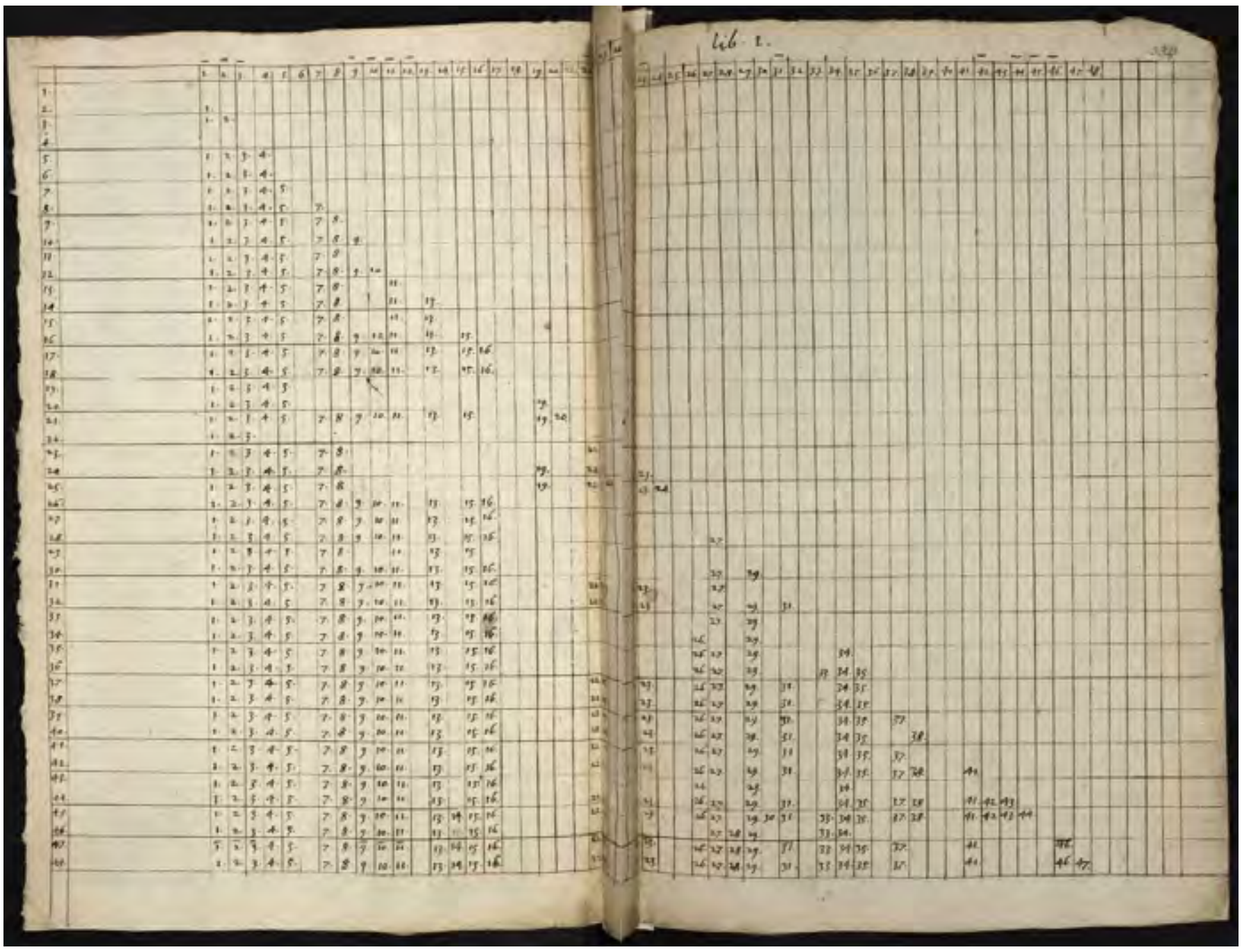

It's easy for us to reproduce this now, and, yes, he did make a few mistakes:

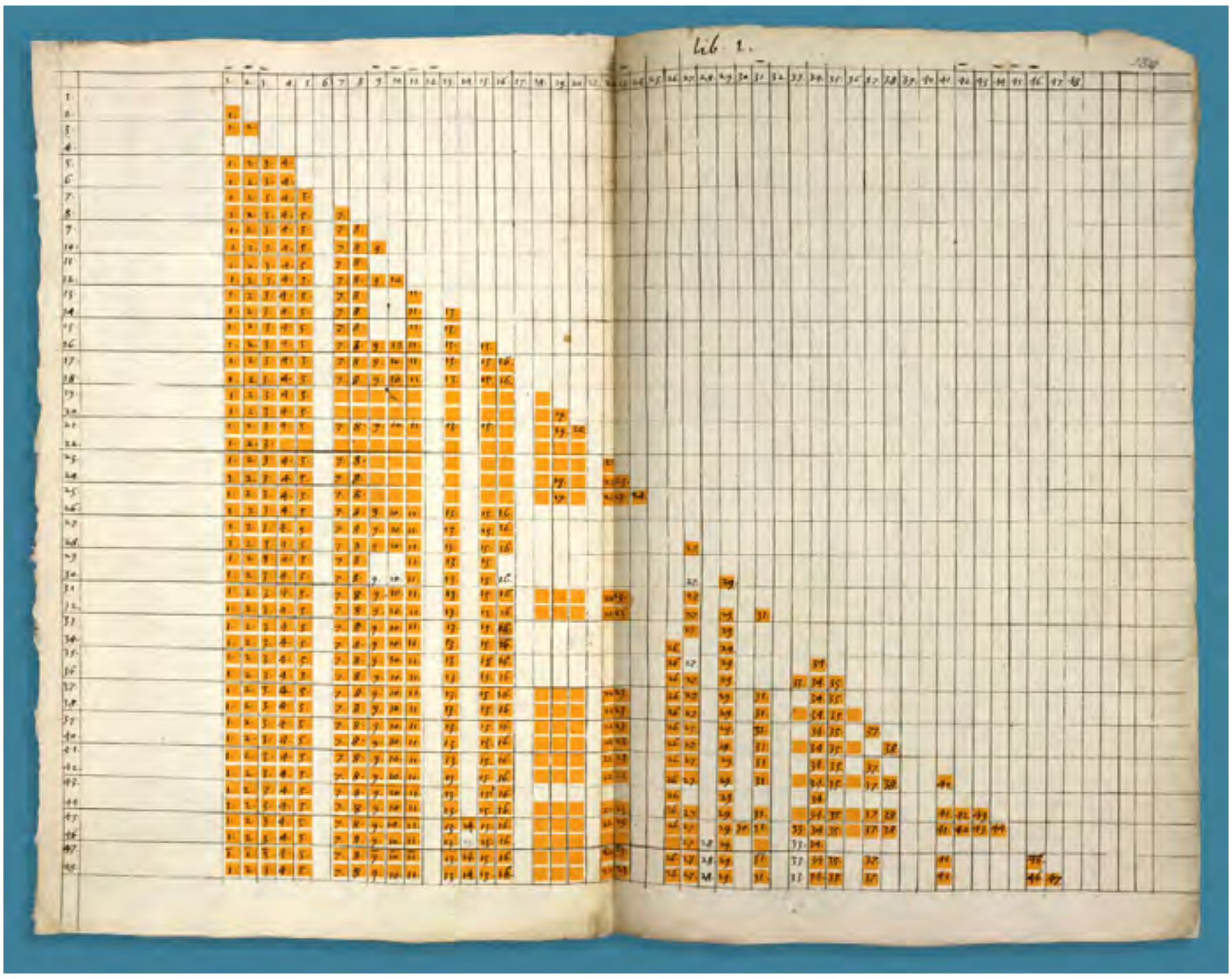

Studying the empirical metamathematics of Euclid seems (to me) like an obvious thing to do, and it's good to know I'm not the first one doing it. And actually I'm now wondering if someone actually already did it not "just" 400 years ago, but perhaps 2000 (or more) years ago...

## References

*Links to references are included within the body of this document.*